\documentclass[12pt,a4paper]{amsart}%
\usepackage{imakeidx}
\usepackage{amsmath}
\usepackage{fullpage}
\usepackage{comment}
\usepackage{MyCommA}
\usepackage[final]{pdfpages}
\usepackage{tikz-cd}

\newcommand{\EitanB}[1]{{{#1}}}
\newcommand{\Rami}[1]{{{#1}}}
\newcommand{\RamiA}[1]{{{#1}}}
\newcommand{\RamiB}[1]{{{#1}}}
\newcommand{\RamiC}[1]{{{#1}}}
\newcommand{\RamiD}[1]{{{#1}}}

\newcommand{\RamiF}[1]{{{#1}}}

\newcommand{\DimaA}[1]{{{#1}}}
\newcommand{\DimaB}[1]{{{#1}}}

\newcommand{\EitanD}[1]{{{#1}}}
\newcommand{\EitanE}[1]{{{#1}}}

\newcommand{\noleft}{\left.\kern-\nulldelimiterspace}

\begin{document}
	
	\author[Aizenbud]{Avraham Aizenbud}
	\address{Avraham Aizenbud,
		Faculty of Mathematical Sciences,
		Weizmann Institute of Science,
		76100
		Rehovot, Israel}
	\email{aizenr@gmail.com}
	\urladdr{https://www.wisdom.weizmann.ac.il/~aizenr/}
	
	\author[Gourevitch]{Dmitry Gourevitch}
	\address{Dmitry Gourevitch,
		Faculty of Mathematical Sciences,
		Weizmann Institute of Science,
		76100
		Rehovot, Israel}
	\email{dimagur@weizmann.ac.il}
	\urladdr{https://www.wisdom.weizmann.ac.il/~dimagur/}

\author[Kazhdan]{David Kazhdan}
	\address{David Kazhdan,
    Einstein Institute of Mathematics, Edmond J. Safra Campus, Givat Ram, The
Hebrew University of Jerusalem, Jerusalem, 91904, Israel}
	\email{david.kazhdan@mail.huji.ac.il}
	\urladdr{https://math.huji.ac.il/~kazhdan/}
	
\author[Sayag]{Eitan Sayag}
	\address{Eitan Sayag,
    Department of Mathematics, Ben Gurion University of the Negev, P.O.B. 653,
Be’er Sheva 84105, Israel}
	\email{eitan.sayag@gmail.com}
	\urladdr{www.math.bgu.ac.il/~sayage}

	\date{\today}
	

	%
	%
	%
	%
	%
	%
	%
	%


    

\title{Effective results in representation theory of reductive groups over local fields}
\maketitle
\begin{abstract} 
We prove effective versions of several classical results in representation theory of reductive groups over non-Archimedean local fields  of arbitrary characteristic. {More specifically, we deal with the support of matrix coefficients of cuspidal representations
and the recursions satisfied by matrix coefficients of  representations \RamiF{of finite length}.
We provide effective generation results for irreducible representations.
Finally, we study the
effective smoothness of characters on the set of regular semi-simple elements.}
\end{abstract}

\tableofcontents 
\section{Introduction}
In this paper we present a few effective versions of well-known results of Casselman and Harish-Chandra in complex representation theory of reductive groups over non-Archimedean local fields of arbitrary characteristic. 

\subsection{Classical non-effective results}
Let us briefly enumerate the classical results \RamiA{for which we provide an effective version in} this paper:
\begin{enumerate}
    \item A matrix coefficient of a cuspidal representation has compact support modulo the center, see \DimaA{e.g.} \cite[VI.2.1]{Ren}.
    \item An irreducible representation is generated by a single vector. This is a tautology, but its effective version is non-trivial.
    \item A matrix coefficient of an irreducible representation satisfies a recursion relation when restricted to the monoid generated by a hyperbolic element. This follows from Casselman's canonical pairing, \DimaA{see e.g.} \cite[VI.9.6]{Ren}.
    \item The character of an irreducible representation, when restricted to the regular semi-simple locus, is given by a smooth function. See \cite[\S4, Corollary]{HC_sub}.
\end{enumerate}
\subsection{Effective results} 
In order to formulate the effective results rigorously one has to introduce some technical notation. We will do it in \S\ref{ssec:intro.not} below. We will now give an informal formulation.

Let $\pi$ be an irreducible (complex) representation of a reductive group $G$ defined over a local field $F$. Then
\begin{enumerate}
    \item Let $v\in \pi$ and $w\in \tilde \pi$ be of depth $n$ (i.e. they are fixed by the Moy--Prasad group of depth $n$). Assume that $\pi$ is cuspidal.
    
    Then the logarithm of the diameter of the support (modulo the center)  of the matrix coefficient  attached to $(v,w)$ is proportional to $n$. Note that the proportionality coefficient might depend on $\pi$. See \Cref{thm:suppm}.
    \item There exists a vector $v\in \pi$ s.t.  \RamiA{for any integer $n>0$, any depth $n$ vector} $w\in \pi$  can be obtained from $v$ by applying an element of $\C[G]$ whose support \RamiA{lies} inside a ball (around $1$) of logarithmic  diameter proportional to $n$. See \Cref{thm:supph} and \S\ref{ssec:app.supph}.
    
    \item Let $v\in \pi$ and $w\in \tilde \pi$ be of depth $n$. Let $m$ be the corresponding matrix coefficient. Let $g\in G$ be a hyperbolic element\footnote{Formally speaking, in our formulation we require a  slightly different condition on $g$, which will be explained below.}. Then 
    \RamiA{there exist a  recurrence relation and a constant $c$ depending only on $\pi$ and $g$ (and not on $v,w$), such that}    
    the sequence $k\mapsto m(g^k)$
    satisfies \RamiA{this} recurrence relation starting with $k=c\cdot n$.
    See \Cref{thm:mcexp}.

    \item \Rami{The character of $\pi$ is constant in a ball around every regular semi-simple element $g$. The logarithm of the radius of this ball depends linearly} on the logarithm of the norm of $g$ and the valuation of the Weyl discriminant of $g$. See \Cref{thm:smooth}.
\end{enumerate}

\subsection{Notation}\label{ssec:intro.not}
Throughout the paper $\mdef{\bfG}$ is a connected reductive group defined over a local field $\mdef{F}$. Let $\mdef{G}:=\bfG(F)$.
For a group  (or an algebraic group) $H$ denote by  $\mdef{\cZ(H)}$ its center.    For  a smooth (complex) representation $\pi$ of $G$ and $v\in \pi$,  $w\in \mdef{\tilde\pi}$ -- the contragredient representation -- we denote by $\mdef{m_{v,w}}:G\to \C$ the matrix coefficient. We also denote by $\mdef{\chi_\pi}$ the Harish-Chandra character of $\pi$ which is a generalized function.

In order to formulate our results we have to choose some additional data related to $\bfG$. In \S\ref{ssec:not} we specify how we make these choices.\footnote{\RamiA{Though the dependence of our results on these choices is rather weak, the choices cannot be made completely arbitrarily.}} For now we will just introduce the notation \Rami{needed} to formulate our results.

We choose:
\begin{itemize}
\item a uniformizer $\mdef{\varpi}\in F $
\item a minimal parabolic $\mdef{\bfB}\subset \bfG$
\item a Levi subgroup $\mdef{\bfT}<\bfB$. \DimaA{Denote by $\mdef{\bfA_\bfT}<\bfT$ its  maximal split torus.} 
\item  a maximal compact subgroup $\mdef{G_0}<G$
\item  a basis of the topology of $G$ at $1$ consisting of open compact subgroups $\mdef{G_{-n}}$, with $G_{-n-1}<G_{-n}<G_0$ for $n\in \Z_{\geq 0}$
    \item open compact subsets $\mdef{G_n^{ad}}\subset G^{ad}_{n+1}\subset \mdef{G^{ad}}:=G/\cZ(G)$ \Rami{for $n\in \Z_{> 0}$, that} cover $G^{ad}$
    \item open compact subsets $\mdef{G_n}\subset G_{n+1}\subset G$ for $n\in \Z_{> 0}$,  \Rami{that} cover $G$.
\end{itemize}

\subsection{Main results}

\begin{notation}
For an open  compact subgroup $K<G$ denote by
$\mdef{\cH_{K}(G)}$ the {Hecke algebra} consisting of compactly supported measures on $G$ that are bi-$K$-invariant.
\end{notation}
\begin{introtheorem}[\S\ref{ssec:app.supph}]\label{thm:suppm}
        Let $\pi$ be a cuspidal representation of $G$ of finite length.        
        Then there is a constant $\mdef{c_{supp}(\pi)}\in \Z_{>0}$ s.t. for any 
        \begin{itemize}
            \item $n\in \Z_{>0}$ 
            \item $v\in \pi^{G_{-n}}$,  $w\in \tilde\pi^{G_{-n}}$,
        \end{itemize}
we have $\Supp(m_{v,w})\subset \cZ(G)G_{ c_{supp}(\pi) \cdot n}$.
\end{introtheorem}
\begin{introtheorem}[\S \ref{sec:gen}]\label{thm:supph}
Let $\pi$ be an irreducible  representation. 
Then there exist  $\mdef{c_{gen}(\pi)}\in \Z_{>0}$ and $\mdef{v_{gen}(\pi)}\in \pi$ s.t. for any:
 \begin{itemize}
     \item  $n\in  \Z_{>0}$
     \item $w\in \pi^{G_{-n}}$ 
 \end{itemize}
there exists
$h\in \cH_{G_{-n}}(G)$ s.t.
    $\supp(h)\subset G_{c_{gen}(\pi)\cdot n}$ and
    $hv_{gen}(\pi)=w$.
\end{introtheorem}
For the next results we will need some more notation.

\begin{notation}\label{not:Lam.A}
   
Denote
    \begin{itemize}
        \item $\mdef{\Lambda({A_{\bfT}})}$ -- the lattice of co-characters of $\mdef{A_{\bfT}}:=\bfA_{\bfT}(F)$.  

                \item Denote by $\mdef{\cE}:=\mdef{\cE_{A_{\bfT}}}:\Lambda({A_{\bfT}})\to {A_{\bfT}}$ the natural map given by the choice of $\varpi$. We also denote $\mdef{\varpi^\lambda}:=\cE_{A_{\bfT}}(\lambda)$.
        \item $\mdef{{A_{\bfT}}^+}$ -- the dominant part of
        $A_{\bfT}$
        w.r.t. $\bfB$. See \cite[\S V.3.19]{Ren}.      
        \item $\mdef{\Lambda^+(A_{\bfT})}:=\cE^{-1}({A_{\bfT}}^+)$.
    \end{itemize}
\end{notation}

   \begin{introtheorem}[\S \ref{sec:asm}]\label{thm:mcexp}
        Let $\pi$ be a \Rami{(smooth)} representation of $G$ of finite length and $\lambda\in \Lambda^+(A_\bfT)$.
                Then there is  
                a constant 
                $\mdef{c_{asm}(\pi,\lambda)}\in \Z_{>0}$ and a monic 
                polynomial\footnote{\RamiA{In fact, this polynomial depends only on the Jacquet module of $\pi$ with respect to the standard parabolic corresponding to $\lambda$ (see \Cref{not:cP} below). The degree of this polynomial is \RamiF{at most} the length of this Jacquet module. Thus this degree is bounded by the size of the Weyl group of $G$ (see \cite[Cor. 2.13]{BZ1}).}}
                $$\mdef{\cP_{\pi,\lambda}}=\sum_{i=0}^{\deg\cP_{\pi,\lambda}} a_i x^i\in \C[x]$$ 
                s.t. for any 
        \begin{itemize}
            \item $n\in \Z_{>0}$, 
            \item $v\in \pi^{G_{-n}}$,  $w\in \tilde\pi^{G_{-n}}$,
            \item $k>c_{asm}(\pi,\lambda)n$,
        \end{itemize}
we have
$$\sum_{i=0}^{\deg\cP_{\pi,\lambda}}  a_im_{v,w}(\varpi^{(k+i)\lambda})=0.$$        
\end{introtheorem}

For the last result we introduce:

\begin{notation}\label{not:rss}
Denote 
\begin{itemize}
    \item $\mdef{\fD}$ -- the Weyl discriminant
    \DimaA{considered as  a  regular function on $\bfG$}
    \item $\mdef{\bfG^{rss}}:= \fD^{-1}(\bA^1\smallsetminus 0)$
    \item $\mdef{G^{rss}}:=\bfG^{rss}(F)$
    \item $\mdef{ov_{G^{rss}}}:G^{rss}\to \R$ -- the function defined by:
    $$ov_{G^{rss}}(x)=\max\left(\min(\{n|x\in G_n\}),\operatorname{val}(\fD(x))\right).$$
\end{itemize}
    
\end{notation}

\begin{introtheorem}[\S \ref{ThmSm}]\label{thm:smooth}
Let $\pi$ be an irreducible representation of $G$. Then there is $\mdef{c_{reg}(\pi)}\in \Z_{>0}$ s.t. for any
        $g\in G^{rss}$,
\Rami{the character $\chi_\pi$ is constant on}
$$
{g G_{-c_{reg}(\pi) \cdot (ov_{G^{rss}}(g)+1)}}
.$$    
\end{introtheorem}

\subsection{Related results}
Bounds on the support of matrix coefficients of cuspidal representations were studied in \cite{Ber_tam}. There it was proven that there is a bound on the support of such matrix coefficients of a given depth which is uniform not only in the specific matrix coefficients but also in the cuspidal representation.

However, this bound on the size of the support is exponential in the depth (in the logarithmic scale that we use in this paper) and not linear as in \Cref{thm:suppm}. On the other hand, the bound in \Cref{thm:suppm} is not uniform in the \DimaA{representation.}

In \cite[Theorem 3]{FL_bnd} a uniform linear bound is proven for a certain class of cuspidal representations, namely representations that are induced from a compact subgroup (modulo the center). Moreover their result is uniform within this class. {A further work in this direction is \cite{Lap}}.

For many reductive groups the representations covered by \cite[Theorem 3]{FL_bnd} are known to be all the cuspidal ones. But in general it is not known.

\Rami{In characteristic zero, the main result of \cite{Kor} provides effective
smoothness of characters near regular semi-simple elements under some mild assumption. These results are more explicit than  \Cref{thm:smooth}. However, they rely crucially on the characteristic zero assumption, which is absent in \Cref{thm:smooth}.
}

\subsection{Ideas of the proofs}
We start with proving the cuspidal case of \Cref{thm:supph}. Here the starting point is \Cref{prop:van.vec}\Rami{. It} states that given a vector $v$ in a cuspidal representation, elements in the Hecke algebra of a given depth $n$, which are supported far away from the center, annihilate $v$. Moreover, the required distance from the center depends linearly (in the logarithmic scale we are working in) on the depth $n$. Note that there is no effective estimate on the dependence on $v$. \RamiA{The proof of \Cref{prop:van.vec} follows the standard proof of its non-effective counterpart.}

\Cref{prop:van.vec} implies \Rami{the cuspidal case of} \Cref{thm:supph} rather easily.

Now we can deduce \Cref{thm:suppm}. 
Note that the classical approach is inherently non-effective. The effective approach of \cite{Ber_tam} gives an exponential bound and we could not improve it. So we use a different approach: We start with a single matrix coefficient $m_0$, which is known to be compactly supported (mod center). Here we do not need to bound its support since it is a single matrix coefficient.
Then we use \Cref{thm:supph} in order to show that any matrix coefficient $m$ can be obtained from a combination of shifts of $m_0$. Moreover, the magnitude of the shifts depends linearly on the depth of $m$. This implies the assertion.\footnote{In order to avoid duplications, we deduce \Cref{thm:suppm} from the general case of \Cref{thm:supph} and not from the cuspidal case. This is purely technical, and one can deduce \Cref{thm:suppm} from the cuspidal case of \Cref{thm:supph} in the same way.}

We deduce the general case of \Cref{thm:supph} using an explicit construction of a generator for a parabolically induced representation, in terms of a generator for the representation of the Levi subgroup.

 The general idea of the proof of \Cref{thm:mcexp} is again to use the non-effective version of \Cref{thm:mcexp} for a single matrix coefficient $m_0$, then to use the fact that any matrix coefficient  $m$ can be obtained from (a combination of) shifts of $m_0$. However this time, equipped with \Cref{thm:supph} we can bound the magnitude of the shifts in terms of the depth of the matrix coefficient $m$. Finally we can try to use it to give the desired effective result for $m$.

However we could not make this idea to work as is. The reason is that we could not obtain an effective bound on the recursion satisfied by  the shift of $m_0$ in terms of the recursion satisfied by $m_0$ and the magnitude of the shifts. 

So we modified this strategy. \Rami{We first prove} 
the following result:

   \begin{introtheorem}[\Cref{prop:mcexp}]\label{thm:eff.jak.cass}
        Let $\pi$ be a representation of $G$ of finite length and $\lambda\in \Lambda^+(A_\bfT)$. Let $U_\lambda$ be the unipotent radical of \Rami{the} standard parabolic group corresponding to $\lambda$.
        
                Then there is  
                a constant 
                $\mdef{c_{asm0}(\pi,\lambda)}\in \Z_{>0}$ and a monic polynomial $$\mdef{\cP_{\pi,\lambda}}=\sum_{i=0}^{\deg\cP_{\pi,\lambda}} a_i x^i\in \C[x]$$ 
                s.t. for any 
        \begin{itemize}
            \item $n\in \Z_{>0}$ 
            \item $v\in \pi^{G_{-n}}$
            \item $k>c_{asm0}(\pi,\lambda)n$,
        \end{itemize}
we have
$$\sum_{i=0}^{\deg\cP_{\pi,\lambda}}  a_i \varpi^{i\lambda} \cdot v= \sum_{i=1}^{N}(u_iw_i-w_i)$$
for some $u_i\in U_\lambda \cap G_k$ and $w_i\in \pi$.  
\end{introtheorem}
It turns out that such a statement behaves well under shifts of $v$ so the strategy above works for it. Finally we deduce \Cref{thm:mcexp} from \Cref{thm:eff.jak.cass}.

For \Cref{thm:smooth} we carefully analyze Harish-Chandra's original proof and see that in view of \cite{AGKS0}, it actually gives \Cref{thm:smooth}.
\subsection{Motivation}
For us, the main motivation for this  paper is our project of extending Harish-Chandra's Integrability theorem to reductive groups over local fields of positive characteristic (\cite{AGKS2,AAGKS,AGKS_jet,AGKS_orb,AGKS_eff}). 

Specifically, we plan to use this paper to generalize the results of \cite{AGKS2,AAGKS} to arbitrary irreducible representations.

In \cite{AGKS2} we use Harish-Chandra's formula for characters of cuspidal representations in terms of matrix coefficients. This formula is also available for discrete series representations, but even in that case it is not effective enough for our purposes.

The results of the current paper \Rami{are} used (together with \cite{KYtempered}) 
in a sequel paper \cite{AGKS_KY} to obtain  an explicit and effective version of Harish-Chandra's  formula for \Rami{characters of} tempered representations.

We plan to use \cite{AGKS_KY} in order to generalize the results of 
\cite{AGKS2,AAGKS} to tempered representations, see \cite{AGKS_int}.  This should imply the results \Rami{of 
\cite{AGKS2,AAGKS}} for arbitrary representations. 

\subsection{Structure of the paper}
In \S\ref{sec:prelim} we fix our conventions, we make the choices of the data listed in \S\ref{ssec:intro.not}, we introduce the notations attached to them, and we record the basic facts on these data that are used throughout the paper.

In \S\ref{sec:cusp} we prove the cuspidal case of \Cref{thm:supph}, namely \Cref{thm:supph.cusp}. The main step is \Cref{prop:van.vec}.

In \S\ref{sec:gen} we prove \Cref{thm:supph}. We first bound the effective generation of a parabolic induction $i_M^G(\rho)$ in terms of the effective generation of $\rho$, and then deduce the general case of \Cref{thm:supph} from its cuspidal case. In \S\ref{ssec:app.supph} we derive consequences of \Cref{thm:supph}: effective generation by finitely many vectors for representations of finite length (\Cref{cor:gen.fin.len.alt}), and \Cref{thm:suppm}.

In \S\ref{sec:asm} we prove \Cref{thm:mcexp}. We introduce the notion of an $n$-admissible vector, we show that admissibility is preserved, with a controlled loss, under the action of $G$, we deduce \Cref{thm:eff.jak.cass} from \Cref{cor:gen.fin.len.alt}, and finally we deduce \Cref{thm:mcexp} from \Cref{thm:eff.jak.cass}.

In \S\ref{sec:Norms} we recall the theory of norms on algebraic varieties over local fields developed in \cite[\S18]{Kot}. In particular, we prove in \Cref{prop:G.norms} that the size functions $ov_G$ and $ov_{G^{rss}}$, used in the formulation of our results, are the logarithms of norms.

In \S\ref{sec:rect} we recall the theory of rectified algebraic varieties, and of the balls and the measures attached to them, introduced in \cite{AGKS0}, and we adapt it to our purposes.

In \S\ref{ThmSm} we prove \Cref{thm:smooth}, by analyzing Harish-Chandra's proof of the smoothness of characters on $G^{rss}$ and making it effective using \S\ref{sec:Norms} and \S\ref{sec:rect}.

In \Cref{sec:app} we prove the supporting statements on the data fixed in \S\ref{ssec:not}. In \S\ref{ssec:app.Theta} we study the function $\Theta$ on the dominant cone, in \S\ref{ssec:app.norms} we compare the two size functions on $G$ with the Cartan parameter, and in \S\ref{ssec:app.ZG} we bound an element of $G$ modulo the center in terms of $\Theta$. In \S\ref{ssec:app.root} we record the behavior of the filtration of the root subgroups, and in \S\ref{ssec:app.ball} we relate the balls of the rectified variety $\bfG^{rss}$ to the cosets of the groups $G_{-n}$.

\subsection{Acknowledgments}
\RamiA{We are grateful to  \EitanE{Yakov} Varshavsky and Jessica Fintzen for answering our questions regarding the Bruhat-Tits building.}

During the preparation of this paper, A.A., D.G. and E.S. were partially supported by the ISF grant no. 1781/23. 
D.K. was partially supported by an ERC grant 101142781.

\section{Preliminaries}\label{sec:prelim}
\subsection{Conventions}
Here we fix general notation and conventions that we will use throughout the paper.

\begin{enumerate}
    \item By a \tdef{variety} we mean a reduced scheme of finite type over $F$. Similarly for an algebraic group.

    \item We will usually denote algebraic varieties by bold face letters (such as $\bfX$) and the spaces of their $F$-points by the corresponding usual face letters (such as $X:=\bfX(F)$). We use the same conventions when we want to interpret vector spaces as algebraic varieties.

    \item We will use the same letter to denote a morphism between algebraic varieties and the corresponding map between the sets of their $F$-points.
\item By a representation of an $l$-group we will mean a complex valued smooth representation.
    \item {When no ambiguity is possible  we will denote the adjoint action simply by $``\cdot"$.}    
    \item For an $l$-space $X$ we will denote by $\mdef{\Sc(X)}:=C_c^\infty(X)$ the space of locally constant compactly supported functions on $X$ and by $\mdef{\Sc^*(X)}$ its dual space.
    For a closed subset $Z\subset X$ we denote  by $\mdef{\Sc_Z^*(X)}\subset \Sc^*(X)$  the subspace of distributions supported in $Z$.
    \item \RamiA{For a group (or more generally a monoid) $G$ we will denote by $\C[G]$ its \href{https://grokipedia.com/page/Group_ring}{group algebra}, i.e. the algebra of \EitanE{complex valued} functions \EitanE{on $G$} with finite support (where the product is the convolution product).}
    \item {We will use the symbol \tdef{$<$} to denote the (not necessarily proper) containment relation for groups}.
\item We will use the words parabolic/Levi to describe subgroups both of algebraic reductive groups and of the groups of their points.

\item 

Several statements in this paper will concern the existence of objects that satisfy some conditions. Sometimes we will want to fix such an object. In these cases 
 we assign a distinct name for the corresponding object inside the formulation of the statement. It is implied that after each such statement we fix such an object and we will refer to it later by this name. Nothing significant will depend on the choices of these objects.

\end{enumerate}

\subsection{Notation and choices of the data}\label{ssec:not}
Let us fix the choices described in \S\ref{ssec:intro.not} and introduce some notation related to it.
\begin{notation}\label{not:settings}$ $
    \begin{enumerate}
        \item\label{it:choices} We make the choices of $\varpi\in F$, $\bfB\subset \bfG$, $\bfT<\bfB$, $\bfA_\bfT<\bfT$ from \S \ref{ssec:intro.not} in an arbitrary way.
        \item \RamiA{Let $\mdef{T}=\bfT(F), \mdef{B}=\bfB(F), \mdef{A_\bfT}:=\bfA_\bfT(F)$.}
        \item Denote by $\mdef{\cB\cT(G)}$ the (extended) {Bruhat--Tits building} of $G$.
        \item\label{it:vertex} Fix a special vertex $\mdef{\fx}\in  \cB\cT(G)$ that lies in the apartment that corresponds to $\bfA_{\bfT}$. Such a vertex exists by \cite[Proposition 1.3.43]{KP}.
        \item We set $G_0:=G_\fx$
        \item \RamiC{For $n\in\Z_{>0}$,} we set $G_{-n}:=\mdef{G_{\fx,n}}$ to be the {Moy--Prasad group}. (See \cite{MP96} or \cite[Definition 13.2.1]{KP}.)
                \item\label{it:0H} For a connected reductive group $\bfH$  with $H=\bfH(F)$ denote by $\mdef[^0G,\, ^0H,\, ^0M,\, ^0T,\, ^0A]{^0H}$ the 
        intersection of kernels of all unramified characters of $H$, see e.g. \cite[V.2.3]{Ren}. Every compact subgroup of $H$ is contained in $^0H$ (\cite[V.2.3]{Ren}). \RamiC{I}n particular, if $H$ has no proper parabolics, so that $^0H$ is compact, then $^0H$ is the unique maximal compact subgroup of $H$.
        \item For a connected reductive group  $\bfH$ denote $\mdef{\Lambda(\bfH)}:=\mdef{\Lambda(H)}:=H/^0H$ where $H=\bfH(F)$. Note that it is compatible with \Cref{not:Lam.A} by \cite[V.2.6]{Ren}. For non-connected $\bfH$ we denote $\Lambda(\bfH)=\Lambda(\bfH^\circ)$, where $\bfH^\circ$ is the identity component of $\bfH$.
        \item\label{it:Phi} Denote by $\mdef{\Phi}$ the set of non-divisible roots of the root system defined by $(\bfG,\bfA_\bfT)$, i.e. those roots $\beta$ of $\bfA_\bfT$ in $\Lie(\bfG)$ for which $\beta/2$ is not a root. It is itself a root system. We consider its elements as $\Q$-valued functionals on $\Lambda(T)$. See \cite[V.3.7, V.3.17, V.2.6]{Ren}.
        \item A \tdef{standard parabolic subgroup} of $\bfG$ is a parabolic that contains $\bfB$. 
        A  \tdef{standard Levi subgroup} of $\bfG$ is a Levi subgroup  of a standard parabolic subgroup of $\bfG$ that contains $\bfT$. We apply analogous terminology for parabolic and Levi subgroups of $G$.
\item For a standard parabolic $\bfP<\bfG$ we denote by $\mdef{\bfP^-}$ the unique parabolic s.t. $\bfP\cap \bfP^{-}$ is a standard Levi. We call it the \tdef{opposite parabolic} of $\bfP$.
    \item \DimaA{    
    For an open compact subgroup $K<G$ and a standard parabolic $P<G$ with a standard Levi $M$, we denote by $\mdef{K_M}:=pr_M(K\cap P)$ and    $\mdef{K^{-}}:=K\cap U^{-}$, where $pr_M:P\to M$ is the projection and $U^{-}$ is the unipotent radical of $P^-$.}
     \item For a standard Levi $M<G$ denote by $\mdef{i_{M}^G}$ the \RamiA{(unnormalized)} {parabolic induction} w.r.t. the standard parabolic that includes $M$.   
        \item For a connected reductive group $\bfH$ denote by $\mdef{\bfA_\bfH}$ the maximal split subtorus of $\cZ(\bfH)$.

\item\label{it:Phi.P} For $\bfP<\bfG$ a standard parabolic or the opposite of one, with unipotent radical $\bfU_\bfP$, denote by $\mdef{\Phi_\bfP}\subset\Phi$ the set of $\beta\in\Phi$ occurring in $\mdef{\bfU_\bfP}$, i.e. those whose weight space for $\bfA_\bfT$ lies in $\Lie(\bfU_\bfP)$. Then $\Phi=\Phi_\bfB\sqcup\Phi_{\bfB^-}$ and $\Phi_{\bfP^-}=-\Phi_\bfP$. If $P=\bfP(F)$ we use the same notation with $\bfP$ replaced by $P$ to denote the same set.
        \item\label{it:Delta} Denote by $\mdef{\Delta}\subset\Phi_\bfB$  the collection of {simple roots} of the root system defined by $(\bfB,\bfA_{\bfT})$. See \cite[V.3.10]{Ren}.
        \item\label{it:M+} Let $$\mdef{T^+}=\{t\in T|\forall \beta\in \Phi_\bfB \text{ we have } \beta(t)\geq 0  \}$$ See \cite[\S V.3.19]{Ren}.
        \item Similarly define $\mdef{T^{++}}$ to be the strictly dominant part \cite[\S V.3.19]{Ren}.
        \item\label{it:Lam+} Denote $\mdef{\Lambda^+(T)}:= T^+/^0T$.
        \item Similarly define $\mdef{\Lambda^{++}(T)}:= T^{++}/^0T$.
        \item Pulling back to $A_{\bfT}$ we get \tdef{$\Lambda^+(A_{\bfT})$} and $\mdef{\Lambda^{++}(A_{\bfT})}$. Note that this is compatible with \Cref{not:Lam.A} above.
\item\label{it:N.Theta}        Choose $\mdef{\fN}\in\Z_{>0}$ such that
$\fN\beta(\lambda)\in\Z$ for every $\lambda\in\Lambda(T)$ and every $\beta\in\Phi$, the elements
of $\Phi$ being viewed as $\Q$-valued functionals on $\Lambda(T)$ as in \eqref{it:Phi}. Such an
$\fN$ exists since $\Phi$ is finite and $\Lambda(T)$ is finitely generated.
Denote
$$
\mdef{\Theta}=\fN\sum_{\beta\in\Phi_\bfB}\beta .
$$
        \item Let $\mdef{G_0^{ad}}$ be the image of $G_0$ in $G^{ad}$. 
        \item\label{it:cartan} $T^+\to G$ induces a bijection 
        $\Lambda^+(T)\to G_0 \backslash G/ G_0 $ (\tdef{Cartan decomposition}, \cite[V.5.1]{Ren}).
        \item\label{it:ZG} Note that $\cZ(G)\subset T^+$. This induces an embedding of $\Lambda(\cZ(G))$ into $\Lambda^+(T)$.
        \item\label{it:ov1} The map in \eqref{it:cartan} induces a bijection $\Lambda^+(T)/\Lambda(\cZ(G))\to G_0^{ad} \backslash G^{ad}/ G_0^{ad}$. 
         Note also that $\Theta$ is trivial on $\Lambda(\cZ(G))$. So it gives a function on $\Lambda^+(T)/\Lambda(\cZ(G))$ and hence on $G_0^{ad} \backslash G^{ad}/ G_0^{ad}$ and on $G^{ad}$. We denote this function by $\mdef{ov_{G^{ad}}}$.
        \item\label{it:ov2} Let $\mdef{\bfG^{ab}}:=\bfG/[\bfG,\bfG]$ and let  $\mdef{G^{ab}}:=\bfG^{ab}(F)$. Note that $\Lambda(G^{ab})$ is a lattice. Choose a basis for it. Let $\mdef{ov_{\Lambda(G^{ab})}}:\Lambda(G^{ab})\to \Z_{\geq 0}$ be  the $L^\infty$ norm w.r.t. this basis. Lift it to the function $\mdef{ov_{G^{ab}}}:G^{ab}\to \Z_{\geq 0}$.
        \item\label{it:ov3} Define $\mdef{ov_G}:G\to \Z_{\geq 0}$ by $\max(ov_{G^{ab}}(g^{ab}),ov_{G^{ad}}(g^{ad})),$ where $g^{ab}$ and $g^{ad}$ are the images of $g \in G$ in $G^{ab}$ and $G^{ad}$. 
        \item\label{it:Gn} 
        \RamiC{For $n\in\Z_{\geq 0}$,}
        let $\mdef{G^{ad}_n}:=ov_{G^{ad}}^{-1}([0,n])$ and
        $\mdef{G_n}:=ov_{G}^{-1}([0,n])$. 
        \item For a standard Levi $M<G$,                we will use the notations $\mdef{ov_M}, \mdef{ov_{M^{ad}}}, \mdef{M_n}, \mdef{M_{-n}}, \mdef{M^{ad}_n}$. We define them analogously to the case $M=G$. Note that they depend on some choices as above. We make these choices in the following way:
\begin{itemize}
    \item In \eqref{it:choices} we choose the same $\bfT$ and ${\bfA_\bfT}$. 
    \item  In \eqref{it:choices} the minimal parabolic that we choose is $\bfB\cap \bfM$.
    \item There is a natural map 
    $\mdef{\beta_M} $ between the apartment $\mdef{\cB\cT(G)_{\bfA_\bfT}}$ that corresponds to $A_\bfT$ inside $\cB\cT(G)$  and the apartment $\cB\cT(M)_{\bfA_\bfT}$; it is the restriction of the map $\pi$ of
    \cite[Proposition 9.7.2]{KP}.\footnote{\RamiA{Note that \cite[Proposition 9.7.2]{KP} deals with the reduced building, but it induces a map for the extended building in a natural way.}}  In \eqref{it:vertex} we choose the point $\beta_M(\fx)\in \cB\cT(M)_{\bfA_\bfT}$.
    \item  In \eqref{it:N.Theta} we choose the same integer $\fN$. Note that $\fN$ also  satisfies the condition of \eqref{it:N.Theta} for all standard Levi subgroups: every root of $(\bfM,\bfA_\bfT)$ is a root of $(\bfG,\bfA_\bfT)$, hence lies in $\Phi$ or is twice an element of $\Phi$.
\end{itemize}
     \end{enumerate}
\end{notation}
\subsection{Basic facts on the data}
\begin{lemma}\label{lem:Gn.product}
For every
$r,s\in\Z_{\RamiA{\geq}0}$, we have
$$
G_rG_s\subset G_{r+s}.
$$
\end{lemma}

\begin{proof}
This follows from the definition of $G_r$ by the function $ov_G$.
For the abelian part it is the subadditivity of the fixed norm on
$\Lambda(G^{ab})$. For the adjoint part it is the usual subadditivity of the Cartan projection \EitanE{(see
\Cref{not:settings}
(\ref{it:cartan}))}\RamiA{:
$$G\to G_0 \backslash G/G_0\cong \Lambda^+(T).$$}
\end{proof}

\begin{notation}\label{not:sigma}
Fix a closed embedding $\mdef{\sigma}:\bfG\to \GL(\mdef{\cW})$ of algebraic groups over $F$ and an
$\cO_F$-lattice $\mdef{\cL}\subset \cW$ with $\sigma(G_0)\cL=\cL$. Such a lattice exists: for any
$\cO_F$-lattice $\cL_0\subset \cW$ the matrix entries of $\sigma(k)$, $k\in G_0$, are bounded, so
$\cL:=\sum_{k\in G_0}\sigma(k)\cL_0$ is an $\cO_F$-submodule with
$\cL_0\subset \cL\subset \varpi^{-N}\cL_0$ for some $N$, hence a lattice, and it is $\sigma(G_0)$-stable.

Fix an $\cO_F$-basis of $\cL$, write $\sigma(g)_{ij}$ for the matrix entries of $\sigma(g)$ in this
basis, and define $\mdef{ov_\sigma}:G\to\Z_{\geq0}$ by
$$
ov_\sigma(g):=\max_{i,j}\max\bigl(0,-\operatorname{val}(\sigma(g)_{ij}),
-\operatorname{val}((\sigma(g)^{-1})_{ij})\bigr).
$$
\end{notation}

\begin{lemma}[\Cref{ssec:app.norms}]\label{cor:ov.sigma.G}
There is $c\in\Z_{>0}$ such that
$$c^{-1}\bigl(ov_\sigma(g)-c\bigr)\leq ov_G(g)\leq c\bigl(ov_\sigma(g)+c\bigr)$$
for every $g\in G$.
\end{lemma}

\begin{lemma}[\Cref{ssec:app.ZG}]\label{lem:Theta.ZG.bound}
There is $\mdef{c_Z}\in\Z_{>0}$ with the following property. Let
$$
g=k_1ak_2,\qquad k_1,k_2\in G_0,\quad a\in T^+ .
$$
If $\Theta(a)\leq R$, then
$$
g\in \cZ(G)G_{c_Z(R+1)} .
$$
\end{lemma}

\begin{lemma}[\Cref{ssec:app.root}]\label{lem:U.contract}
There is $\mdef{c_{cnj}}\in\Z_{>0}$ s.t. for every
\begin{itemize}
\item $Q<G$ a standard parabolic or the opposite of one, with Levi subgroup $M_Q$ and
unipotent radical $U_Q$,
\item $k\in\Z_{>0}$,
\item $a\in T$ with $\beta(a)\geq c_{cnj}k$ for every $\beta\in\Phi_Q$,
\end{itemize}
we have
\begin{enumerate}
\item\label{lem:U.contract:1} $a\,(U_Q\cap G_k)\,a^{-1}\subset U_Q\cap G_{-k}$.
\item\label{lem:U.contract:2} If $a\in A_{\bfM_Q}$ then for every $n\in\Z_{>0}$ and every
$g\in G_k\cap M_Q$,
$$
a\,(U_Q\cap G_{-n})\,a^{-1}<g^{-1}(U_Q\cap G_{-n})\,g .
$$
\end{enumerate}
\end{lemma}

\begin{lemma}[\RamiC{\cite[Proposition 13.2.5(5)]{KP}}
]\label{lem:iva.decomp}
    Let $P<G$ be a standard parabolic. Let $M<P$ be its standard Levi subgroup and $U<P$ be its unipotent radical. Let $U^-$ be the unipotent radical of $P^-$. Then
    for $n>0$ we have 
    $$G_{-n}=(U^-\cap G_{-n})(M\cap G_{-n})(U\cap G_{-n}).$$
\end{lemma}

\begin{lemma}[\Cref{ssec:cGM}]\label{lem:cGM}
\RamiA{Let $M<G$ be a standard Levi.}
Then 
there exists $\mdef{c^G_M}\in \Z_{>0}$ s.t.  for any $n\in \Z_{>0}$
\begin{enumerate}
    \item\label{lem:cGM:1} $M_n\subset G_{c^G_M n}$
    \item\label{lem:cGM:2} $M_{-n}= (G_{-n})_M$.
\end{enumerate}
\end{lemma}

\section{Effective generation of cuspidal representations -- the cuspidal case of \Cref{thm:supph}}\label{sec:cusp}

\begin{notation}
\RamiF{For a compact $l$-group $K$, denote by $\mdef{e_K}$ the probability Haar measure on $K$. If
     $K<G$ is compact open then we can regard ${e_K}$ as an element in} the Hecke algebra $\cH_K(G)$. For $g\in G$ denote by $\mdef{\delta_g}\in \Sc^*(G)$ the delta-distribution supported in $g$.
\end{notation}

\begin{lemma}\label{lem:eff.cusp}
Let $\pi$ be a cuspidal representation of $G$ and let $v\in\pi$. Then there exists
$\mdef{n_{cusp}(v)}\in\Z_{>0}$ such that for every proper parabolic subgroup $P<G$, with
unipotent radical $U$, there exists an open compact subgroup $\mdef{V_{cusp}(v,P)}<U$ with
$$
V_{cusp}(v,P)\subset G_{n_{cusp}(v)}\qquad\text{and}\qquad e_{V_{cusp}(v,P)}v=0 .
$$
\end{lemma}
\RamiA{
\begin{proof}[Idea of the proof]
    For a fixed \EitanE{parabolic subgroup} $P\EitanE{<G}$ the result is well known. We use the compactness of the flag variety to get a uniform result.
\end{proof}
}

\begin{proof}
Fix first a proper parabolic subgroup $P<G$ with unipotent radical $U$. Since $\pi$ is cuspidal
there is an open compact
subgroup $V<U$ with $e_Vv=0$. See \cite[III.2.9, VI.2.1]{Ren}.

It remains to bound $V$ uniformly as $P$ varies. Let $K_v<G$ be an open compact subgroup fixing
$v$, \DimaA{let $Q_1,\dots,Q_s$ be representatives for the finitely many $K_v$-conjugacy classes of proper
parabolic subgroups, and let $U_i$ be the unipotent radical of $Q_i$.
For every \DimaB{$i$ we }apply the previous paragraph to the proper parabolic subgroup
$Q_i$, \DimaB{and} obtain an open
compact subgroup $V_{i}<U_i$ with $e_{V_{i}}v=0$. The set
$\bigcup_{i}K_vV_{i}K_v$ is compact.
\RamiA{By \EitanE{their definition the subsets} $G_n$ (\Cref{not:settings}\eqref{it:Gn}) form an open increasing cover of $G$. 
Hence} there is
$n_{cusp}(v)\in\Z_{>0}$ with
$$
\bigcup_{i}K_vV_{i}K_v\subset G_{n_{cusp}(v)} .
$$

Now let $P<G$ be an arbitrary proper parabolic subgroup, with unipotent radical $U$. Then
$P=hQ_ih^{-1}$ for some $i$ and some $h\in K_v$.
\RamiA{Set
$$V_{cusp}(v,P):=hV_{i}h^{-1}.$$}
Note that $V_{cusp}(v,P)$ is an open compact subgroup of $U$, and it is contained in
$K_vV_{i}K_v\subset G_{n_{cusp}(v)}$. Finally: 
$$
e_{V_{cusp}(v,P)}v=\RamiA{h\,e_{V_{i}}h^{-1}v=}h\,e_{V_{i}}v=0.$$}
\end{proof}

\begin{lemma}\label{lem:K.cov.U}
There exists $\mdef{c_{uni}}\in\Z_{>0}$ such that for every
$n,m\in\Z_{>0}$ and every
$$
g\in G\setminus \cZ(G)G_{c_{uni}(n+m+1)},
$$
there exists a proper parabolic subgroup $P<G$, with unipotent radical
$U$, such that
$$
g^{-1}G_{-n}g\supset U\cap G_m .
$$
\end{lemma}
\RamiA{
\begin{proof}[Idea of the proof]
We first use the Cartan decomposition to reduce the statement to the case when $g\in T^+$. Then we use the lower bound on the size of $g$ in order to find a simple root that is large on $g$. This defines a parabolic $P$ that satisfies the desired condition thanks to \Cref{lem:U.contract}\eqref{lem:U.contract:1}.
\end{proof}

}
\begin{proof}
If $\Delta=\varnothing$, then $G/\cZ(G)$ is
compact, so after increasing $c_{uni}$ the set
$G\setminus \cZ(G)G_{c_{uni}(n+m+1)}$ is empty and there is nothing
to prove. Assume $\Delta\neq\varnothing$.

By the Cartan decomposition write $g=k_1ak_2$ with $k_1,k_2\in G_0$ and
$a\in T^+$. Since $G_{-n}$ is normalized by $G_0$ and $G_m$ is
$G_0$-bi-invariant, it suffices to find a standard proper parabolic
subgroup $Q$, with unipotent radical $U_Q$, such that
$$
U_Q\cap G_m\subset a^{-1}G_{-n}a ;
$$
then $P:=k_2^{-1}Qk_2$ and $U:=k_2^{-1}U_Qk_2$ satisfy
$U\cap G_m=k_2^{-1}(U_Q\cap G_m)k_2\subset k_2^{-1}a^{-1}G_{-n}ak_2
=g^{-1}G_{-n}g$.

Let $c_{cnj}$ be the constant of \Cref{lem:U.contract} and put
$$
L:=c_{cnj}(n+m).
$$
Let $c_\Theta\in\Z_{>0}$ satisfy $\Theta(a)\leq
c_\Theta\max_{\alpha\in\Delta}\alpha(a)$ for every $a\in T^+$. Such a $c_\Theta$
exists: every element of $\Phi_\bfB$ is a non-negative integral combination of $\Delta$, hence so
is $\Theta=\fN\sum_{\beta\in\Phi_\bfB}\beta$, say $\Theta=\sum_{\alpha\in\Delta}c_\alpha\alpha$
with $c_\alpha\in\Z_{\geq0}$; since $\alpha(a)\geq 0$ for every $\alpha\in\Delta$ and every
$a\in T^+$, any $c_\Theta\geq\sum_{\alpha\in\Delta}c_\alpha$ satisfies the
inequality. Set $R:=c_\Theta L$ and choose $c_{uni}$ so large
that $c_Z(R+1)\leq c_{uni}(n+m+1)$ for all $n,m$, where $c_Z$ is
the constant from \Cref{lem:Theta.ZG.bound}. Since $g\notin
\cZ(G)G_{c_{uni}(n+m+1)}$, \Cref{lem:Theta.ZG.bound} gives
$\Theta(a)>R$, hence there is $\alpha\in\Delta$ with $\alpha(a)>L$.

Let $Q=Q_\alpha$ be the standard maximal parabolic subgroup attached to
$\Delta\smallsetminus\{\alpha\}$ under the bijection between standard parabolic subgroups
and subsets of $\Delta$ (\cite[V.4.2]{Ren}); then $\Phi_Q$ consists of those
$\beta\in\Phi_\bfB$ whose expansion in $\Delta$ has a positive coefficient at
$\alpha$. As $a\in T^+$, every element of $\Delta$ is non-negative on $a$, so
$$
\beta(a)\geq\alpha(a)>L=c_{cnj}(n+m).
$$
Applying \Cref{lem:U.contract}\eqref{lem:U.contract:1} with $k=n+m$ we get
$$
a\,(U_Q\cap G_m)\,a^{-1}\subset a\,(U_Q\cap G_{n+m})\,a^{-1}\subset U_Q\cap G_{-(n+m)}
\subset G_{-n},
$$
that is $U_Q\cap G_m\subset a^{-1}G_{-n}a$, which is what we needed.
\end{proof}

\begin{prop}\label{prop:van.vec}
    Let $\pi$ be  an irreducible cuspidal representation of $G$. Let $v\in \pi$. Then there exists $\mdef{c_{supp,v}}\in \Z_{>0}$ s.t. for any $n\in \Z_{>0}$ and any $g\in G\smallsetminus G_{nc_{supp,v}} \cZ(G)$ we have:
$$e_{G_{-n}} g v=0.$$   
\end{prop}
\RamiA{
\begin{proof}[Idea of the proof]
Note that for the unipotent radical $U$ of a parabolic $P<G$ we have 
$$e_{G_{-n}} g v=g(e_{g^{-1}G_{-n}g})  v=g(e_{g^{-1}G_{-n}g})(e_{U\cap g^{-1}G_{-n}g}) v.$$

By \Cref{lem:K.cov.U} we can find $P$ s.t. 
$U\cap g^{-1}G_{-n}g$ will be large enough. By \Cref{lem:eff.cusp} this implies the vanishing of the r.h.s.
\end{proof}
}
\begin{proof}
Write $m:=n_{cusp}(v)$ for the constant of \Cref{lem:eff.cusp}, let $c_{uni}$ be the constant of
\Cref{lem:K.cov.U}, and put $c_{supp,v}:=c_{uni}(m+2)$.

Let $n\in\Z_{>0}$ and $g\in G\smallsetminus G_{nc_{supp,v}}\cZ(G)$. Since $n\geq1$ we have
$c_{supp,v}n\geq c_{uni}(n+m+1)$, so $g\notin \cZ(G)G_{c_{uni}(n+m+1)}$ and \Cref{lem:K.cov.U}
provides a proper parabolic subgroup $P<G$, with unipotent radical $U$, such that
$$
K:=g^{-1}G_{-n}g\supset U\cap G_m .
$$
Put $V:=V_{cusp}(v,P)$. By \Cref{lem:eff.cusp} it is an open compact subgroup of $U$ contained in
$G_m$, hence $V\subset U\cap G_m\subset K$, and $e_Vv=0$.


\RamiA{Finally, 
$$
e_{G_{-n}}gv=ge_Kv=g\int_K\int_V\pi(ku)v\,du\,dk=g\int_K\pi(k)\bigl(e_Vv\bigr)\,dk=0. 
$$}
\end{proof}
\begin{notation}
    For $X\subset G$ and an open compact $K<G$, denote by
    $\cH_{K}^X(G)\subset \cH_{K}(G)$  the subspace of measures supported in $X$.
\end{notation}
\begin{cor}\label{cor:van.act}
    Let $\pi$ be an irreducible cuspidal representation of $G$.  Let $v\in \pi$. Then there is $\mdef{c_{van,v}}\in \Z_{>0}$ s.t.  
    for any $n\in \Z_{>0}$ and any 
    $$h\in \cH_{G_{-n}}^{G  \smallsetminus \cZ(G)G_{c_{van,v}\cdot n}}(G)$$  we have $h v=0$.
\end{cor}
\begin{proof}
Let $n_0$ be s.t. $G_{-n_0}v=v$
$ $
\begin{enumerate}[{Case} 1.]
\item $h=e_{G_{-n}}\delta_{g} e_{G_{-n}}$, for $g\in G$ and $n>n_0$.\\
Follows immediately from the previous Proposition (\Cref{prop:van.vec}), with
$c_{van,v}:=c_{supp,v}$.
\item $n>n_0$.\\
Follows from the previous case.
\item General case.\\
Follows from the previous case. Indeed, let $c$ be the constant $c_{van,v}$ from the previous case. Let $c_{van,v}:=c(n_0+1)$.
For any $n\in \Z_{>0}$ we have: 
\begin{equation}\label{cor:van.act:1}
c_{van,v}n \geq \max(cn,c(n_0+1))    
\end{equation}

Now take 
any $n\in \Z_{>0}$ and any $$h\in \cH_{G_{-n}}^{G  \smallsetminus \cZ(G)G_{c_{van,v}\cdot n}}(G)$$

If $n> n_0$ we are done by the previous case. Otherwise we have 
$$h\in \cH_{G_{-n}}^{G  \smallsetminus \cZ(G)G_{c_{van,v}\cdot n}}(G)\subset \cH_{G_{-(n_0+1)}}^{G  \smallsetminus \cZ(G)G_{c_{van,v}\cdot n}}(G)\overset{\eqref{cor:van.act:1}}{\subset} \cH_{G_{-(n_0+1)}}^{G  \smallsetminus \cZ(G)G_{c\cdot ({n_0+1})}}(G).$$   Again we are done by previous case.
\end{enumerate}
    
\end{proof}

We are now ready to prove the cuspidal case of \Cref{thm:supph}. 
\begin{cor}\label{thm:supph.cusp}
Let $\pi$ be an irreducible cuspidal representation. 
Then there exist  $\mdef{v_{cgen}(\pi)}\in \pi$ and $\mdef{c_{cgen}(\pi)}\in\Z_{>0}$ with
$v_{cgen}(\pi)\in \pi^{G_{-k}}$ for every $k\in\Z_{>0}$ such that $\pi^{G_{-k}}\neq 0$, and s.t. for any:
 \begin{itemize}
     \item  $n\in  \Z_{>0}$
     \item $w\in \pi^{G_{-n}}$ 
 \end{itemize}
there exists
$h\in \cH_{G_{-n}}(G)$ s.t.
    $\supp(h)\subset G_{c_{cgen}(\pi)\cdot n}$ and
    $hv_{cgen}(\pi)=w$.
\end{cor}
\begin{proof}
Take $v_{cgen}(\pi)\in\pi\smallsetminus 0$ s.t. for all $k$ satisfying $\pi^{G_{-k}}\neq 0$ we have $v_{cgen}(\pi)\in \pi^{G_{-k}}$, and let $c_{cgen}(\pi):=c_{van,v_{cgen}(\pi)}$ be the constant of \Cref{cor:van.act} for this vector.

If $\pi^{G_{-n}}= 0$ then we are done.  Otherwise we have 

\begin{align*}
w\in \pi^{G_{-n}}
&= \cH_{G_{-n}}(G)v_{cgen}(\pi)\\
&=\bigl(\cH_{G_{-n}}^{\cZ(G)G_{c_{cgen}(\pi)n}}(G)\oplus
\cH_{G_{-n}}^{G\smallsetminus \cZ(G)G_{c_{cgen}(\pi)n}}(G)\bigr)v_{cgen}(\pi)\\
&\overset{\text{Cor~\ref{cor:van.act}}}{=}
\cH_{G_{-n}}^{\cZ(G)G_{c_{cgen}(\pi)n}}(G)\,v_{cgen}(\pi)\\
&=\cH_{G_{-n}}^{G_{c_{cgen}(\pi)n}}(G)\,\bigl(\cZ(G)v_{cgen}(\pi)\bigr)\\
&=\cH_{G_{-n}}^{G_{c_{cgen}(\pi)n}}(G)\,v_{cgen}(\pi) .
\end{align*}
\end{proof}

\section{Effective generation of irreducible representations -- the proof of \Cref{thm:supph}}\label{sec:gen}
We  deduce  the general case of \Cref{thm:supph} from the cuspidal case using the fact that any representation is a quotient of a parabolic induction of a cuspidal one.

\subsection{Generation of induction}
In this subsection we study the effective  generation of the parabolic induction  of \EitanE{an irreducible} representation $\rho$ in terms of the effective generation of $\rho$.

Throughout the subsection we fix a  standard parabolic $\mdef{P}<G$ with a standard Levi decomposition $P=\mdef{M}\mdef{U}$. We denote by $\mdef{pr_M}:P\to M$ the projection.
Denote by $\mdef{U^-}$ the unipotent radical of $P^-$.

We also fix an irreducible representation $\mdef{\rho}$ of $M$.

\RamiA{In this subsection we construct a collection of maps $\iota_L:\rho\to i_M^G(\rho)$ that depend\EitanE{s} on a compact subgroup $L<G$ satisfying certain conditions (see \Cref{not:iota}\eqref{it:Madm} below). We prove two main statements on these maps:
\begin{itemize}
    \item The image of $\iota_{\EitanE{L}}$ effectively generates the \EitanE{target} (\Cref{lem:iota.gen}).
    \item Certain compatibility of $\iota_{\EitanE{L}}$ with the action of $M$ and with changing $L$ (\Cref{cor:c.iota}).
\end{itemize}
}

\begin{notation}\label{not:iota}
$ $
\begin{enumerate}

    \item\label{it:Madm} We call a group $L<G$ \tdef{$M$-admissible} if:
    \begin{itemize}
        \item $LP$ is open. 
        \item $L$ is compact. 
        \item $L\cap P=1$.
    \end{itemize}    
    \item  For such $L$, define $\mdef{\iota_{M,L,\rho}}:\rho\to i_M^G(\rho)$ by
    $$\RamiA{\iota_{L}:=}\iota_{M,L,\rho}(v)(x)=\begin{cases}
    p^{-1}v &\text{ if } x=kp \text{ where } k\in L,p\in P\\
    0 &\text{ otherwise }
    \end{cases}.$$    
\end{enumerate}
\end{notation}

\Cref{lem:iva.decomp} gives us:
\begin{cor}\label{lem:Gn.M.adm}
    For $n>0$ we have that $(G_{-n})^-$ is $M$-admissible.
\end{cor}
\begin{lemma}\label{lem:iota.inv}
    For every $n\in \Z_{>0}$ and every vector $v\in \rho^{(G_{-n})_M}$, the vector
    $\iota_{(G_{-n})^-}(v)$ is invariant under $G_{-n}$.
\end{lemma}
\begin{proof}
Write $K:=G_{-n}$ and $f:=\iota_{K^-}(v)$; we have to show that $f(gx)=f(x)$ for every $g\in K$ and every $x\in G$. By \Cref{lem:iva.decomp} we have $K=K^-(K\cap P)$, hence $KK^-P=K^-P$; since $f$ is supported on $K^-P$, both sides vanish unless $x\in K^-P$. Let $g\in K$ and let $x=lp$ with $l\in K^-$ and $p\in P$. Write $gl=l_1p_1$ with $l_1\in K^-$ and $p_1\in K\cap P$. Then $pr_M(p_1)\in K_M$, so
$$f(gx)=f(l_1p_1p)=(p_1p)^{-1}v=p^{-1}\bigl(pr_M(p_1)^{-1}v\bigr)=p^{-1}v=f(x).$$
\end{proof}
\RamiA{The first main property of $\iota_L$ is the following:}
\begin{lemma}\label{lem:iota.gen}
Let $n\in \Z_{>0}$. Then
    $$\cH^{G_0}_{G_{-n}}(G)\cdot ( \iota_{(G_{-n})^-}(\rho^{(G_{-n})_M}))=i_M^G(\rho)^{G_{-n}}.$$
\end{lemma}
\begin{proof}
    For $x\in G_0$ we write $[x]$ for its class in $\cF:=G_0/(G_{-n}P\cap G_0)$.
    Let $f\in i_M^G(\rho)^{G_{-n}}$. We can write
\begin{align*}
f
&= 1_{G} \cdot f
 = \sum_{[g]\in \cF} 1_{gG_{-n}P }\cdot f =\\
&= \sum_{[g]\in \cF} 1_{gG_{-n}P}\cdot gg^{-1}f
 = \sum_{[g]\in \cF} g\left(1_{G_{-n}P}\cdot g^{-1}f\right) \overset{\text{Lem~\ref{lem:iva.decomp}}}{=}\\
& =\sum_{[g]\in \cF}
g\left(1_{(G_{-n})^-P}\cdot g^{-1}f\right) =\\
&= \sum_{[g]\in \cF}
g\,\iota_{(G_{-n})^-}\left((g^{-1}f)(1)\right) =\\
&= \sum_{[g]\in \cF}
g\,\iota_{(G_{-n})^-}\left(f(g)\right) =\\
&= \sum_{[g]\in \cF}
(e_{G_{-n}}\delta_g e_{G_{-n}})\,\iota_{(G_{-n})^-}\left(f(g)\right).
\end{align*}
\end{proof}
\RamiA{The second main property of $\iota_L$ is the following:}
\begin{prop}\label{cor:c.iota}
There is $\mdef{c_{\iota}}\in \Z_{>0}$ such that for
any $n\leq m\leq k\in \Z_{>0}$ and any $v\in \rho^{(G_{-n})_M}$, we have
\begin{enumerate}
    \item\label{cor:c.iota:1}
    $
    \iota_{(G_{-m})^-}(\cS^*_{M_k}(M)v)
    \subset
    \cS^*_{G_{c_{\iota}k}}(G)
    \left(\iota_{(G_{-n})^-}(v)\right).
    $

    \item\label{cor:c.iota:2}
    $
    i_M^G(\rho)^{G_{-m}} \cap
    \iota_{(G_{-m})^-}
    \left(\cH^{M_k}_{M_{-m}}(M)v\right)
    \subset
    \cH_{G_{-m}}^{G_{c_{\iota}k}}(G)
    \left(\iota_{(G_{-n})^-}(v)\right).
    $
\end{enumerate}
\end{prop}
\RamiA{
\begin{proof}[Idea of the proof]
 \eqref{cor:c.iota:2} follows easily from \eqref{cor:c.iota:1}. So it is enough to prove \eqref{cor:c.iota:1}.
    We observe the following compatibility:
    \begin{equation}\label{eq:naiv.comp}
        g\iota_{L}(v)=\iota_{gLg^{-1}}(gv),
    \end{equation}
    for $g\in M$. See \Cref{lem:iota.M}.

    So we \EitanE{are} reduce\EitanE{d} to showing that for $g\in M_k$ the vector 
    $\iota_{g^{-1}(G_{-m})^-g}(v)$  is in the r.h.s. (after changing the constant $c_\iota$). 

    We prove \Cref{lem:iota.enlarge} that states that for $L<L'$ the map $\iota_{L'}$  can be expressed as a combination of shifts of $\iota_{L}$. Thus we can replace $g^{-1}(G_{-m})^-g$ by any smaller group. 

    This allows us to assume  that $g\in A_\bfM$ (by virtue of \Cref{lem:A.cont}\eqref{lem:A.cont:1} below). Replacing $g$ with another (possibly somewhat larger) element of $A_\bfM$, we can replace $m$ with $n$. Applying \eqref{eq:naiv.comp} again (in the opposite direction) we deduce the statement from the fact that $A_\bfM$ acts by scalars on $\rho$.
\end{proof}
For the actual proof we will need some preparations.
}

\begin{lemma}\label{lem:iota.enlarge}
Let $L<L'$ be $M$-admissible subgroups. Then $L$ is open and of finite index in $L'$. Moreover, for any set $R\subset L'$ of representatives for the left cosets $L'/L$ and any $v\in \rho$,
$$
\iota_{L'}(v)=\sum_{r\in R} r\circ \iota_{L}(v).
$$
\end{lemma}
\begin{proof}
We have $L'\cap LP=L$: if $l'=lp$ with $l\in L$ and $p\in P$, then
$p=l^{-1}l'\in L'\cap P=1$. Since $LP$ is open in $G$, the subgroup $L$ is
therefore open in $L'$, and $L'$ is compact, so $L'/L$ is finite.

Let $R\subset L'$ be a set of representatives for the left cosets $L'/L$.
We prove the equality pointwise. Let $x\in G$. If $x\notin L'P$, then
$x\notin rLP$ for every $r\in R$, so both sides vanish at $x$.

Assume now that $x\in L'P$. Since $L'\cap P=1$, we can write uniquely
$x=l'p$ with $l'\in L'$ and $p\in P$. Write uniquely $l'=rl$ with
$r\in R$ and $l\in L$. Then the only non-zero summand on the right-hand side
is the one indexed by $r$, and
$$
(r\circ \iota_{L}(v))(x)
=
\iota_{L}(v)(r^{-1}x)
=
\iota_{L}(v)(lp)
=
p^{-1}v.
$$
On the other hand, by definition,
$$
\iota_{L'}(v)(x)=\iota_{L'}(v)(l'p)=p^{-1}v.
$$
Therefore the two sides agree for every $x\in G$.
\end{proof}

\begin{lemma}\label{lem:iota.M}
    Let $L$ be an $M$-admissible subgroup and $m\in M$.
    Then $mLm^{-1}$ is $M$-admissible, and
    $$m\iota_{L}(v)=\iota_{mLm^{-1}}(mv).$$
\end{lemma}
\begin{proof}
Set $L'=mLm^{-1}$. Since conjugation by $m$ is an analytic homeomorphism, $L'$ is compact and analytic. Also, since $m$ normalizes $P$,
$$
L'P=mL(m^{-1}P)=mLP
$$
is open, and
$$
L'\cap P
=
m(L\cap m^{-1}Pm)m^{-1}
=
m(L\cap P)m^{-1}
=
1.
$$
Thus $L'$ is $M$-admissible.

Let $f=\iota_{L}(v)$. Both $mf$ and $\iota_{L'}(mv)$ are supported on $L'P=mLP$. Let $x\in L'P$, and write
$$
x=l'p
$$
with $l'=mlm^{-1}$, $l\in L$, and $p\in P$. Then
$$
(mf)(x)
=
f(m^{-1}x)
=
f(lm^{-1}p)
=
(m^{-1}p)^{-1}v
=
p^{-1}mv.
$$
This is exactly $\iota_{L'}(mv)(x)$. Hence
$$
m\iota_{L}(v)=\iota_{mLm^{-1}}(mv).
$$
\end{proof}

\begin{lemma}\label{lem:A.cont}
There is $\mdef{c_{cont}}\in \Z_{>0}$ such that for every $n\in\Z_{>0}$:
\begin{enumerate}
\item\label{lem:A.cont:1} For every $m\in M$ there is $a\in A_\bfM$ with
$ov_G(a)\leq c_{cont}(ov_G(m)+n)$ and
$$
a^{-1}(G_{-n})^-a\subset m^{-1}(G_{-n})^-m .
$$
\item\label{lem:A.cont:2} For every $m\in\Z_{>0}$ with $m>n$ there is $a\in A_\bfM$ with
$ov_G(a)\leq c_{cont}m$ and
$$
a(G_{-n})^-a^{-1}<G_{-m} .
$$
\end{enumerate}
\end{lemma}

\begin{proof}
If $U^-=1$ then $(G_{-n})^-=1$ and $a:=1$ satisfies both assertions. Assume $U^-\neq1$. Then $P$
is proper, so there is $\lambda\in\Lambda(A_\bfM)$ with $\langle\beta,\lambda\rangle>0$ for
every $\beta\in\Phi_P$; fix such a $\lambda$. Then $\langle\beta,\lambda\rangle<0$ for every
$\beta\in\Phi_{P^-}$, since $\Phi_{P^-}=-\Phi_P$; put
$$
b:=\min_{\beta\in\Phi_{P^-}}\bigl(-\langle\beta,\lambda\rangle\bigr)>0 .
$$
For $r\in\Z_{\geq0}$ we have $\beta(\varpi^{-r\lambda})=-r\langle\beta,\lambda\rangle\geq rb$ for
every $\beta\in\Phi_{P^-}$, and, $\lambda$ being fixed, there is $c'\in\Z_{>0}$ with
$ov_G(\varpi^{r\lambda})\leq c'r$ and $ov_G(\varpi^{-r\lambda})\leq c'r$. Note that $P^-$ is the opposite of a standard parabolic, that
$U^-$ is its unipotent radical, and that its Levi subgroup is $M$, so \Cref{lem:U.contract}
applies to $P^-$.

\eqref{lem:A.cont:1} Let $k\in\Z_{>0}$ be minimal with $m\in G_k$, set
$r:=\lceil c_{cnj}k/b\rceil$ and $a:=\varpi^{r\lambda}\in A_\bfM$. Then
$\beta(a^{-1})\geq rb\geq c_{cnj}k$ for every $\beta\in\Phi_{P^-}$, so
\Cref{lem:U.contract}\eqref{lem:U.contract:2}, applied to $a^{-1}$ and to $g=m$, gives
$$
a^{-1}(G_{-n})^-a\subset m^{-1}(G_{-n})^-m .
$$
Since $k\leq ov_G(m)+1\leq ov_G(m)+n$, we get $ov_G(a)\leq c'r\leq c_{cont}(ov_G(m)+n)$ after
choosing $c_{cont}$ large enough.

\eqref{lem:A.cont:2} Set $r:=\lceil c_{cnj}m/b\rceil$ and $a:=\varpi^{-r\lambda}\in A_\bfM$. Then
$\beta(a)\geq rb\geq c_{cnj}m$ for every $\beta\in\Phi_{P^-}$, and
$(G_{-n})^-<U^-\cap G_m$ because $G_{-n}<G_0\subset G_m$, so
\Cref{lem:U.contract}\eqref{lem:U.contract:1}, applied with $k=m$, gives
$$
a(G_{-n})^-a^{-1}\subset U^-\cap G_{-m}<G_{-m} .
$$
Finally $ov_G(a)\leq c'r\leq c_{cont}m$ after enlarging $c_{cont}$, since $m\geq1$.
\end{proof}

\begin{lemma}\label{lem:iota.m.gen}
    There is $\mdef{c_{\iota cov}}\in \Z_{>0}$ s.t. for any $n,k\in \Z_{>0}$, $m\in M_k$ and $v\in \rho$ we have 
    $$\iota_{(G_{-n})^-}(mv)\in \cS^*_{G_{c_{\iota cov}(k+n)}}(G)(\iota_{(G_{-n})^-}(v)).$$
\end{lemma}
\begin{proof}
Set $c_{\iota cov}= 4 c_{cont} c_M^G$. By \Cref{lem:A.cont}\eqref{lem:A.cont:1} we can choose $a\in A_\bfM$ s.t.
\begin{equation}\label{eq:iota.ov.a}
ov_{G}(a)\leq c_{cont} (ov_{G}(m)+n)
\end{equation}
and
    $$a^{-1}(G_{-n})^- a \subset m^{-1}(G_{-n})^- m.$$

Denote $R:=(m^{-1}(G_{-n})^-m)/(a^{-1}(G_{-n})^-a)$. Now we have:
\begin{align*}
\iota_{(G_{-n})^-}(mv)
&\overset{\text{Lem~\ref{lem:iota.M}}}{=}m\, \iota_{m^{-1}(G_{-n})^-m}(v)\\
&\overset{\text{Lem~\ref{lem:iota.enlarge}}}{=}m \sum_{[l]\in R} l\cdot  \iota_{a^{-1}(G_{-n})^-a}(v)\\
&=m \sum_{[l]\in R} l\cdot a^{-1}a\,  \iota_{a^{-1}(G_{-n})^-a}(v)\\
&\overset{\text{Lem~\ref{lem:iota.M}}}{=}m \sum_{[l]\in R} l\cdot a^{-1}  \iota_{(G_{-n})^-}(av)\\
&=c_\rho(a) \sum_{[l]\in R} m\cdot l\cdot a^{-1}  \iota_{(G_{-n})^-}(v),
\end{align*}
where $c_\rho$ is the central character of $\rho$.

Now, for all $l\in m^{-1}(G_{-n})^-m$ we have
\begin{equation}\label{eq:iota.ov.l}
ov_G(l)\leq 2ov_G(m).
\end{equation}
Indeed
$(G_{-n})^-\subset G_0$ and $ov_G(m^{-1})=ov_G(m)$, so, by the
$G_0$-bi-invariance of the $G_s$,
$$
l\in G_{ov_G(m)}G_0G_{ov_G(m)}=G_{ov_G(m)}G_{ov_G(m)}
\overset{\text{Lem~\ref{lem:Gn.product}}}{\subset} G_{2ov_G(m)}.
$$
Hence
\begin{align*}
ov_G(m\cdot l\cdot a^{-1})
&\overset{\text{Lem~\ref{lem:Gn.product}}}{\leq} ov_G(m)+ ov_G(l)+ov_G(a)\\
&\overset{\eqref{eq:iota.ov.l}}{\leq} 3ov_G(m)+ov_G(a)\\
&\overset{\eqref{eq:iota.ov.a}}{\leq} 3ov_G(m)+c_{cont} (ov_G(m)+n)\\
&\overset{\text{Lem~\ref{lem:cGM}\eqref{lem:cGM:1}}}{\leq}
3 c_M^G k+c_{cont} (c_M^Gk+n)\leq c_{\iota cov} (n+k).
\end{align*}
\end{proof}

\begin{cor}\label{cor:iota.M.gen}  
    For any $n,k\in \Z_{>0}$ we have $$\iota_{(G_{-n})^-}(\cS^*_{M_k}(M)v)\subset \cS^*_{G_{c_{\iota cov}(k+n)}}(G)(\iota_{(G_{-n})^-}(v)).$$
\end{cor}

\begin{lemma}\label{lem:iota.cont}
For any $n\leq m\in \Z_{>0}$ and any $v\in \rho$ we have
$$
\iota_{(G_{-m})^-}(v)
\in
\cS^*_{G_{c_{cont}m}}(G)
\left(\iota_{(G_{-n})^-}(v)\right).
$$
\end{lemma}

\begin{proof}
By \Cref{lem:A.cont}\eqref{lem:A.cont:2}, there is $a\in A_\bfM$ with 
\begin{equation}\label{eq:mag.a}
 ov_G(a)\leq c_{cont}m   
\end{equation}
 such that
$$
a\RamiA{(G_{-n})^-}a^{-1}\subset G_{-m}.
$$
Since $a\in A_\bfM$ normalizes $U^-$, we also have
$$
a\RamiA{(G_{-n})^-}a^{-1}\subset U^-.
$$
Hence
$$
a\RamiA{(G_{-n})^-}a^{-1}\subset G_{-m}\cap U^-=\RamiA{(G_{-m})^-}.
$$

By \Cref{lem:Gn.M.adm,lem:iota.M}, both $a\RamiA{(G_{-n})^-}a^{-1}$ and $\RamiA{(G_{-m})^-}$ are
$M$-admissible. Let
$$
R\subset \RamiA{(G_{-m})^-}
$$
be a set of representatives for the left cosets
\begin{equation}\label{eq:decomp.a}
\RamiA{(G_{-m})^-}/(a\RamiA{(G_{-n})^-}a^{-1}).
\end{equation}
By \Cref{lem:iota.enlarge},
$$
\iota_{\RamiA{(G_{-m})^-}}(v)
=
\sum_{r\in R} r\cdot \iota_{a\RamiA{(G_{-n})^-}a^{-1}}(v).
$$

Since $a\in A_\bfM\subset \cZ(M)$ and $\rho$ is irreducible, Schur's lemma
implies that $a$ acts on $\rho$ by a scalar. Denote this scalar by
$c_\rho(a)\in \C^\times$. 
\RamiA{
We get
\begin{align*}
\iota_{(G_{-m})^-}(v)
&\overset{\eqref{eq:decomp.a}}{=}
\sum_{r\in R} r\cdot \iota_{a(G_{-n})^-a^{-1}}(v)=\\
&=c_\rho(a)^{-1}\sum_{r\in R} r\cdot\iota_{a{(G_{-n})^-}a^{-1}}(av)
\overset{\text{Lem }\ref{lem:iota.M}}{=}
\\
&=
c_\rho(a)^{-1}\sum_{r\in R} r a\cdot \iota_{{(G_{-n})^-}}(v)
\overset{\text{Lem~\ref{lem:Gn.product}},\ \eqref{eq:mag.a}}{\in}
\\
&
\in
\cS^*_{G_{c_{cont}m}}(G)\left(\iota_{(G_{-n})^-}(v)\right),
\end{align*}
as required.
}
\end{proof}

\begin{proof}[Proof of \Cref{cor:c.iota}]
We first prove \eqref{cor:c.iota:1}. The convolution of distributions supported in $X$ and $Y$ is supported in
$XY$; this is used in the last step below. We have
\begin{align*}
\iota_{(G_{-m})^-}(\cS^*_{M_k}(M)v)
&\overset{\text{Cor~\ref{cor:iota.M.gen}}}{\subset}
\cS^*_{G_{c_{\iota cov}(k+m)}}(G)
\left(\iota_{(G_{-m})^-}(v)\right)\\
&\overset{\text{Lem~\ref{lem:iota.cont}}}{\subset}
\cS^*_{G_{c_{\iota cov}(k+m)}}(G)
\cS^*_{G_{c_{cont}m}}(G)
\left(\iota_{(G_{-n})^-}(v)\right)\\
&\overset{\text{Lem~\ref{lem:Gn.product}}}{\subset}
\cS^*_{G_{c_{\iota cov}(k+m)+c_{cont}m}}(G)
\left(\iota_{(G_{-n})^-}(v)\right),
\end{align*}
where the first step is applied with $m$ in the role of $n$. Since $m\leq k$ we have
$c_{\iota cov}(k+m)+c_{cont}m\leq (2c_{\iota cov}+c_{cont})k$, so taking
$$
c_{\iota}\geq 2c_{\iota cov}+c_{cont}
$$
proves \eqref{cor:c.iota:1}.

For \eqref{cor:c.iota:2} let $u$ be an element of the left hand side and write, using
\eqref{cor:c.iota:1}, $u=\mu\,\iota_{(G_{-n})^-}(v)$ with $\mu\in \cS^*_{G_{c_\iota k}}(G)$. By
\Cref{lem:iota.inv} the vector $\iota_{(G_{-n})^-}(v)$ is $G_{-n}$-invariant, hence
$G_{-m}$-invariant since $G_{-m}<G_{-n}$. As $u$ is $G_{-m}$-invariant as well, we get
$$
u=e_{G_{-m}}\,\mu\, e_{G_{-m}}\bigl(\iota_{(G_{-n})^-}(v)\bigr),
$$
and $e_{G_{-m}}\mu e_{G_{-m}}\in \cH_{G_{-m}}^{G_{c_\iota k}}(G)$ because $G_{-m}<G_0$ and
$G_{c_\iota k}$ is $G_0$-bi-invariant.
\end{proof}

\subsection{The general case of  \Cref{thm:supph}}

\begin{proof}[Proof of \Cref{thm:supph}]
$ $
\begin{enumerate}[{Case} 1.]

    \item $\pi=i_M^G(\rho)$  where $M<G$ is a standard Levi and  $\rho$ is an irreducible cuspidal representation of $M$. Here we do not assume that $\pi$ is irreducible.\\

    Let $v_{cgen}(\rho)$ be as in
    \Cref{thm:supph.cusp}.
    Let $n_0$ be the minimal integer  such that     
    $$v_{cgen}(\rho)\in \rho^{M_{-n_0}}\overset{\text{Lem \ref{lem:cGM}\eqref{lem:cGM:2}}}{=}\rho^{(G_{-n_0})_M}.$$  
Note that   \Cref{thm:supph.cusp} assures that for any $n<n_0$ we have
\begin{equation}\label{eq:rho.M.van}
\rho^{M_{-n}}=0.    
\end{equation}

    Let $c_{gen}(\pi)=c_{\iota}c_{cgen}(\rho)$
and $v_{gen}(\pi):=\iota_{M,(G_{-n_0})^-,\rho}(v_{cgen}(\rho))$.

For $n<n_0$ we have 

\begin{align*}
i_M^G(\rho)^{G_{-n}}
&\overset{\text{Lem } \ref{lem:iota.gen}}{=}\cH^{G_0}_{G_{-n}}(G)\cdot ( \iota_{M,(G_{-n})^-,\rho}(\rho^{(G_{-n})_M}))\\
&\overset{\text{Lem \ref{lem:cGM}\eqref{lem:cGM:2}}}{=} \cH^{G_0}_{G_{-n}}(G)\cdot ( \iota_{M,(G_{-n})^-,\rho}(\rho^{M_{-n}}))\\
&\overset{\eqref{eq:rho.M.van}}{=}0
\end{align*}
and there is nothing to prove. Let $n\geq n_0$. Then:
\begin{align*}
\pi^{G_{-n}}
&\overset{\text{Lem~\ref{lem:iota.gen}}}{=}
\cH_{G_{-n}}^{G_0}(G)\Bigl(\pi^{G_{-n}}\cap\iota_{M,(G_{-n})^-,\rho}\bigl(\rho^{(G_{-n})_M}\bigr)\Bigr)\\
&\overset{\text{Lem~\ref{lem:cGM}\eqref{lem:cGM:2}}}{=}
\cH_{G_{-n}}^{G_0}(G)\Bigl(\pi^{G_{-n}}\cap\iota_{M,(G_{-n})^-,\rho}\bigl(\rho^{M_{-n}}\bigr)\Bigr)\\
&\overset{\text{Thm~\ref{thm:supph.cusp}}}{\subset}
\cH_{G_{-n}}^{G_0}(G)\Bigl(\pi^{G_{-n}}\cap\iota_{M,(G_{-n})^-,\rho}\bigl(\cH_{M_{-n}}^{M_{c_{cgen}(\rho)n}}(M)\,v_{cgen}(\rho)\bigr)\Bigr)\\
&\overset{\text{Cor~\ref{cor:c.iota}\eqref{cor:c.iota:2}}}{\subset}
\cH_{G_{-n}}^{G_0}(G)\,\cH_{G_{-n}}^{G_{c_{\iota}c_{cgen}(\rho)n}}(G)\bigl(\iota_{M,(G_{-n_0})^-,\rho}(v_{cgen}(\rho))\bigr)\\
&\overset{\text{Lem~\ref{lem:Gn.product}}}{\subset}
\cH_{G_{-n}}^{G_{c_{gen}(\pi)n}}(G)\bigl(v_{gen}(\pi)\bigr),
\end{align*}
where \Cref{cor:c.iota}\eqref{cor:c.iota:2} is applied with $n_0$, $n$ and $c_{cgen}(\rho)n$ in the
roles of $n$, $m$ and $k$.

    \item General case.\\
     Follows from the previous case using the fact that any irreducible representation is  a quotient  of  a parabolic induction of a cuspidal representation.
\end{enumerate}
\end{proof}

\subsection{Applications of \Cref{thm:supph} \DimaB{and the proof of \Cref{thm:suppm}}}\label{ssec:app.supph}

{
\begin{cor}\label{cor:supph.fl.alt}
      Let $\pi$ be a representation of $G$ of finite length.
Then
there exist $\mdef[v_{gen,i}]{v_{gen,1}}(\pi),\dots, v_{gen,N}(\pi)\in \pi$
and  $c_{gen}(\pi)$ s.t.
for any:
\begin{itemize}
    \item $n\in \Z_{>0}$
    \item $w\in \pi^{G_{-n}}$
\end{itemize}
there exist
$h_1,\dots,h_N\in \cH_{G_{-n}}(G)$ s.t.
\begin{enumerate}
    \item $\supp(h_i)\subset G_{c_{gen}(\pi)n}$ for every $i$
    \item $\sum h_iv_{gen,i}(\pi)=w$.
\end{enumerate}
\end{cor}
\begin{proof}
\RamiB{The proof is} by induction on the length $L$ of $\pi$. For $L=0$ there is nothing to prove.

Let $L>0$. Choose a maximal proper subrepresentation $\pi'<\pi$; then $\pi'':=\pi/\pi'$ is
irreducible and $\pi'$ is of length $L-1$. Let $q:\pi\to\pi''$ be the projection, let
$v_{gen}(\pi'')$ and $c_{gen}(\pi'')$ be as in \Cref{thm:supph} for $\pi''$, and fix $v_0\in\pi$
with $q(v_0)=v_{gen}(\pi'')$. Let $v'_1,\dots,v'_{N'}\in\pi'$ and $c'$ be given by the induction
hypothesis for $\pi'$.

Let $n\in\Z_{>0}$ and $w\in\pi^{G_{-n}}$. Then $q(w)\in(\pi'')^{G_{-n}}$, so by
\Cref{thm:supph} there is $h_0\in\cH_{G_{-n}}(G)$ with $\supp(h_0)\subset G_{c_{gen}(\pi'')n}$
and $h_0v_{gen}(\pi'')=q(w)$. As $h_0$ is left $G_{-n}$-invariant we have
$e_{G_{-n}}h_0=h_0$, so $h_0v_0\in\pi^{G_{-n}}$, and
$$
q(h_0v_0)=h_0q(v_0)=h_0v_{gen}(\pi'')=q(w),
$$
whence
$$
w-h_0v_0\in(\pi')^{G_{-n}} .
$$

By the induction hypothesis $w-h_0v_0=\sum_kh'_kv'_k$ with $h'_k\in\cH_{G_{-n}}(G)$ and
$\supp(h'_k)\subset G_{c'n}$. Thus $w=h_0v_0+\sum_kh'_kv'_k$, and the assertion holds for $\pi$
with the vectors $v_0,v'_1,\dots,v'_{N'}$ and the constant $c_{gen}(\pi):=\max(c',c_{gen}(\pi''))$.
\end{proof}
\RamiB{
The last corollary implies the following analog, in which we replace $\cH(G)$ with $\C[G]$.
}
\begin{cor}\label{cor:gen.fin.len.alt}
      Let $\pi$ be a representation of $G$ of finite length.
\RamiB{Then
for any $n\in \Z_{>0}$ and   $w\in \pi^{G_{-n}}$}
there exist
$h_1,\dots,h_N\in \C[G]$ s.t.
\begin{enumerate}
    \item $\supp(h_i)\subset G_{c_{gen}(\pi)n}$ for every $i$
    \item $\RamiB{\sum_{i=1}^N}h_iv_{gen,i}(\pi)=w$. \EitanE{Here  $v_{gen,i}$ are the vectors from \Cref{cor:supph.fl.alt}.}
\end{enumerate}
\end{cor}
\RamiB{
\begin{proof}
\RamiB{Denote  $v_i:=v_{gen,i}(\pi)$.
F}ix an open compact subgroup $\RamiB{K}<\RamiB{G_{-n}}$ fixing all of
$v_1,\dots,v_N$.

Fix $n\in\Z_{>0}$ and $w\in\pi^{G_{-n}}$. By \Cref{cor:supph.fl.alt} there are
$\tilde h_1,\dots,\tilde h_N\in\cH_{G_{-n}}(G)$ with
\begin{equation}\label{eq:h.v.w}
\sum_i\tilde h_iv_i=w
\quad \text{ and } \quad 
\supp(\tilde h_i)\subset G_{c_{gen}(\pi)n}.
\end{equation}
Fix $i$. Choosing representatives $g_{i1},\dots,g_{ir_i}\in G_{c_{gen}(\pi)n}$ for the double
cosets in $\RamiB{K}\backslash\supp(\tilde h_i)/\RamiB{K}$, we may write
\begin{equation}\label{eq:h.dec}
\tilde h_i=\sum_{j=1}^{r_i} a_{ij}\,\RamiB{e_{K}}\delta_{g_{ij}}\RamiB{e_{K}},
\qquad a_{ij}\in\C.
\end{equation}
The smooth vector
$$w':=\sum_{i=1}^N\sum_{j=1}^{r_i}a_{ij}g_{ij}v_i$$ is fixed by an open subgroup $K'<K$ of finite
index, so

\begin{align*}
w&\overset{\eqref{eq:h.v.w}}{=}\sum_{i=1}^{N}\tilde h_iv_i\\
&\overset{\eqref{eq:h.dec}}{=}\sum_{i=1}^{N}\sum_{j=1}^{r_i} a_{ij}\,\RamiB{e_{K}}\delta_{g_{ij}}\RamiB{e_{K}}v_i\\
&={e_{K}}\sum_{i=1}^{N}\sum_{j=1}^{r_i} a_{ij}{g_{ij}}v_i\\
&=\frac{1}{[K:K']}\sum_{[k]\in K/K'}
{e_{K'}}\sum_{i=1}^{N}\sum_{j=1}^{r_i} a_{ij}k{g_{ij}}v_i\\
&=\frac{1}{[K:K']}\sum_{[k]\in K/K'}
\sum_{i=1}^{N}\sum_{j=1}^{r_i} a_{ij}k{g_{ij}}v_i\\
&=\sum_{i=1}^{N}\frac{1}{[K:K']}\sum_{[k]\in K/K'}
\sum_{j=1}^{r_i} a_{ij}k{g_{ij}}v_i
.
\end{align*}
Set
$$
h_i:=\frac{1}{[K:K']}\sum_{[k]\in K/K'}
\sum_{j=1}^{r_i} a_{ij}k{g_{ij}}
\in
\C[G]
$$
We have
$\supp(h_{i})\subset 
G_{-n}
G_{c_{gen}n}
G_{-n}=G_{c_{gen}n}$
and 
$w=\sum_{i=1}^{N}h_iv_i$,
as required.
\end{proof}
}
}
\begin{proof}[Proof of \Cref{thm:suppm}]
The contragredient $\tilde\pi$ is again of finite length, and it is cuspidal.
Let $v_{gen,1}(\pi),\dots,v_{gen,N}(\pi)$ and $c_{gen}(\pi)$, respectively
$v_{gen,1}(\tilde\pi),\dots,v_{gen,N'}(\tilde\pi)$ and $c_{gen}(\tilde\pi)$, be as in
\Cref{cor:gen.fin.len.alt} for $\pi$ and for $\tilde\pi$. For all $i,j$ the function
$m_{v_{gen,i}(\pi),\,v_{gen,j}(\tilde\pi)}$ is a matrix coefficient of $\pi$, which is cuspidal,
and so has compact support modulo $\cZ(G)$ (\DimaA{see e.g.} \cite[VI.2.1]{Ren}). As the pairs $(i,j)$ are finite
in number, there is $N_0\in\Z_{>0}$ with
\begin{equation}\label{eq:supp.gen}
\Supp\bigl(m_{v_{gen,i}(\pi),\,v_{gen,j}(\tilde\pi)}\bigr)\subset\cZ(G)G_{N_0}
\qquad\text{for all }i,j .
\end{equation}
Put $c_{supp}(\pi):=c_{gen}(\pi)+c_{gen}(\tilde\pi)+N_0$.

Fix $n\in\Z_{>0}$, $v\in\pi^{G_{-n}}$ and $w\in\tilde\pi^{G_{-n}}$, and write, by
\Cref{cor:gen.fin.len.alt},
$$
v=\RamiA{\sum_{i=1}^N}h_iv_{gen,i}(\pi),\qquad w=\RamiA{\sum_{j=1}^{N'}}\tilde h_jv_{gen,j}(\tilde\pi),
$$
with $h_i=\sum_\alpha a_{i\alpha}\delta_{x_{i\alpha}}$ and
$\tilde h_j=\sum_\beta b_{j\beta}\delta_{z_{j\beta}}$ in $\C[G]$, where
\begin{equation}\label{eq:gen.xz}
x_{i\alpha},\ x_{i\alpha}^{-1}\in G_{c_{gen}(\pi)n},
\qquad z_{j\beta}\in G_{c_{gen}(\tilde\pi)n},
\end{equation}
the membership of $x_{i\alpha}^{-1}$ because $ov_G(y^{-1})=ov_G(y)$. Since
$\langle \pi(y)a,\tilde\pi(z)b\rangle=\langle\pi(z^{-1}y)a,b\rangle$ we get, for every $g\in G$,
\begin{equation}\label{eq:mvw.exp}
m_{v,w}(g)=\sum_{i,j,\alpha,\beta}a_{i\alpha}b_{j\beta}\;
m_{v_{gen,i}(\pi),\,v_{gen,j}(\tilde\pi)}\bigl(z_{j\beta}^{-1}\,g\,x_{i\alpha}\bigr) .
\end{equation}
\RamiA{
We get:
\begin{align*}
\supp(m_{v,w})
&\overset{\eqref{eq:mvw.exp}}{\subset}
\bigcup_{i,j,\alpha,\beta} z_{j\beta}
\left(\supp(m_{v_{gen,i}(\pi),\,v_{gen,j}(\tilde\pi)})\right)x_{i\alpha}^{-1}
\overset{\eqref{eq:supp.gen}}{\subset}\\
&\subset
\bigcup_{i,j,\alpha,\beta} z_{j\beta} \left(\cZ(G)G_{N_0}\right)x_{i\alpha}^{-1}
\overset{\eqref{eq:gen.xz}}{\subset}\\
&\subset
\bigcup_{i,j,\alpha,\beta} G_{c_{gen}(\tilde \pi)n} \left(\cZ(G)G_{N_0}\right)G_{c_{gen}(\pi)n}
=\\
&=\cZ(G)\,G_{c_{gen}(\tilde\pi)n}G_{N_0}G_{c_{gen}(\pi)n}
\overset{\text{Lem~\ref{lem:Gn.product}}}{\subset}\\
&\subset
\cZ(G)G_{(c_{gen}(\pi)+c_{gen}(\tilde\pi))n+N_0} \overset{n\geq 1}{\subset}\\
&\subset\cZ(G)G_{c_{supp}(\pi)n},
\end{align*}
as required.}
\end{proof}

\section{Asymptotics of matrix coefficients of representations of finite length -- the proof of \Cref{thm:mcexp}}\label{sec:asm}
Throughout this section we will fix $\pi$ and $\lambda$ as in \Cref{thm:mcexp}.

We denote by $\mdef{\bfM_\lambda} \subset \mdef{\bfP_\lambda}\subset \bfG$ the corresponding Levi and parabolic subgroups. Let $\mdef{\bfU_\lambda}$ be the unipotent radical of $\bfP_\lambda$. Let $\mdef{M_\lambda}, \mdef{P_\lambda}$ and $\mdef{U_\lambda}$ be the corresponding groups of $F$-points.
\subsection{\RamiA{Proof of \Cref{thm:eff.jak.cass}}}
\RamiA{
Let us first introduce \EitanE{some} terminology which will allow us to streamline and explicate the statement of \Cref{thm:eff.jak.cass}.

}

\begin{notation}\label{not:cP}
Let $\mdef{\cP}$ be the (monic) minimal polynomial of the action of $\varpi^\lambda$ on the unnormalized Jacquet module $\mdef{r_{M_\lambda}^G}(\pi)$.
\end{notation}

\begin{defn}
We say that $v\in \pi$ is 
\tdef{$n$-admissible} if 
$$\cP(\delta_{\varpi^\lambda})(v)=\sum_{i=1}^{N}(u_iw_i-w_i)$$
for some $u_i\in U_\lambda \cap G_n$ and $w_i\in \pi$.

Here we consider $\delta_{\varpi^{\lambda}}$ as an element in the group algebra $\C[G]$.
\end{defn}

\RamiA{
\begin{proof}[Idea of the proof of \Cref{thm:eff.jak.cass}]
    We have to find a constant $c$ s.t. \DimaB{every} vector $v$ of depth $n$ is $nc$-admissible. A classical argument shows that any vector is $k$-admissible for some $k$ (depending on the vector), see \Cref{lem:v.adm} below. We show that acting on a $k$-admissible vector $v$ with an element $g$ of a given magnitude makes it $k'$-admissible with a linear control on $k'$, see \Cref{cor:sh}. Finally we use  \Cref{thm:supph} (more specifically \Cref{cor:gen.fin.len.alt}) to obtain \DimaB{every} vector of depth $n$ from the fixed generators  using elements of the group of a bounded magnitude (linearly in $n$). 

    The proof of \Cref{cor:sh} is based on an effective version of the Iwasawa decomposition (see \Cref{lem.Iv}). Namely, we prove it separately for  $g\in M_\lambda$ (see \Cref{lem:adm.shift.M}) and for $g\in U_\lambda$ (see \Cref{lem:adm.shift.U}). We also use the fact that $G_0v$ is finite to deduce a uniform bound on the admissibility on the vectors in $G_0v$, see \Cref{cor:v.adm}.
\end{proof}
}

The following is obvious:
\begin{lemma}\label{lem:adm.lin}
Let $n\leq n'\in\Z_{>0}$, let $a_1,a_2\in\C$ and let $v_1,v_2\in \pi$ be $n$-admissible
vectors. Then
\begin{enumerate}
\item\label{lem:adm.lin:1} $v_1$ is $n'$-admissible.
\item\label{lem:adm.lin:2} The vector $a_1v_1+a_2v_2$ is $n$-admissible.
\end{enumerate}
\end{lemma}

\begin{lemma}\label{lem:v.adm}
    For any $v\in \pi$ there is $\mdef{\operatorname{Adm}_{v}}\in \Z_{>0}$ s.t. $v$ is $\operatorname{Adm}_{v}$-admissible.
\end{lemma}
\RamiA{This is a classical statement, for completeness we include  the proof here.}
\begin{proof}
Write $r_{M_\lambda}^G(\pi)=\pi/\pi(U_\lambda)$ for the (unnormalized) $U_\lambda$-coinvariants, where $\pi(U_\lambda):=\operatorname{span}_\C\{uw-w:\ u\in U_\lambda,\ w\in\pi\}$, and let $j:\pi\twoheadrightarrow r_{M_\lambda}^G(\pi)$ be the projection. 
\RamiA{By definition of $\cP$ we get that}
\RamiA{
$$
\cP(\delta_{\varpi^\lambda})j(v)=0.
$$
Thus}
$$
\cP(\delta_{\varpi^\lambda})(v)\in\ker j=\pi(U_\lambda).
$$
Write $\cP(\delta_{\varpi^\lambda})v=\sum_{i=1}^N(u_iw_i-w_i)$ with $u_i\in U_\lambda$ and $w_i\in\pi$. 
\RamiA{There} is $k$ with $u_1,\dots,u_N\in G_k$.
\RamiA{So, taking} $\operatorname{Adm}_v:=k$
\RamiA{implies the assertion}.
\end{proof}
\begin{cor}\label{cor:v.adm}
    For any $v\in \pi$ there is $\mdef{\operatorname{Adm}'_{v}}\in \Z_{>0}$ s.t. any $w\in G_0 v$ is $\operatorname{Adm}'_{v}$-admissible.
\end{cor}

\begin{lemma}\label{lem:adm.shift.M}
For every $n\in\Z_{>0}$, every $n$-admissible $v\in\pi$ and every $m\in M_\lambda$, the vector $mv$ is $(n+2\,ov_G(m))$-admissible.
\end{lemma}
\begin{proof}
Since $\varpi^\lambda\in A_{\bfM_\lambda}\subset \cZ(M_\lambda)$, the operator $\cP(\delta_{\varpi^\lambda})$ commutes with the action of $m$. Writing $\cP(\delta_{\varpi^\lambda})v=\sum_i(u_iw_i-w_i)$ with $u_i\in U_\lambda\cap G_n$,
$$
\cP(\delta_{\varpi^\lambda})(mv)=m\,\cP(\delta_{\varpi^\lambda})v=\sum_i\bigl((mu_im^{-1})(mw_i)-(mw_i)\bigr).
$$
\RamiA{It is enough to show that $mu_im^{-1}\in U_\lambda\cap G_{n+2ov_G(m)}$.}

\RamiA{Note that} $mu_im^{-1}\in U_\lambda$ (as $M_\lambda$ normalizes $U_\lambda$).

\RamiA{By} \Cref{lem:Gn.product}
$$
ov_G(mu_im^{-1})\le 2\,ov_G(m)+ov_G(u_i)\le 2\,ov_G(m)+n,
$$
\RamiA{We obtain} $mu_im^{-1}\in U_\lambda\cap G_{n+2ov_G(m)}$ \RamiA{as required}. 
\end{proof}

\begin{lemma}\label{lem:adm.shift.U}
There exists $\mdef{c_{shu}}\in\Z_{>0}$ s.t.
for every $n\in\Z_{>0}$, every $n$-admissible $v\in\pi$ and every $u\in U_\lambda$, the vector $uv$ is 
$\bigl(n+ov_G(u)\bigr) c_{shu}$-admissible.
\end{lemma}
\begin{proof}
Set $d:=\deg\cP$ and take \RamiA{$$c_{shu}:= 
1+2d\,ov_G(\varpi^\lambda).$$
}
\RamiA{Since $v$ is $n$-admissible, we  can write
$$\cP(\delta_{\varpi^\lambda})v=\sum_{j=1}^{r}(u_jw_j-w_j)$$ 
with $u_j\in U_\lambda\cap G_n$.
}

Write $\cP=\sum_{i=0}^d a_ix^i$.
\RamiA{ 
We have 
\begin{align*}
\cP(\delta_{\varpi^\lambda})(uv)&=\cP(\delta_{\varpi^\lambda})v+\cP(\delta_{\varpi^\lambda})(uv)-\cP(\delta_{\varpi^\lambda})v\\
&=\sum_{j=1}^{r}(u_jw_j-w_j)+\sum_{i=0}^d \bigl(a_i\varpi^{i\lambda}uv-a_i\varpi^{i\lambda}v\bigr)\\
&=\sum_{j=1}^{r}(u_jw_j-w_j)+\sum_{i=0}^d\bigl(\bigl(\varpi^{i\lambda}u\varpi^{-i\lambda}\bigr)(a_i\varpi^{i\lambda}v)-(a_i\varpi^{i\lambda}v)\bigr).
\end{align*}
Also, for $j\in \{1,\dots,r\}, i\in \{0,\dots,d\} $  we have:
\begin{itemize}
    \item 
$u_j, \varpi^{i\lambda}u\varpi^{-i\lambda}\in U_\lambda$ 
    \item $w_j, a_i\varpi^{i\lambda} v \in \pi$
    \item
$ov_G(u_j)\leq n\leq \bigl(n+ov_G(u)\bigr) c_{shu}$
    \item  $ov_G(\varpi^{i\lambda}u\varpi^{-i\lambda})
\leq 
\bigl(n+ov_G(u)\bigr) c_{shu}$. Indeed:
\begin{align*}
ov_G(\varpi^{i\lambda}u\varpi^{-i\lambda})
&\overset{\text{Lem \ref{lem:Gn.product}}}{\leq}  2iov_G(\varpi^{\lambda})+ov_G(u)\\
&\leq  2dov_G(\varpi^{\lambda})+ov_G(u)
\leq 
\bigl(n+ov_G(u)\bigr) c_{shu}
\end{align*}
\end{itemize}
So we are done.
}
\end{proof}

\begin{lemma}[\Cref{ssec:eff.Iwa}]\label{lem.Iv}
    There exists $\mdef{c_I}\in \Z_{>0}$ s.t. for any $n\in \Z_{>0}$ we have
    $$G_n\subset (G_{c_In}\cap U_{\lambda})(G_{c_In}\cap M_\lambda) G_{0}.$$
\end{lemma}

\begin{cor}\label{cor:sh}
For any $v\in \pi$ there is a constant $\mdef{c_{sh,v}}$ s.t. for any $g\in G$ the vector $gv$ is $c_{sh,v}(1+ov_{G}(g))$-admissible.
\end{cor}

\begin{proof}
Let $n_0:=\operatorname{Adm}'_v$ be the constant from \Cref{cor:v.adm}, let $c_I$ and $c_{shu}$ be the constants from \Cref{lem.Iv} and \Cref{lem:adm.shift.U}, and set $c_{sh,v}:=(n_0+3c_I)c_{shu}$.

Fix $g\in G$ and put $n:=ov_G(g)+1\ge 1$, so $g\in G_{ov_G(g)}\subset G_n$. By \Cref{lem.Iv} we may write $g=u\,m\,k$ with $u\in U_\lambda\cap G_{c_In}$, $m\in M_\lambda\cap G_{c_In}$ and $k\in G_0$; in particular $ov_G(u),ov_G(m)\le c_In$.

By \Cref{cor:v.adm} the vector $kv\in G_0v$ is $n_0$-admissible. By \Cref{lem:adm.shift.M} the vector $m(kv)$ is $(n_0+2\,ov_G(m))$-admissible. By \Cref{lem:adm.shift.U} the vector $gv=u\bigl(m(kv)\bigr)$ is $\bigl(n_0+2\,ov_G(m)+ov_G(u)\bigr)c_{shu}$-admissible. Finally
\begin{align*}
\bigl(n_0+2\,ov_G(m)+ov_G(u)\bigr)c_{shu}
&\le (n_0+3c_In)c_{shu}\\
&\overset{n\ge1}{\le}(n_0+3c_I)c_{shu}\,n=c_{sh,v}(1+ov_G(g)),
\end{align*}
so by \Cref{lem:adm.lin}\eqref{lem:adm.lin:1} the vector $gv$ is $c_{sh,v}(1+ov_G(g))$-admissible.
\end{proof}
Now we are ready to prove \Cref{thm:eff.jak.cass}:
\begin{thm}\label{prop:mcexp}
There is $\mdef{c_{asm0}}\in\Z_{>0}$ s.t. for any $n\in\Z_{>0}$ every $v\in\pi^{G_{-n}}$ is $c_{asm0}n$-admissible.
\end{thm}
\begin{proof}
Let $v_{gen,1}(\pi),\dots,v_{gen,N}(\pi)$ and $c_g:=c_{gen}(\pi)$ be the generators and the constant of \Cref{cor:gen.fin.len.alt}, let $c_0:=\max_ic_{sh,v_{gen,i}(\pi)}$ be the largest of the constants of \Cref{cor:sh} for these finitely many vectors, and set $c_{asm0}:=c_0(1+c_g)$.

Fix $n\in\Z_{>0}$ and $v\in\pi^{G_{-n}}$. By \Cref{cor:gen.fin.len.alt} there are $h_i=\sum_\alpha b_{i\alpha}\delta_{g_{i\alpha}}\in\C[G]$ with $g_{i\alpha}\in G_{c_gn}$ and $\sum_ih_iv_{gen,i}(\pi)=v$; that is,
$$v=\sum_{i,\alpha} b_{i\alpha}\, g_{i\alpha} v_{gen,i}(\pi),\qquad b_{i\alpha}\in\C,\ g_{i\alpha}\in G_{c_gn}.$$

For each $i,\alpha$ we have $ov_G(g_{i\alpha})\le c_gn$, so by \Cref{cor:sh} the vector $g_{i\alpha} v_{gen,i}(\pi)$ is $c_{sh,v_{gen,i}(\pi)}(1+ov_G(g_{i\alpha}))$-admissible, hence $c_0(1+c_gn)$-admissible by \Cref{lem:adm.lin}\eqref{lem:adm.lin:1}. By \Cref{lem:adm.lin}\eqref{lem:adm.lin:2} the vector $v$ is $c_0(1+c_gn)$-admissible, and
$$c_0(1+c_gn)\overset{n\ge1}{\le}c_0(1+c_g)n=c_{asm0}n,$$
so $v$ is $c_{asm0}n$-admissible.
\end{proof}
\subsection{\RamiA{Proof of \Cref{thm:mcexp}}}
\RamiB{
\begin{proof}[Idea of the proof]
    By \Cref{thm:eff.jak.cass} it is enough to show that $$\langle \varpi^{k\lambda}(u v' -v') \varpi^{-k\lambda},w\rangle=0,$$
    for $u\in U_\lambda \cap G_{m}$ where $m$ is bounded linearly in the depth $n$ of $w$, and $k$ big enough (linearly in $n$).

    This follows from the fact that 
    $$\varpi^{k\lambda}\,(U_\lambda\cap G_n)\,\varpi^{-k\lambda}\subseteq G_{-n}.$$
    See \Cref{lem:Ulam.contract} below.
\end{proof}
}
\begin{lemma}\label{lem:Ulam.contract}
There is $\mdef{c_{ucont}}\in\Z_{>0}$ such that for all $m,n\in\Z_{>0}$ with $n\geq m c_{ucont}$ we have
$$\varpi^{n\lambda}\,(U_\lambda\cap G_m)\,\varpi^{-n\lambda}\subseteq G_{-m}.$$
\end{lemma}
\begin{proof}
By \Cref{not:settings}\eqref{it:Phi.P}, $\Phi_{\bfP_\lambda}$ is the set of
$\beta\in\Phi_\bfB$ with $\langle\beta,\lambda\rangle>0$.
\RamiA{If $\varpi^\lambda$ is central then $U_\lambda$ is trivial and statement is obvious. Otherwise $\Phi_\bfB$ is non-empty.}
Set
$d_\lambda:=\min_{\beta\in\Phi_{\bfP_\lambda}}\langle\beta,\lambda\rangle>0$. Choose
$c_{ucont}\in\Z_{>0}$ with $c_{ucont}d_\lambda\geq c_{cnj}$, where $c_{cnj}$ is the constant of
\Cref{lem:U.contract}. For $n\geq mc_{ucont}$ and every $\beta\in\Phi_{\bfP_\lambda}$,
$$
\beta(\varpi^{n\lambda})=n\langle\beta,\lambda\rangle\geq mc_{ucont}d_\lambda\geq c_{cnj}m ,
$$
so \Cref{lem:U.contract}\eqref{lem:U.contract:1}, applied to $Q=\bfP_\lambda$ with $k=m$ and
$a=\varpi^{n\lambda}$, gives
$$\varpi^{n\lambda}(U_\lambda\cap G_m)\varpi^{-n\lambda}\subseteq U_\lambda\cap G_{-m}\subseteq G_{-m}.$$
\end{proof}

\begin{proof}[Proof of \Cref{thm:mcexp}]
Set $\cP_{\pi,\lambda}:=\cP$  and $\cP=\sum_{i=0}^{\deg\cP}a_ix^i$; let $c_{asm0}$ be the constant of \Cref{prop:mcexp}, $c_{ucont}$ that of \Cref{lem:Ulam.contract}, and put $c_{asm}(\pi,\lambda):=c_{ucont}c_{asm0}$.

Fix $n\in\Z_{>0}$, $v\in\pi^{G_{-n}}$, $w\in\tilde\pi^{G_{-n}}$ and $k>c_{asm}(\pi,\lambda)n$. 
\RamiA{We have:}
\begin{equation}\label{eq:mc.P}
\sum_{i=0}^{\deg\cP}a_i\,m_{v,w}(\varpi^{(k+i)\lambda})
=\RamiA{\sum_{i=0}^{\deg\cP}a_i\big\langle\varpi^{k\lambda}\varpi^{i\lambda}v,\,w\big\rangle
=}\big\langle\varpi^{k\lambda}\,\cP(\delta_{\varpi^\lambda})v,\,w\big\rangle .
\end{equation}

By \Cref{prop:mcexp} the vector $v$ is $c_{asm0}n$-admissible, so
\begin{equation}\label{eq:P.unip}
\cP(\delta_{\varpi^\lambda})v=\sum_{j=1}^{N}(u_jw_j-w_j),\qquad u_j\in U_\lambda\cap G_{c_{asm0}n},\ w_j\in\pi.
\end{equation}
Put $u_j^{(k)}:=\varpi^{k\lambda}u_j\varpi^{-k\lambda}$. As $k>c_{asm}(\pi,\lambda)n=c_{ucont}(c_{asm0}n)$, 
\RamiA{we get:}
\begin{equation}\label{eq:G-n}    
u_j^{(k)}\in\varpi^{k\lambda}(U_\lambda\cap G_{c_{asm0}n})\varpi^{-k\lambda}\overset{\text{Lem \ref{lem:Ulam.contract}}}{\subseteq} G_{-c_{asm0}n}\subseteq G_{-n}.
\end{equation}

Hence:
\RamiA{
\begin{equation}\label{eq:pair.inv}
\begin{split}
\big\langle \varpi^{k\lambda} u_j w_j,\,w\big\rangle
&=\big\langle u_j^{(k)}  \varpi^{k\lambda} w_j,\,w\big\rangle\\
&=\big\langle \varpi^{k\lambda} w_j,\,(u_j^{(k)})^{-1} w\big\rangle\\
&\overset{\eqref{eq:G-n},\, w\in\tilde\pi^{G_{-n}}}{=}\big\langle \varpi^{k\lambda} w_j,\,w\big\rangle.
\end{split}
\end{equation}
}
\RamiA{Therefore:
\begin{align*}
\sum_{i=0}^{\deg\cP_{\pi,\lambda}}a_i\,m_{v,w}(\varpi^{(k+i)\lambda})
&\overset{\eqref{eq:mc.P}}{=}\big\langle \varpi^{k\lambda}\,\cP(\delta_{\varpi^\lambda})v,\,w\big\rangle
\overset{\eqref{eq:P.unip}}{=}\\
&=\big\langle \varpi^{k\lambda}\,\sum_{j=1}^{N}(u_jw_j-w_j),\,w\big\rangle=\\
&=\sum_{j=1}^N\left(\big\langle \varpi^{k\lambda}  u_j w_j,\,w\big\rangle-
\big\langle \varpi^{k\lambda}   w_j,\,w\big\rangle\right)
\overset{\eqref{eq:pair.inv}}{=}0.
\end{align*}
}
\end{proof}

\section{Norms}\label{sec:Norms}
In this section, we recall basic parts of the theory of norms developed in \cite[\S18]{Kot}.
We will use the following notions from \cite[\S18]{Kot}.
\begin{enumerate}
    \item An \tdef{abstract norm} on a set $Z$ is a real-valued
function $\|\cdot \|_Z$ on $Z$ such that $\|x\|_Z\geq 1$ for all $x \in Z$. 

\item For two abstract norms $\|\cdot\|_Z^1,\|\cdot\|_Z^2$ on $Z$ we say that $\|x\|_Z^1\mdef{\prec}\|x\|_Z^2$ if there is a constant $c>1$ such that $\|x\|_Z^1<c(\|x\|_Z^2)^c$.

    \item We say that two abstract norms $\|\cdot\|_Z^1,\|\cdot\|_Z^2$ on $Z$ are equivalent if  $\|x\|_Z^1\prec\|x\|_Z^2\prec\|x\|_Z^1$.
We denote this as 
     $$\|\cdot\|_Z^1 \mdef{\sim} \|\cdot\|_Z^2.$$
    \item Let $\bfM$ be an algebraic variety. In \cite[\S18]{Kot} there is a definition of a canonical equivalence class of abstract norms on $M=\bfM(F)$. The abstract norms in this class are called \tdef{norms} on $M$.    
    \item For a map $\phi:Z_1\to Z_2$ between sets and abstract norms $\|\cdot\|_{Z_1}, \|\cdot\|_{Z_2}$ on these sets, \cite[\S18]{Kot} defines norms $\phi_*(\|\cdot\|_{Z_1})$ and $\phi^*(\|\cdot\|_{Z_2})$ on $\operatorname{Im}(\phi)$ and $Z_1$  respectively by $\mdef[\phi_*(\vert\vert\cdot\vert\vert)]{\phi_*}(\|\cdot\|_{Z_1})(z):=\inf_{y\in \phi^{-1}({z})}\|y\|_{Z_1}$, and $\mdef[\phi^*(\vert\vert\cdot\vert\vert)]{\phi^*}(\|\cdot\|_{Z_2})(y):=\|\phi(y)\|_{Z_2}$. 
\end{enumerate}
\begin{notation}\label{not:ov.and.ball}

As a rule, we will put the domain of definition of a norm in a subscript in the notation for that norm.
Given an abstract norm $\|\cdot\|_X$ on a set $X$, we denote:
\begin{enumerate}
    \item \label{not:ov.and.ball:1} By $\mdef{ov}_X:X\to \R$ the map given by $ov_X(x)=\log_{\ell}(\|x\|_X)$, where $\mdef{\ell}$ is the size of the residue field of $F$. 
\item For $i \in \mathbb{Z}$, $\mdef[\bullet_i]{X_i}:=\{x\in X\, \vert ov_X(x)\leq  i\}.$
\end{enumerate}
If we consider more than one norm on the same set $X$, we will distinguish between them and the related notation $ov_X$ and $X_i$ using superscripts.
\end{notation}

\begin{notation}
    For two sets $Z_1,Z_2$, and abstract norms $\|\cdot \|_{Z_1},\|\cdot \|_{Z_2}$ on them we define an abstract norm $\mdef{\vert\vert\cdot\vert\vert_{Z_1}\times \vert\vert\cdot\vert\vert_{Z_2}}:=\|\cdot\|_{Z_1\times {Z_2}}$, where 
    $$\|(a,b)\|_{{Z_1}\times {Z_2}}:=\max(\|a\|_{Z_1},\|b\|_{Z_2}).$$
\end{notation}
\RamiA{
\begin{lemma}[{\cite[Proposition 18.1(2)]{Kot}}]\label{lem:nrm.cl.emb}
    Let $\phi:\bfX\to \bfY$ be  a closed embedding of algebraic varieties. Let $X:=\bfX(F)$ and  $Y:=\bfY(F)$.    
    Then $\phi^*(\|\cdot\|_{Y})$  is a norm on $X$.
\end{lemma}
\begin{lemma}[{\cite[Proposition 18.1(4)]{Kot}}]\label{lem:nrm.prd}
    Let $\bfX,\bfY$ be algebraic varieties. 
    Then $\|\cdot\|_{X} \times \|\cdot\|_{Y}$  is a norm on $X\times Y$.
\end{lemma}
}

\begin{lemma}\label{lem:norm.sigma}
The map $g\mapsto \ell^{ov_\sigma(g)}$ is a norm on $G$, where $ov_\sigma$ is as in
\Cref{not:sigma}.
\end{lemma}

\begin{proof}
Put $n:=\dim \cW$ and identify $\GL(\cW)$ with the principal open subset of the affine space
$\operatorname{End}(\cW)$ defined by $\det$.
By \cite[Example 18.1]{Kot} the map
$x\mapsto\max(1,\max_{i,j}|x_{ij}|)$ is a norm on $\operatorname{End}(\cW)$, so by
\cite[Proposition 18.1(5)]{Kot}
$$
\|x\|':=\max\bigl(1,\max_{i,j}|x_{ij}|,|\det x|^{-1}\bigr)
$$
is a norm on $\GL(\cW)$. Write $\|x\|:=\max(1,\max_{i,j}|x_{ij}|,\max_{i,j}|(x^{-1})_{ij}|)$, so that
$\ell^{ov_\sigma(g)}=\|\sigma(g)\|$. As $\det(x^{-1})$ is a sum of products of $n$ entries of
$x^{-1}$, we get $|\det x|^{-1}\leq\|x\|^{n}$ and hence $\|x\|'\leq\|x\|^{n}$; conversely each
$(x^{-1})_{ij}$ is $\pm(\det x)^{-1}$ times a product of $n-1$ entries of $x$, so
$\|x\|\leq(\|x\|')^{n}$. Thus $\|\cdot\|$ and $\|\cdot\|'$ are equivalent, so $\|\cdot\|$ is a norm
on $\GL(\cW)$. Since $\sigma$ is a closed embedding, \RamiA{\Cref{lem:nrm.cl.emb} shows} that the
restriction of $\|\cdot\|$ to $G$ is a norm on $G$.
\end{proof}

\begin{prop}\label{prop:G.norms}
$ $
    \begin{enumerate}
        \item\label{prop:G.norms:1} The map $g\mapsto \ell^{ov_G(g)}$ is a norm on $G$
        \item\label{prop:G.norms:2} The map $g\mapsto \ell^{ov_{G^{rss}}(g)}$ is a norm on $G^{rss}$
        \item\label{prop:G.norms:3} The map $g\mapsto \ell^{ov_{G^{ab}}(g)}$ is a norm on $G^{ab}$.
    \end{enumerate}
\end{prop}
\begin{proof}
\eqref{prop:G.norms:1} By \Cref{lem:norm.sigma} it is enough to show that
the abstract norms $\ell^{ov_G}$ and $\ell^{ov_\sigma}$ are equivalent, that is, that there is
$c\in\Z_{>0}$ with $ov_G\leq c\,ov_\sigma+c$ and $ov_\sigma\leq c\,ov_G+c$. This is \Cref{cor:ov.sigma.G}.

\eqref{prop:G.norms:3} The group $\bfG^{ab}$ is a torus; let $\bfS$ be its maximal split quotient
and $\chi_1,\dots,\chi_r$ a $\Z$-basis of $X^*(\bfS)=X^*(\bfG^{ab})_F$. By
\cite[Example 18.2]{Kot} the map $s\mapsto\max_i\max(|\chi_i(s)|,|\chi_i(s)|^{-1})$ is a norm on
$S(F)$, and by \cite[Lemma 18.6]{Kot} its pullback
$$
h\mapsto \ell^{\,\max_i|\operatorname{val}\chi_i(h)|}
$$
is a norm on $G^{ab}$. By \Cref{not:settings}\eqref{it:0H} the map
$h\mapsto(\operatorname{val}\chi_i(h))_i$ has kernel ${}^0G^{ab}$, so it induces an injection
$\Lambda(G^{ab})\to\Z^r$, whose image is of finite index by \cite[V.2.6]{Ren}. Hence
$\max_i|\operatorname{val}\chi_i(\cdot)|$ and $ov_{\Lambda(G^{ab})}$ are two norms on the
$r$-dimensional real vector space $\Lambda(G^{ab})\otimes_\Z\R$, so they are bounded by constant
multiples of each other, and the two abstract norms on $G^{ab}$ are equivalent.

\eqref{prop:G.norms:2} By \Cref{not:rss}, $\bfG^{rss}$ is the principal open subset of $\bfG$
defined by $\fD$, so by \eqref{prop:G.norms:1} and \cite[Proposition 18.1(5)]{Kot} the map
$$
x\mapsto\max\bigl(\ell^{ov_G(x)},|\fD(x)|^{-1}\bigr)
$$
is a norm on $G^{rss}$. By \Cref{not:settings}\eqref{it:Gn} \RamiA{and \Cref{not:rss} this map is $$x \mapsto \ell^{ov_{G^{rss}}(x)}$$ as required.}
\end{proof}

\begin{notation}
    
We fix the following norms:

\begin{enumerate}
    \item We denote by $\mdef[\vert\vert\cdot\vert\vert_{G}]{\|\cdot\|_{G}}$ the norm given by \Cref{prop:G.norms}\eqref{prop:G.norms:1}
        \item We denote by $\mdef[\vert\vert\cdot\vert\vert_{G^{rss}}]{\|\cdot\|_{G^{rss}}}$ the norm given by \Cref{prop:G.norms}\eqref{prop:G.norms:2}
\item We denote by $\mdef[\vert\vert\cdot\vert\vert_{G^{ad}}]{\|\cdot\|_{G^{ad}}}$ the abstract norm given by $g\mapsto \ell^{ov_{G^{ad}}(g)}$

\item For any other variety $\bfV$ we choose a norm $\mdef[\vert\vert\cdot\vert\vert_{V}]{\|\cdot\|_{V}}$
on $V:=\bfV(F).$  
\end{enumerate}

\end{notation}
Note that this notation makes  
\Cref{not:rss}  and \Cref{not:settings}\eqref{it:ov1}-\eqref{it:ov2}  compatible with
\Cref{not:ov.and.ball}\eqref{not:ov.and.ball:1}.

\section{Rectified algebraic varieties}\label{sec:rect}
In this section, we recall the theory of rectified algebraic varieties and balls and measures on them, introduced in \cite{AGKS0}, and adapt it to our purposes.
\begin{definition}
    Let $\bfX$ be a smooth algebraic variety over $F$.  
    
    \begin{enumerate}
        \item 
    A \tdef{rectification}   of $\bfX$ is a finite open cover $\bfX \subset \bigcup_{\alpha \in I} {\bfU_{\alpha}}$ 
    with 
    closed embeddings ${i_\alpha}:\bfU_{\alpha}\to \bA^M$.    
    \item We will call a rectification \tdef[simple rectification]{simple} if $|I|=1$.
    \item     By a \tdef{rectified variety} we will mean a smooth algebraic variety equipped with a rectification. {By a map or a morphism of such we just mean a morphism of the underlying algebraic varieties.}
    \item     A \tdef{$\mu$-rectification} of  $\bfX$  is a rectification of $\bfX$  together with invertible top differential forms ${\omega_{\alpha}}\in \Omega^{top}(\bfU_{\alpha})$.
    \item We define similarly the notion of a \tdef{$\mu$-rectified} variety, and 
     simple $\mu$-rectification.
    \end{enumerate}
\end{definition}

\RamiC{
\begin{definition}\label[definition]{def:simp balls}
$ $
    Let $(\bfX,i)$ be a simply rectified variety. Then, for any $m\in \Z$ define:  \begin{enumerate}
    \item \label{def:simp balls:1}  
$\mdef{B_m^{\bfX}}:=i^{-1} \left(\varpi^{-m}\cO_F^{M}\right){\subset \bfX(F)}.$
        \item\label{def:simp balls:2} For $x\in \bfX(F)$ define {a ball around $x$:}
        $$\mdef{B_m^{\bfX}(x)}:= i^{-1} \left(i(x)+\varpi^{-m}\cO_F^{M}\right).$$
\end{enumerate}
\end{definition}
}
\RamiC{

\begin{lem}[{\cite[Lemma 3.5]{AGKS0}}]\label{lem:KL}
Let $(\bfX,\bfU_\alpha,i_\alpha)_{\alpha\in I}$ be a rectified variety. Then there
is $s\in\N$ with the following property. For every  $\alpha,\beta\in I$,
every $m\in\N_{>0}$, every
\[
x\in B^{\bfU_\alpha}_{m}\cap B^{\bfU_\beta}_{m},
\]
and every integer $r\ge1$ 
we have 
$$B_{-r-sm}^{\bfU_\alpha}(x)\subset  B_{-r}^{\bfU_\beta}(x).$$
\end{lem}
}

\RamiC{
We choose for every rectified variety an $s$ as in the above lemma and denote it by $\mdef{s_\bfX}$. We choose it in a way that if the rectification is simple then $s=0$.

We are now ready to define balls and measures on our rectified varieties.}
\RamiC{
\begin{definition}\label[definition]{def: balls and measures}
$ $
\begin{enumerate}
    \item \label{def: balls and measures:1}
    Let $(\bfX,\bfU_\alpha,i_\alpha)$ be a rectified variety. Then, for any $m\in \Z$ define:    
            \begin{enumerate}
            \item\label{def: balls and measures:a}
$\mdef{B_m^{\bfX}}:=\bigcup_{\alpha} B_m^{\bfU_\alpha}{\subset \bfX(F)}.$
        \item\label{def: balls and measures:c} For $x\in \bfX(F)$ define {a ball around $x$:}
        $$\mdef{B_m^{\bfX}(x)}:= \bigcup_{\alpha, a \text{ s.t. } x\in B_{a}^{\bfU_\alpha}; a>0} B_{m-as_\bfX  }^{\bfU_\alpha}(x) .$$
        \item    For $\bfZ\subset \bfX$  we define:  $$\mdef{B_m^{\bfX}(\bfZ)}:=\bigcup_{z\in {B^{\bfZ}_\infty}}   B_m^{\bfX}(z).$$ 
    \end{enumerate}
        
        \item 
    Let $(\bfX,\bfU_\alpha,i_\alpha,\omega_\alpha)$ be a \RamiB{$\mu$-rectified} variety. Then, for any {integer $m\in \Z$}  we define
        a measure $\mu_m^{\bfX}$ on $\bfX(F)$ supported on $B_m^{\bfX}$ by: $$\mdef{\mu_m^{\bfX}}:=\sum_\alpha |(\omega_\alpha)_{F}| \cdot 1_{i_\alpha^{-1} \left(\varpi^{-m}\cO_F^{M}\right)}.$$
\end{enumerate}
\end{definition}
}

\begin{rem}\label{rem:ball.vs.norm}
Let $\bfX$ be a rectified variety and $X:=\bfX(F)$.  
    Note that by definition there is a constant $c\in \Z_{>0}$ s.t. for any $n\in \Z_{>0}$ we have $$X_n\subset B_{cn}^{\bfX} \subset X_{c^2n}.$$
\end{rem}

\begin{definition}
    Let $\bfX$ be a rectified variety. Let $m\in\Z_{>0}$. We say that 
    $f\in C^\infty(\bfX(F))$
    is \tdef{$m$-smooth} if for any 
    $x\in \bfX(F)$
    the function $f|_{B_{-m}^{\bfX}(x)}$ 
    is constant.
\end{definition}

The proof of \cite[Theorem 5.7]{AGKS0} can be easily adapted to prove the following:
\begin{thm}\label{lem:push.mes.sm}
    Let $\gamma:\bfX_1 \to \bfX_2$ be a smooth map of {$\mu$-rectified} varieties. Then, there is $\mdef{c_{sub,\gamma}}$ s.t. for any $m\in \Z_{>0}$ and any $m$-smooth function
    $g\in C^\infty_c({\bfX_1}(F))$, there is an $mc_{sub,\gamma}$-smooth function
    $f\in C^\infty_c(B_{m c_{sub,\gamma}}^{\bfX_2})$ such that
    $$
    \gamma_*(g\mu_m^{\bfX_1})=f\cdot \mu_{m c_{sub,\gamma}}^{\bfX_2}.
    $$
\end{thm}

The following lemma compares the balls of
$\bfG^{rss}$ with the cosets of the groups $G_{-n}$; it is used in \Cref{ThmSm}.

\begin{lemma}[\Cref{ssec:app.ball}]\label{lem:coset.ball}
Equip $\bfG^{rss}$ with a simple rectification. Then there is $\mdef{c_{ball}}\in\Z_{>0}$ such that
for every $x\in G^{rss}$ and every $m\in\Z_{>0}$
$$
x\,G_{-c_{ball}(m+ov_{G^{rss}}(x))}\subset B^{\bfG^{rss}}_{-m}(x).
$$
\end{lemma}

\section{Effective version of Harish-Chandra's smoothness theorem  -- the proof of \Cref{thm:smooth}}\label{ThmSm}
Throughout this section we fix an irreducible representation $\pi$ of $G$.

\RamiB{In \cite{HC_sub}, Harish-Chandra proved that the map 
$\mdef{\phi_x}:\bfG\times \bfB\to \bfG$ given by $\phi_x(h,b)=hxh^{-1}b$ is submersive, for every $x \in G^{rss}$. Moreover, he gave a formula for the value $\chi_\pi(x)$ in terms of the measure:
$${f^{K}_x}:=(\phi_x)_*(1_{G_0\times K}\cdot \mu_{G}\boxtimes \mu_B),$$ where $K<B$ is a small enough open compact subgroup.
As a result the character of \DimaB{every} irreducible representation is locally constant on $G^{rss}$. In this
section we study the submersive map $\psi:G \times B \times G^{rss} \to G \times G^{rss}$ 
given by $$\psi\bigl((h,b),x\bigr)=(hxh^{-1}b,\,x).$$ Using the methods of \cite{AGKS0}, more specifically \Cref{lem:push.mes.sm} of the present paper, we quantify Harish-Chandra's argument. This allows us to obtain uniformity in Harish-Chandra's argument and as a consequence an effective smoothness result for the restriction of the character to $G^{rss}$ .
}
\begin{notation}
    For any $x\in G^{rss}$ define  $\mdef{\phi_x}:\bfG\times \bfB\to \bfG$ to be the map given by $\phi_x(h,b)=hxh^{-1}b$.
\end{notation}

\begin{lemma}[{\cite[Theorem 1]{HC_sub}, \cite[Appendix~B, by G.~Prasad]{AD02}}]\label{lem:HC.submersion}
    $\phi_x$ is a submersion.
\end{lemma}

\begin{notation}[{cf. \cite[\S 4]{HC_sub}}]\label{thm:HC.trace}
Let $x\in G^{rss}$
\begin{itemize}
    \item 
Denote by $\mdef{\cT_x}:\pi\to\pi$ the linear map $$\cT_{x}:=\int_{k\in G_{0}}\pi(kxk^{-1})dk$$ where $dk$ is the normalized Haar measure on $G_{0}$.
 \item  
 For an open compact subgroup $K<B$ denote by
    $$\mdef{f^{K}_x}:=(\phi_x)_*(\RamiF{ e_{G_0}\boxtimes e_K})$$ 
    where $\mu_{G}$ and $\mu_B$ are some (left) Haar measures.
\end{itemize}
\end{notation}
\EitanB{We recall the following qualitative result of Harish-Chandra:}
\begin{theorem}[{\cite[\S 4]{HC_sub}}]\label{thm:HC.char}
    For any $x\in G^{rss}$ the operator $\mdef{\cT_x}:\pi\to\pi$ is of finite rank and the function
    $x\mapsto \operatorname{tr}(\cT_x)$ represents $\chi_\pi|_{G^{rss}}$.
\end{theorem}
\begin{lemma}[{\cite[page 99, 5th display formula]{HC_sub}}]\label{lem : T_x}
    For any  $n\in \Z_{>0}$ there is an open compact $K<B$ s.t. for any
    \begin{itemize}
        \item $x\in G^{rss}$ 
        \item $m\in T^+$ 
        \item $k\in G_0$ 
        \item $v\in \pi^{G_{-n}}$
    \end{itemize}
    We have 
    $$\cT_x(km\cdot v)=k \cdot f^{K}_x \cdot m\cdot  v.$$    
\end{lemma}
\begin{cor}\label{cor:fK.det}
    There is an open compact $K<B$ s.t. for any 
    $x_1,x_2\in G^{rss}$,
    if $f^{K}_{x_1}=f^{K}_{x_2}$, then:
    \begin{enumerate}
        \item $\cT_{x_1}=\cT_{x_2}$.   \item $\chi_{\pi}(x_1)=\chi_{\pi}(x_2)$.     
    \end{enumerate}
\end{cor}
\begin{proof}
\EitanB{
\begin{enumerate}
    \item \label{it:Tv}
Take a non-zero vector $v \in V.$ It is enough to show the equality
$\cT_{x_1}(gv)=\cT_{x_2}(gv)$ for any $g \in G$. 
Take $n_0 \in \Z_{>0}$ such that $v \in V^{G_{-n_0}}$. Observe that
then for any $k \in G_0$ we have $kv \in  V^{G_{-n_0}}$.

Now take $K$ corresponding to $n_0$ from \Cref{lem : T_x}. 
Let $g \in G.$ By Cartan decomposition  (see \Cref{not:settings}\eqref{it:cartan}) $G=G_0T^{+}G_0$ we can write $g=k_1mk_2$ with $k_1,k_2 \in G_0$ and $m \in T^{+}.$ Now, as $k_2v \in V^{G_{-n_0}}$, by \Cref{lem : T_x} for $i=1,2$ we have:
$$\cT_{x_{i}}(gv)=\cT_{x_{i}}(k_1m (k_2v))=k_1 \cdot f_{x_{i}}^{K} \cdot m \cdot k_2v$$

The assumption $f^{K}_{x_1}=f^{K}_{x_2}$ now implies 
$$\cT_{x_{1}}(gv)=k_1 \cdot f_{x_{1}}^{K} \cdot m \cdot k_2v=k_1 \cdot f_{x_{2}}^{K} \cdot m \cdot k_2v=\cT_{x_{2}}(gv)$$

\item Follows from \eqref{it:Tv} in view of \Cref{thm:HC.char}.
\end{enumerate}
}
\end{proof}

\begin{lemma}\label{lem:fiber.const}
    Let $\bfX,\bfY,\bfZ$ be smooth algebraic varieties, equipped with simple $\mu$-rectifications. Let $\psi:\bfX\times\bfY\to\bfZ\times\bfY$ be a submersion s.t. the following diagram commutes
    $$\begin{tikzcd}
\bfX\times \bfY \arrow[rr, "\psi"] \arrow[dr, "p"'] & & \bfZ\times \bfY \arrow[dl, "q"] \\
& \bfY &
\end{tikzcd}$$
    where $p,q$ are the projections. For $y\in Y:=\bfY(F)$ let $\psi_y:\bfX\to\bfZ$ be s.t. $\psi(x,y)=(\psi_y(x),y)$; note that $\psi_y$ is again a submersion.

    Let $N\in\Z_{>0}$ and $f\in C_c^\infty(B_N^{\bfX})$. Then there is $c\in\Z_{>0}$ s.t. for any $y\in Y$ and any $y'\in B_{-c(ov_Y(y)+1)}^{\bfY}(y)$ we have
    $$(\psi_y)_*(f\mu_N^{\bfX})=(\psi_{y'})_*(f\mu_N^{\bfX}).$$
\end{lemma}
\RamiA{
\begin{proof}
$ $
\begin{enumerate}[{Step} 1.]
    \item Setup.\\    
    We equip $\bfX\times\bfY$ and $\bfZ\times\bfY$  \RamiD{with the simple $\mu$-rectifications induced from the $\mu$-rectifications on $\bfX$, $\bfY$ and $\bfZ$}. For all $m,s\in\Z$ and \RamiD{all $(x ,y)\in B_\infty^{\bfX\times \bfY}$ we have:}    \begin{equation}\label{eq:prod.rect}
    \mu_m^{\bfX\times\bfY}=\mu_m^{\bfX}\boxtimes\mu_m^{\bfY},\qquad
    B_{-s}^{\bfX\times\bfY}(x,y)=B_{-s}^{\bfX}(x)\times B_{-s}^{\bfY}(y),
    \end{equation}
    and likewise for $\bfZ\times\bfY$. \RamiD{Thus},  $\phi\boxtimes\chi$ is $\max(a,b)$-smooth whenever $\phi$ is $a$-smooth and $\chi$ is $b$-smooth.

    \RamiA{Let}  $m_0$ \RamiA{be s.t. $f$ is $m_0$-smooth. Let $\RamiA{c_{sub,\psi}}\RamiA{\in\Z_{>0}}$ be the constant from  \Cref{lem:push.mes.sm}}. Let $c_R\in\Z_{>0}$ be the constant of \Cref{rem:ball.vs.norm} for $\bfY$, and put $c:=c_{sub,\psi}c_R(m_0+N)$. Fix $y\in Y$ and set
    $$k:=\max\bigl(m_0,N,\lceil ov_Y(y)\rceil\bigr),\qquad m:=c_Rk .$$
    Then $m\geq k\geq\max(m_0,N)$,
    and $y\in Y_k\subset B^{\bfY}_{c_Rk}=B_m^{\bfY}$.
\item Application of  \Cref{lem:push.mes.sm} for $\psi$.\\
    Let $\RamiD{\tilde f}:=f\boxtimes 1_{B_m^{\bfY}}$.
    Since $1_{B_m^{\bfY}}$ is $m$-smooth (as $m\geq 0$), ${\RamiD{\tilde f}}$ is $\max(m_0,m)=m$-smooth.

     As $\operatorname{supp}f\subseteq B_N^{\bfX}$ and $\mu_N^{\bfX},\mu_m^{\bfX}$ have equal densities on $B_N^{\bfX}$ (recall $m\geq N$), we get $f\mu_N^{\bfX}=f\mu_m^{\bfX}$; hence by \eqref{eq:prod.rect}
    $$(f\mu_N^{\bfX})\boxtimes\mu_m^{\bfY}={\RamiD{\tilde f}}\,\mu_m^{\bfX\times\bfY}.$$
    Applying \Cref{lem:push.mes.sm} to $\psi$ and ${\RamiD{\tilde f}}$ yields an $m\RamiA{c_{sub,\psi}}$-smooth $g\in C_c^\infty(B_{m\RamiA{c_{sub,\psi}}}^{\bfZ\times\bfY})$ with
    \begin{equation}\label{eq:fiber.push}
    \psi_*\bigl(f\mu_N^{\bfX}\boxtimes\mu_m^{\bfY}\bigr)=g\,\mu_{m\RamiA{c_{sub,\psi}}}^{\bfZ\times\bfY}.
    \end{equation}
\item Expressing $(\psi_{y'})_*(f\mu_N^{\bfX})$ in terms of $g$.\\
    For $\varphi\in C_c^\infty(\bfZ(F))$ and $y'\in B_m^{\bfY}$ put
    \begin{align*}    h_\varphi(y')&:=\bigl\langle(\psi_{y'})_*(f\mu_N^{\bfX}),\varphi\bigr\rangle
    =\int_{\bfX}\varphi(\psi_{y'}(x))\,f(x)\,\mu_N^{\bfX}(x),\\
    g_\varphi(y')&:=\langle g(\cdot,y')\,\mu_{m\RamiA{c_{sub,\psi}}}^{\bfZ} ,\varphi \rangle =\int_{\bfZ}\varphi\,g(\cdot,y')\,\mu_{m\RamiA{c_{sub,\psi}}}^{\bfZ}.
    \end{align*}
We will show that $h_\varphi(y')=g_\varphi(y')$.
    
    As $\varphi\circ\psi$ is locally constant and $\operatorname{supp}f$ is compact, $h_\varphi$ is locally constant on $Y$; by \eqref{eq:prod.rect} and the $m\RamiA{c_{sub,\psi}}$-smoothness of $g$, the function $g_\varphi$ is $m\RamiA{c_{sub,\psi}}$-smooth. Since $\psi$ preserves the $\bfY$-coordinate, for every $\eta\in C_c^\infty(B_m^{\bfY})$, using $\mu_m^{\bfY}=\mu_{m\RamiA{c_{sub,\psi}}}^{\bfY}$ on $B_m^{\bfY}$ (as $m\RamiA{c_{sub,\psi}}\geq m$),
    \begin{align*}
    \int_{\bfY} h_\varphi\,\eta\,\mu_m^{\bfY}
    &=\int_{\bfZ\times\bfY}(\varphi\boxtimes\eta)\,\psi_*\bigl(f\mu_N^{\bfX}\boxtimes\mu_m^{\bfY}\bigr)\\
    &\overset{\eqref{eq:fiber.push}}{=}\int_{\bfY} g_\varphi\,\eta\,\mu_{m\RamiA{c_{sub,\psi}}}^{\bfY}
    =\int_{\bfY} g_\varphi\,\eta\,\mu_m^{\bfY}.
    \end{align*}
    The density of $\mu_m^{\bfY}$ is strictly positive on $B_m^{\bfY}$, so $h_\varphi=g_\varphi$ almost everywhere on $B_m^{\bfY}$; as both are locally constant, $h_\varphi=g_\varphi$ on all of $B_m^{\bfY}$.
\item Completion of the proof.\\
Let $y'\in B_{-c(ov_Y(y)+1)}^{\bfY}(y)$. We have      \begin{equation*}\label{eq:m.bound}
     \begin{split}
     m{c_{sub,\psi}}&=k c_R {c_{sub,\psi}}=
     \max\bigl(m_0,N,\lceil ov_Y(y)\rceil\bigr) c_R {c_{sub,\psi}}\\
     &\leq c_{sub,\psi}c_R(m_0+N)\,(ov_Y(y)+1) =  c\,(ov_Y(y)+1).
     \end{split}
    \end{equation*}
Thus $B_{-c(ov_Y(y)+1)}^{\bfY}(y) \subset  B_{-m\RamiA{c_{sub,\psi}}}^{\bfY}(y) \subseteq B_m^{\bfY}$. Hence $y'\in B_m^{\bfY}$, and since $g_\varphi$ is $m\RamiA{c_{sub,\psi}}$-smooth,
    $$h_\varphi(y')=g_\varphi(y')=g_\varphi(y)=h_\varphi(y).$$
    As $\varphi\in C_c^\infty(\bfZ(F))$ was arbitrary, $(\psi_{y'})_*(f\mu_N^{\bfX})=(\psi_y)_*(f\mu_N^{\bfX})$.
\end{enumerate}
\end{proof}
}

\begin{proof}[Proof of \Cref{thm:smooth}]
\RamiB{Recall that $G^{rss} \subset G$ is an open subset.}    
    By \Cref{thm:HC.char}, $\chi_\pi|_{G^{rss}}$ is represented by the locally constant function $\theta_\pi(x):=\operatorname{tr}(\cT_x)$. \RamiB{So, it suffices} to find $c_{reg}(\pi)\in\Z_{>0}$ s.t. for every $x\in G^{rss}$ the coset $xG_{-c_{reg}(\pi)\cdot(ov_{G^{rss}}(x)+1)}$ is contained in $G^{rss}$ and the function $\theta_\pi$ is constant on it.

    Fix the open compact $K<B$ of \Cref{cor:fK.det}. Consider the simple $\mu$-rectified varieties $\bfX:=\bfG\times\bfB$, $\bfZ:=\bfG$ and $\bfY:=\bfG^{rss}$, with top forms giving Haar measures, and the map
    $$\psi:\bfX\times\bfY\to\bfZ\times\bfY,\qquad\psi\bigl((h,b),x\bigr)=(hxh^{-1}b,\,x).$$
    The map $\psi$ is a submersion by \Cref{lem:HC.submersion}. Put $f:=1_{G_0\times K}\in C_c^\infty(B_N^{\bfX})$ with $N$ large enough so that $G_0\times K\subseteq B_N^{\bfX}$; on $\operatorname{supp}f$ the Haar measure equals $\mu_N^{\bfX}$, so $f^K_x=(\psi_x)_*(f\mu_N^{\bfX})$ for every $x\in G^{rss}$.

    Let $c$ be the constant of \Cref{lem:fiber.const} for the above data: for every $n\in\Z_{>0}$ and every $x\in B^{\bfG^{rss}}_n$, $f_x^K$ is constant on $B^{\bfG^{rss}}_{-c(n+1)}(x)$. By \Cref{prop:G.norms}\eqref{prop:G.norms:2}, $ov_{G^{rss}}$ is the $ov$ of a norm on $\bfG^{rss}(F)$, so $(G^{rss})_n=\{ov_{G^{rss}}\le n\}$; hence \Cref{rem:ball.vs.norm} gives $c_R\in\Z_{>0}$ with $(G^{rss})_n\subseteq B^{\bfG^{rss}}_{c_R n}$. Let $c_{ball}$ be the constant of \Cref{lem:coset.ball} for the above rectification of $\bfG^{rss}$ and set $$c_{reg}(\pi):=c_{ball}\bigl(c(c_R+1)+1\bigr).$$

    Fix $x\in G^{rss}$, put $r:=ov_{G^{rss}}(x)$, $n:=c_R(r+1)$ and $m:=c(c_R+1)(r+1)$. Then $x\in(G^{rss})_{r+1}\subseteq B^{\bfG^{rss}}_{c_R(r+1)}= B^{\bfG^{rss}}_{n}$, so $f_x^K$ is constant on $B^{\bfG^{rss}}_{-c(n+1)}(x)$, and $m\geq c(n+1)$, so it is constant on $B^{\bfG^{rss}}_{-m}(x)$. By \Cref{lem:coset.ball},
    $$xG_{-c_{reg}(\pi)(r+1)}\subset xG_{-c_{ball}(m+r)}\subset B^{\bfG^{rss}}_{-m}(x)\subset G^{rss} ,$$
    the first inclusion because
    \begin{align*}
    c_{reg}(\pi)(r+1)&=c_{ball}\bigl(c(c_R+1)+1\bigr)(r+1)=
    c_{ball}(c(c_R+1)(r+1)+r+1)\\
    &=c_{ball}(m+r+1)\geq c_{ball}(m+r).
    \end{align*}
    Hence $f^K_{x'}=f^K_x$ for every $x'\in xG_{-c_{reg}(\pi)(r+1)}$, and by \Cref{cor:fK.det}, $\theta_\pi(x')=\theta_\pi(x)$ for every such $x'$.
\end{proof}
\appendix
\section{Supporting statements}\label{sec:app}

\subsection{The function $\Theta$ on the dominant cone}\label{ssec:app.Theta}
The function \EitanD{ $\Theta: \Lambda(T) \to \Q$} of \Cref{not:settings}\eqref{it:N.Theta} is $\fN$ times the sum of the
positive roots. The lemma below records few consequences \EitanD{that are used in the proof of Lemma  \ref{lem:Theta.ZG.bound} in \Cref{ssec:app.ZG}}: $\Theta$ vanishes
on the directions on which every root vanishes, and on the dominant cone it dominates every
positive root.

\begin{lemma}\label{lem:Theta.dom}
Denote $\mdef{V}:=\Lambda(T)\otimes_\Z\R$ and
$\mdef{V_0}:=\bigcap_{\beta\in\Delta}\ker\beta\subset V$. Then
\begin{enumerate}
\item\label{lem:Theta.dom:1} $V=\operatorname{span}_\R(\Delta^\vee)\oplus V_0$, where $\mdef{\Delta^\vee}\subset V$ is
the set of {coroots} 
of the elements of $\Delta$ \RamiC{(see \cite[V.3.9]{Ren})}.
\item \label{lem:Theta.dom:2.1} $\Theta$ vanishes on $V_0$.
\item\label{lem:Theta.dom:2}
$$0\leq \beta(\lambda)\leq \Theta(\lambda)\qquad\text{for every }
\lambda\in\Lambda^+(T)\text{ and every }\beta\in\Phi_\bfB .$$
\end{enumerate}
\end{lemma}

\begin{proof}
By \Cref{not:settings}\eqref{it:Phi} the elements of $\Phi$ are $\R$-valued functionals on $V$.
The elements of $\Delta$ are linearly independent, being a base of the root system $\Phi$, and the
matrix $(\beta(\alpha^\vee))_{\alpha,\beta\in\Delta}$ is invertible. This gives
\eqref{lem:Theta.dom:1}.

Item 
\eqref{lem:Theta.dom:2.1} is clear: Every $\beta\in\Phi_\bfB$ is a non-negative integral combination of $\Delta$, so it vanishes on
$V_0$. Hence so does $\Theta=\fN\sum_{\beta\in\Phi_\bfB}\beta$.

Let $\lambda\in\Lambda^+(T)$. By \Cref{not:settings}\eqref{it:M+} and \eqref{it:Lam+} we have
$\gamma(\lambda)\geq0$ for every $\gamma\in\Phi_\bfB$, which is the left inequality of
\eqref{lem:Theta.dom:2}. All the summands of $\Theta(\lambda)$ being non-negative, for every
$\beta\in\Phi_\bfB$
$$
\Theta(\lambda)=\fN\sum_{\gamma\in\Phi_\bfB}\gamma(\lambda)\geq\fN\beta(\lambda)\geq\beta(\lambda),
$$
the last step because $\fN\geq1$. This is the right inequality of \eqref{lem:Theta.dom:2}.
\end{proof}

\subsection{The Cartan parameter and the two size functions}\label{ssec:app.norms}
Two functions measure the size of an element of $G$: the intrinsic $ov_G$ of
\Cref{not:settings}\eqref{it:ov3}, and the valuation $ov_\sigma$ of the matrix entries in the
fixed embedding of \Cref{not:sigma}. This subsection compares both of them with the norm of the
Cartan parameter, and deduces \Cref{cor:ov.sigma.G}, which compares them with each other.

\begin{notation}\label{not:cartan.lam}
For $g\in G$ denote by $\mdef{\lambda(g)}\in\Lambda^+(T)$ the element corresponding to the double coset
$G_0gG_0$ under \Cref{not:settings}\eqref{it:cartan}. Fix a norm $\mdef[\vert\vert\cdot\vert\vert]{\|\cdot\|}$ on the real vector
space $V=\Lambda(T)\otimes_\Z\R$, \RamiC{that is invariant with respect to the action of the Weyl group of $T$.} 
\end{notation}

\begin{lemma}\label{lem:sigma.basic}
$ $
\begin{enumerate}
\item\label{lem:sigma.basic:1} $ov_\sigma(k_1gk_2)=ov_\sigma(g)$ for all $g\in G$ and
$k_1,k_2\in G_0$, and $ov_\sigma(g^{-1})=ov_\sigma(g)$.
\item\label{lem:sigma.basic:2} $ov_\sigma(gh)\leq ov_\sigma(g)+ov_\sigma(h)$ for all $g,h\in G$.
\end{enumerate}
\end{lemma}

\begin{proof}
\eqref{lem:sigma.basic:1} By the choice of $\cL$ we have $\sigma(k)_{ij}\in\cO_F$ and
$(\sigma(k)^{-1})_{ij}\in\cO_F$ for $k\in G_0$. Since $|\cdot|$ is ultrametric, multiplying by such
a matrix does not increase $\|\cdot\|$, so $ov_\sigma(k_1gk_2)\leq ov_\sigma(g)$; applying this to
$g=k_1^{-1}(k_1gk_2)k_2^{-1}$ gives equality. The second assertion is immediate from the
definition.

\eqref{lem:sigma.basic:2} $|(\sigma(gh))_{ij}|=|\sum_k\sigma(g)_{ik}\sigma(h)_{kj}|\leq
\ell^{ov_\sigma(g)+ov_\sigma(h)}$, and likewise for $\sigma((gh)^{-1})=\sigma(h)^{-1}\sigma(g)^{-1}$.
\end{proof}

\begin{lemma}\label{lem:ov.cartan}
Both $ov_G$ and $ov_\sigma$ are invariant under left and right multiplication by $G_0$. \RamiC{I}n
particular each of them \RamiC{factors through $\lambda:G\to \Lambda^+(T)$}.
\end{lemma}

\begin{proof}
For $ov_\sigma$ this is \Cref{lem:sigma.basic}\eqref{lem:sigma.basic:1}. \RamiC{Recall that
$ov_G(g)=\max(ov_{G^{ab}}(g^{ab}),ov_{G^{ad}}(g^{ad}))$ where $g^{ab}$ and $g^{ad}$ are the images of $g \in G$ in $G^{ab}$ and $G^{ad}$.}
The function
$ov_{G^{ad}}$ is by construction a function on $G_0^{ad}\backslash G^{ad}/G_0^{ad}$
(\Cref{not:settings}\eqref{it:ov1}), and the image of $G_0$ in $G^{ab}$ is a compact subgroup,
hence $|\chi|=1$ on it for every $F$-rational character $\chi$ of $\bfG^{ab}$, so that this image
lies in ${}^0G^{ab}$ and $ov_{G^{ab}}$ is unchanged. This shows the invariance of $ov_G$.

The last assertion follows from \RamiC{Cartan decomposition formulated in}
\Cref{not:settings}\eqref{it:cartan}.
\end{proof}

\begin{lemma}\label{lem:ov.G.comp}
There is $c\in\Z_{>0}$ such that $c^{-1}\|\lambda(g)\|\leq ov_G(g)\leq c\|\lambda(g)\|$ for every
$g\in G$.
\end{lemma}

\begin{proof}
By \Cref{lem:ov.cartan} we may assume $g\in T^+$; 
let
$\RamiC{\lambda(g)}^{ab}\in\Lambda(G^{ab})$ be \RamiC{the image of ${\lambda(g)}$}  under the map induced by $\bfT\to\bfG^{ab}$. Then
\begin{equation}\label{eq:ovG.cartan}
ov_G(g)=\max\bigl(ov_{\Lambda(G^{ab})}(\RamiC{\lambda(g)}^{ab}),\Theta(\RamiC{\lambda(g)})\bigr)
\end{equation}
by \Cref{not:settings}\eqref{it:ov1}--\eqref{it:ov3}. Both
$v\mapsto ov_{\Lambda(G^{ab})}(v^{ab})$ and $\Theta$ are restrictions to $\Lambda(T)$ of continuous
functions on $V$ that are positively homogeneous of degree $1$; taking for $c$ an upper bound of
their maximum on $\{\|v\|=1\}$ gives $ov_G(g)\leq c\|\RamiC{\lambda(g)}\|$.

For the lower bound consider the linear map
$$
p:V\to \Lambda(G^{ab})\otimes_\Z\R\times\R^{\Phi_\bfB},\qquad
p(v)=\bigl(v^{ab},(\beta(v))_{\beta\in\Phi_\bfB}\bigr).
$$
We claim that $p$ is injective. Indeed, let $v\in\ker p$. Then $v\in V_0$ in the notation of
\Cref{lem:Theta.dom}. Identify $V$ with $X_*(\bfA_\bfT)\otimes_\Z\R$ using \cite[V.2.6]{Ren}; under
this identification $V_0=X_*(\bfA_\bfG)\otimes_\Z\R$, since $\bfA_\bfG$ is the identity component of
the intersection of the kernels of the elements of $\Phi$, and $\Phi=\Phi_\bfB\sqcup(-\Phi_\bfB)$
by \Cref{not:settings}\eqref{it:Phi.P}. The group $\bfA_\bfG\cap[\bfG,\bfG]$ is finite, so
$\bfA_\bfG\to\bfG^{ab}$ has finite kernel and its image is a split subtorus of $\bfG^{ab}$;
consequently $X_*(\bfA_\bfG)\otimes\R\to\Lambda(G^{ab})\otimes\R$ is injective. As $v^{ab}=0$ we
get $v=0$, proving the claim.

Since $p$ is injective, $v\mapsto\max\bigl(ov_{\Lambda(G^{ab})}(v^{ab}),\max_\beta|\beta(v)|\bigr)$
is a norm on $V$, hence is bounded below by $c_1^{-1}\|\cdot\|$ for some $c_1$. For
$\lambda\in\Lambda^+(T)$ we have $0\leq\beta(\lambda)\leq \Theta(\lambda)$ by
\Cref{lem:Theta.dom}\eqref{lem:Theta.dom:2}, so
$$
c_1^{-1}\|\lambda\RamiC{(g)}\|\leq
\max\bigl(ov_{\Lambda(G^{ab})}(\lambda^{ab} \RamiC{(g)}),\Theta(\lambda\RamiC{(g)})\bigr)
\overset{\eqref{eq:ovG.cartan}}{=} ov_G(g),
$$
which gives the assertion after enlarging $c$.
\end{proof}

\begin{lemma}\label{lem:ov.sigma.comp}
There is $c\in\Z_{>0}$ such that
$$c^{-1}\|\lambda(g)\|-c\leq ov_\sigma(g)\leq c(\|\lambda(g)\|+1)$$ for every $g\in G$.
\end{lemma}

\begin{proof}
$ $
\begin{enumerate}[{Case} 1.]
    \item \RamiC{$g\in A_\bfT$}.\\
     Let $\chi_1,\dots,\chi_r$ be a $\Z$-basis of $X^*(\bfA_\bfT)$. By
\cite[Example 18.2]{Kot} the map
$a\mapsto\max_i\max(|\chi_i(a)|,|\chi_i(a)|^{-1})=\ell^{\max_i|\operatorname{val}\chi_i(a)|}$ is a
norm on $A_\bfT$, and so is the restriction of $\ell^{ov_\sigma}$, by \Cref{lem:norm.sigma},
\Cref{lem:nrm.cl.emb} and the fact that $\bfA_\bfT$ is a closed subgroup of $\bfG$. Being two
norms on the same set they are equivalent. Denote by $\mu(a)\in\Lambda(A_\bfT)\subset\Lambda(T)$
the class of $a\in A_\bfT$. By \cite[V.2.6]{Ren} the valuation identifies $\Lambda(A_\bfT)$ with
$X_*(\bfA_\bfT)$ and $\Lambda(A_\bfT)$ is of finite index in $\Lambda(T)$, so
$\Lambda(A_\bfT)\otimes_\Z\R=V$ and $\max_i|\langle\chi_i,\cdot\rangle|$ is a norm on $V$;
being two norms on a finite-dimensional real vector space, it and $\|\cdot\|$ are bounded by
constant multiples of each other. Together, there is $c_1\in\Z_{>0}$ \RamiC{such that for any $a\in A_\bfT$ we have:}
\begin{equation}\label{eq:ovsigma.split}
\begin{split}
\RamiC{c_1^{-1}\|\lambda(a)\|-c_1}&=c_1^{-1}\|\mu(a)\|-c_1\\
&\leq ov_\sigma(a)\\
&\leq c_1(\|\mu(a)\|+1)=\RamiC{c_1(\|\lambda(a)\|+1)}.
\end{split}
\end{equation}
\item \RamiC{$g\in T$}.\\
 Choose $t_1,\dots,t_h\in T$ whose classes represent
$\Lambda(T)/\Lambda(A_\bfT)$. Since ${}^0T$ is compact by \Cref{not:settings}\eqref{it:0H}, the
set $\RamiC{C_\bfT}:={}^0T\{t_1,\dots,t_h\}$ is compact and $T=A_\bfT {\RamiC{C_\bfT}}$. Let $c_0$ be a bound for $ov_\sigma$
on ${\RamiC{C_\bfT}}\cup {\RamiC{C_\bfT}}^{-1}$ and $c_2$ a bound for $\|\nu\|$, $\nu$ running over the classes of the elements
of ${\RamiC{C_\bfT}}$.
\RamiC{Put $c_3:=\lceil c_1 c_2+c_1+c_0\rceil$.}

\RamiC{For any $g\in T$ } write $g=ab$ with $a\in A_\bfT$ and $b\in {\RamiC{C_\bfT}}$. By
\Cref{lem:sigma.basic}\eqref{lem:sigma.basic:2},
$$
ov_\sigma(a)-c_0\leq ov_\sigma(g)\leq ov_\sigma(a)+c_0.
$$
\RamiC{
Thus
\begin{align*}
ov_\sigma(g)&\leq ov_\sigma(a)+c_0\\
&\overset{\eqref{eq:ovsigma.split}}{\leq} c_1\|\lambda(a)\|+c_1+c_0\\
&\leq c_1(\|\lambda(g)\|+c_2)+c_1+c_0\\
&\leq c_3\|\lambda(g)\|+c_3 .
\end{align*}
Similarly
\begin{align*}
c_3^{-1}\|\lambda(g)\|-c_3&\leq c_1^{-1}(\|\lambda(g)\|-c_2)-c_1-c_0\\
&\leq c_1^{-1}\|\lambda(a)\|-c_1-c_0\\
&\overset{\eqref{eq:ovsigma.split}}{\leq} ov_\sigma(a)-c_0\\
&\leq ov_\sigma(g) .
\end{align*}
}
\item \RamiC{General case}.\\
By \Cref{lem:ov.cartan} and \Cref{not:settings}\eqref{it:cartan} both $ov_\sigma(g)$ and
$\|\lambda(g)\|$ are unchanged if $g$ is replaced by the element of $T^+$ with the same class
$\lambda(g)$. \RamiC{So, the statement follows from the previous case.}
\end{enumerate}

\end{proof}

\begin{proof}[Proof of \Cref{cor:ov.sigma.G}]
Combine \Cref{lem:ov.G.comp} and \Cref{lem:ov.sigma.comp}.
\end{proof}

\begin{remark}\label{rem:app.Levi}
Let $M<G$ be a standard Levi subgroup, with the data fixed for it in \Cref{not:settings}, and let
$\sigma_M$ be a closed embedding of $\bfM$ and $\cL_M$ a lattice as in \Cref{not:sigma} for $M$.
Nothing in \S\ref{ssec:app.Theta}--\S\ref{ssec:app.norms} uses more about $\bfG$ than
\Cref{not:settings} and \Cref{not:sigma} provide, so every statement of these two subsections
holds for $M$; in particular so does \Cref{cor:ov.sigma.G}. The same applies to
\Cref{lem:norm.sigma}, whose proof uses only that $\sigma$ is a closed embedding, and hence to
\Cref{prop:G.norms}\eqref{prop:G.norms:1}: the map $m\mapsto \ell^{ov_M(m)}$ is a norm on $M$.
\end{remark}

\begin{notation}\label{not:N1}
Denote by $\mdef{N_1}:G\to\R_{\geq0}$ the function $N_{x_0}$ of \cite[\S19.2]{Kot}, for
$\bfX=\bfG$ and $x_0=1$, attached to the generators
$\sigma_{ij}-\delta_{ij},(\sigma^{-1})_{ij}-\delta_{ij}$ of the $F$-algebra $F[\bfG]$ --- they
generate it because $\sigma$ is a closed embedding, and they vanish at $1$, which is the
condition required there. Thus
$$
N_1(g)=\max_{i,j}\max\bigl(|\sigma(g)_{ij}-\delta_{ij}|,|(\sigma(g)^{-1})_{ij}-\delta_{ij}|\bigr),
$$
and by \cite[(19.2.1)]{Kot} the sets $\{g\in G: N_1(g)<\varepsilon\}$, $\varepsilon>0$, form a
basis of neighborhoods of $1$ in $G$. For $s\in\Z_{>0}$ we have
\begin{equation}\label{eq:N1.cong}
\sigma(g)\in 1+\varpi^{s}\End_{\cO_F}(\cL)\ \Longrightarrow\ N_1(g)\leq \ell^{-s},
\end{equation}
because $1+\varpi^{s}\End_{\cO_F}(\cL)$ is a subgroup of $\GL(\cW)$: for
$X\in\End_{\cO_F}(\cL)$ one has $(1+\varpi^sX)^{-1}=\sum_{i\geq0}(-\varpi^sX)^i\in
1+\varpi^{s}\End_{\cO_F}(\cL)$.
\end{notation}
\subsection{Bounding an element by $\Theta$ modulo the center}\label{ssec:app.ZG}

\begin{proof}[Proof of \Cref{lem:Theta.ZG.bound}]
$ $
\begin{enumerate}[Step 1.]
    \item \RamiC{There is $c_1\in \Z_{>0}$ s.t. for any $a\in T^+$ there is $\nu\in \Lambda(\cZ(G))\subset \Lambda^+(T)$ with
    $\|\lambda(a)-\nu\|\leq c_1(\Theta(a)+1).$}\\
\RamiC{Let $\|\cdot\|$ be the norm on
$V=\Lambda(T)\otimes_\Z\R$ fixed in \Cref{not:cartan.lam}, and let $V_0$ be as in
\Cref{lem:Theta.dom}.}
By \Cref{lem:Theta.dom}\eqref{lem:Theta.dom:1}, for $a\in T^+$ write $\lambda(a)=\lambda_1\RamiC{(a)}+\lambda_0\RamiC{(a)}$ with
$\lambda_1\RamiC{(a)}\in\operatorname{span}_\R(\Delta^\vee)$ and $\lambda_0\RamiC{(a)}\in V_0$. 

The elements of
$\Phi_\bfB$ vanish on $V_0$ and separate the points of
$\operatorname{span}_\R(\Delta^\vee)$, so
$\|\lambda_1\RamiC{(a)}\|\leq c'\max_{\beta\in\Phi_\bfB}|\beta(\lambda(a))|$ for a constant $c'$.
\RamiC{By} \Cref{lem:Theta.dom}\eqref{lem:Theta.dom:2},
we have $0\leq\beta(\lambda(a))\leq \RamiC{\Theta(a)}$ for every $\beta$, whence
$\|\lambda(a)-\lambda_0\RamiC{(a)}\|=\|\lambda_1\RamiC{(a)}\|\leq c'\RamiC{\Theta(a)}$. 

The image
of $\Lambda(\cZ(G))$ in $V$ contains the image of $\Lambda(\bfA_\bfG)$, which is a lattice
spanning $V_0$ by \cite[V.2.6]{Ren}, hence there is a
constant $c''$, its covering radius, and $\nu\in\Lambda(\cZ(G))$ with
$\|\lambda_0-\nu\|\leq c''$. \RamiC{Taking $c_1:=\lceil\max(c',c'')\rceil$ completes the step.}
\item \RamiC{Completion of the proof}.\\ 
\RamiC{Let $c_2$ be the constant $c$ from \Cref{lem:ov.G.comp}.
Put $c_Z:=c_1c_2+c_2$.
Let $g=k_1ak_2$ and $R$ be as in the formulation of the lemma.  Let $\nu$ be as in the previous step applied to $a$. 
}
Choose $z\in \cZ(G)$ whose image in $\Lambda^+(T)$ is $\nu$.
The class of $z^{-1}a$ in $\Lambda(T)$ is $\lambda(a)-\nu$, and
$\lambda(z^{-1}a)$ is its image under one of the  elements of the Weyl group, so $$\|\lambda(z^{-1}a)\|=\|\lambda(a)-\nu\|.$$

Now we have 

\RamiC{
\begin{align*}
ov_G(z^{-1}g)&\overset{\text{Lem \ref{lem:ov.cartan}}}{=}ov_G(z^{-1}a)\\
&\overset{\text{Lem \ref{lem:ov.G.comp}}}{\leq}c_2\|\lambda(z^{-1}a)\|+c_2\\
&=c_2\|\lambda(a)-\nu\|+c_2\\
&\leq c_2c_1(\Theta(a)+1)+c_2\\
&\leq c_Z(\Theta(a)+1).
\end{align*}
}
Hence $g\in \cZ(G)G_{c_Z(R+1)}$.
\end{enumerate}
\end{proof}

\subsection{The root subgroups and their filtration}\label{ssec:app.root}
\RamiC{The purpose of this subsection is to prove \Cref{lem:U.contract}. The proof is based on analyzing the action of elements of $T$ on unipotent radicals of parabolic subgroups. This is done using a decomposition of these unipotent radicals into a product of root subgroups and analyzing the action of $T$ on these root subgroups via a certain filtration of these groups.}
\begin{notation}\label{not:root.sub}
$ $
\begin{enumerate}
\item\label{not:root.grp}
For every $\beta\in \Phi$ let $\mdef{\bfU_\beta}$ be the corresponding \tdef{root subgroup}, i.e. the subgroup
whose Lie algebra is the sum of the weight spaces of $\bfA_\bfT$ in $\Lie(\bfG)$ whose weights
are the positive integer multiples of $\beta$. See \cite[V.3.7]{Ren}.
\RamiC{Set $\mdef{U_\beta}:=\bfU_\beta(F)$.}
\item\label{not:root.filtration}
For $\beta\in\Phi$ and $r\in\R$ let
$$
\mdef{U_{\beta,r}}:=U_{\beta,\fx,r}\subset U_\beta
$$
be the Bruhat--Tits root-filtration subgroup at the fixed special vertex
$\fx$, i.e. $U_{\beta,r}=\phi_\beta^{-1}([r,\infty])$ for the valuation $\phi$ of the root
datum of $(\bfG,\bfA_\bfT)$ corresponding to $\fx$ \RamiC{(see \cite[Construction 9.6.1]{KP})}. We use the normalization of
\cite[Definition 6.1.2]{KP}. The
family $(U_{\beta,r})_{r\in\R}$ is decreasing in $r$.
\item\label{not:root.div}
\RamiF{Let $\beta\in\Phi$ be such that $2\beta$ is a root of $(\bfG,\bfA_\bfT)$, let
$\mdef{U_{2\beta}}<U_\beta$ be the root subgroup of $2\beta$, and for $r\in\R$ let
$$
\mdef{U_{2\beta,r}}:=U_{2\beta,\fx,r}=\phi_{2\beta}^{-1}([r,\infty])\subset U_{2\beta}
$$
be the corresponding root-filtration subgroup at $\fx$. 
If $2\beta$ is not a root we set $U_{2\beta}:=\{1\}$ and $U_{2\beta,r}:=\{1\}$, as in
\cite[Remark 7.3.4]{KP}. By axiom V\,4 of \cite[Definition 6.1.2]{KP} we have
$\phi_{2\beta}=2\phi_\beta$ on $U_{2\beta}$, and consequently
$$
U_{2\beta,2r}=U_{2\beta}\cap U_{\beta,r}\qquad(r\in\R).
$$}
\end{enumerate}
\end{notation}

\begin{lemma}[{\cite[Definition 6.1.2, Axiom V\,6]{KP}}]\label{lem:root.conj}
\RamiC{
Let $t\in T,\ \beta\in\Phi,\ r\in\R$. Then
we have: }
$$
tU_{\beta,r}t^{-1}=U_{\beta,\,r+\beta(t)}
$$
\RamiF{
and
$$
tU_{2\beta,r}t^{-1}=U_{2\beta,\,r+2\beta(t)}
$$
}
where $\beta(t)$ is the value at the class of $t$ in $\Lambda(T)$ of $\beta$ viewed as a
$\Q$-valued functional as in \Cref{not:settings}\eqref{it:Phi}.
\RamiC{In particular, $T$ normalizes each $U_\beta$.} 
\end{lemma}

\begin{lemma}\label{lem:sigma.nbhd}
There is $\mdef{c_{nbh}}\in\Z_{>0}$ such that
$$
\bigl\{g\in G\ :\ \sigma(g)\in 1+\varpi^{c_{nbh}}\End_{\cO_F}(\cL)\bigr\}\subset G_{-1},
$$
where $\mdef{\End_{\cO_F}(\cL)}\subset \End_F(\cW)$ is the ring of the endomorphisms preserving the
lattice $\cL$ of \Cref{not:sigma}.
\end{lemma}

\begin{proof}
By \Cref{not:N1} the sets $\{g\in G: N_1(g)<\varepsilon\}$, $\varepsilon>0$, form a basis of
neighborhoods of $1$ in $G$, and $G_{-1}$ is open, so $\{g\in G:N_1(g)<\varepsilon\}\subset
G_{-1}$ for some $\varepsilon>0$. Any $c_{nbh}\in\Z_{>0}$ with $\ell^{-c_{nbh}}<\varepsilon$
then satisfies the assertion, by \eqref{eq:N1.cong}.
\end{proof}

\RamiC{
\begin{lemma}\label{lem:U.filtration.alg}
Let \RamiC{$\bfQ<\bfG$} be a
standard parabolic or the opposite of one, let $\RamiC{\bfU_\bfQ}$ be its unipotent radical. \RamiC{Fix}
an ordering of $\Phi_Q$. Then, with the notation of \Cref{not:root.sub}
the multiplication map
\RamiC{$\prod_{\beta\in\Phi_Q}\bfU_\beta\to \bfU_\bfQ$}
is an isomorphism of $F$-varieties. 
\end{lemma}
}
\begin{proof}
    \RamiC{This follows from \cite[Theorem 2.11.15]{KP}, applied to the
action of the split torus $\bfA_\bfT$ on the smooth connected solvable group $\bfU_\bfQ$, and to
the partition of the set of weights of $\bfA_\bfT$ in $\Lie(\bfU_\bfQ)$ into the sets of its
elements that are positive integer multiples of $\beta$, for $\beta\in\Phi_Q$. Indeed, two
distinct elements of $\Phi_Q$ are not proportional, so these sets generate pairwise disjoint
subsemigroups of $X^*(\bfA_\bfT)$, and the subgroup that \cite[Proposition 2.11.13]{KP} attaches
to the subsemigroup generated by the set corresponding to $\beta$ is $\bfU_\beta$ by
\Cref{not:root.sub}\eqref{not:root.grp}.}
\end{proof}
\begin{lemma}\label{lem:U.level0}
Let $Q<G$ be a standard parabolic or the opposite of one, with standard Levi subgroup $M_Q$ and
unipotent radical $U_Q$, and fix an ordering of $\Phi_Q$. Then the multiplication map
$$
\prod_{\beta\in\Phi_Q}U_{\beta,0}\xrightarrow{\ \sim\ }U_Q\cap G_0
$$
\RamiC{is a bijection}.
\end{lemma}

\begin{proof}
\RamiB{WLOG assume that $Q$ is standard}. Every $\beta\in
\Phi_B\smallsetminus\Phi_Q$ is a root of $\bfM_Q$: indeed $\Phi$ is the disjoint union of
$\Phi_Q$, of $-\Phi_Q$ and of the set of its elements that are roots of $\bfM_Q$, and
$\Phi_B\cap(-\Phi_Q)=\varnothing$ because $\Phi_B\cap(-\Phi_B)=\varnothing$. Denote by
$N$ the normalizer of $\bfA_\bfT$ in $G$.

We apply \cite[Proposition 7.3.12]{KP} at the point $\fx$, to the constant function $0$ on
$\Phi\cup\{0\}$ --- which is concave in the sense of \cite[Definition 7.3.1]{KP} --- and to the
positive system of roots of $(\bfG,\bfA_\bfT)$ consisting of those occurring in
$\bfU_{\bfB^-}$, so that the sets called $\Phi^{+,nd}$ and $\Phi^{-,nd}$ there are
$\Phi_{B^-}$ and $\Phi_B$, and the groups called $U^+$ and $U^-$ there are $U_{\bfB^-}$
and $U_\bfB$. For this function, the group written in bold face as $\mathbf{U}_{\beta,\fx,0}$ in
\cite[Definition 7.3.3(1)]{KP} is $U_{\beta,0}$: it is by definition
$U_{\beta,0}U_{2\beta,\fx,0}$, where $U_{2\beta}<U_\beta$ is the root subgroup of $2\beta$ if
$2\beta$ is a root of $(\bfG,\bfA_\bfT)$ and is trivial otherwise, and
$U_{2\beta,\fx,r}:=\phi_{2\beta}^{-1}([r,\infty])$; and axiom V\,4 of
\cite[Definition 6.1.2]{KP} gives $\phi_{2\beta}=2\phi_\beta$ on $U_{2\beta}$, so
$U_{2\beta,\fx,0}=U_{2\beta}\cap U_{\beta,0}$. Consequently the group $G^\sharp_{\fx,0}$ of
\cite[Definition 7.3.3(2)]{KP} is the subgroup of $G$ generated by the groups $U_{\beta,0}$,
$\beta\in\Phi$, and $G_{\RamiC{\fx,}0}=G^\sharp_{\fx,0}Z$ by \cite[Definition 7.3.3(3)]{KP}, where $Z<T$ is
the group denoted $Z(k)_0$ there. Put $U_{\bfB,0}:=U_\bfB\cap G^\sharp_{\fx,0}$ and
$U_{\bfB^-,0}:=U_{\bfB^-}\cap G^\sharp_{\fx,0}$.
\begin{enumerate}[Step 1.]
    \item 
$N\cap U_\bfB U_{\bfB^-}=\{1\}$.\\
\DimaB{Let}
$n=vw$ with $n\in N$, $v\in U_\bfB$ and $w\in U_{\bfB^-}$. For $s\in A_\bfT$ put
$t_s:=(n^{-1}sn)s^{-1}\in A_\bfT$; then
$$
(svs^{-1})\cdot 1\cdot (sws^{-1})=sns^{-1}=nt_s=v\cdot t_s\cdot (t_s^{-1}wt_s)\RamiC{.}
$$
\DimaB{Since the multiplication map $U_\bfB\times T\times
U_{\bfB^-}\to G$ is injective, this implies
$t_s=1$}. Thus $n$ centralizes $A_\bfT$, i.e. $n\in T$
(\cite[V.4.1]{Ren}), and comparing $v\cdot 1\cdot w$ with $1\cdot n\cdot 1$ gives $n=1$.

\item $U_Q\cap G_{\RamiC{\fx,}0}\subset U_{\bfB,0}$.\\
By \cite[Proposition 7.3.12(1)]{KP} we have
$G^\sharp_{\fx,0}=U_{\bfB,0}\,U_{\bfB^-,0}\,(N\cap G^\sharp_{\fx,0})$, and $Z<N$, so
$G_{\RamiC{\fx,}0}=U_{\bfB,0}\,U_{\bfB^-,0}\,(N\cap G_{\RamiC{\fx,}0})$. Let $u\in U_Q\cap G_{\RamiC{\fx,}0}$ and write $u=vwn$
accordingly. As $u,v\in U_\bfB$ we get $n^{-1}=u^{-1}vw\in U_\bfB U_{\bfB^-}$, so $n=1$ by Step 1; then
$w=v^{-1}u\in U_{\bfB^-}\cap U_\bfB=\{1\}$, and $u=v\in U_{\bfB,0}$.

\item \RamiC{$
\prod_{\beta\in\Phi_Q}U_{\beta,0}\xrightarrow{\ \sim\ }U_Q\cap G_{{\fx,}0}
$
is a bijection}.\\
Order $\Phi_B$ so that the elements of $\Phi_Q$ come first, in the
fixed order. By \cite[Proposition 7.3.12(3)]{KP} the multiplication map
$\prod_{\beta\in\Phi_B}U_{\beta,0}\to U_{\bfB,0}$ is a bijection. Let $u\in U_Q\cap G_{\RamiC{\fx,}0}$. By
Step 2 it factors as $u=u_Qu_M$ with $u_Q\in\prod_{\beta\in\Phi_Q}U_{\beta,0}$ and
$u_M\in\prod_{\beta\in\Phi_B\smallsetminus\Phi_Q}U_{\beta,0}$. Then $u_Q\in U_Q$ and
$u_M\in M_Q$, so $u_M=u_Q^{-1}u\in U_Q\cap M_Q=\{1\}$ and $u=u_Q$. 
\DimaB{So} the map of
the lemma is onto $U_Q\cap G_{\RamiC{\fx,}0}$\RamiC{. I}t is injective as a restriction of the bijection above. 
\item \RamiC{$
\prod_{\beta\in\Phi_Q}U_{\beta,0}\xrightarrow{\ \sim\ }U_Q\cap G_{0}
$
is a bijection}.\\
\RamiC{Let $G^0<G$ be as in \cite[Definition 2.6.23]{KP}. 
It is contained in  $^0G$ by \cite[Notation 2.6.15, Lemma 11.5.2, Lemma 2.6.16(4)]{KP},  
\cite[V.2.3]{Ren}. 
By
\cite[Proposition 7.6.4(1)]{KP} we have 
$G_{0}\cap G^0=G_{\fx,0}$.
So, by Step 3, it is enough to show that
$U_Q\subset G^0$.

Let $\beta\in\Phi_Q$ and $u\in U_\beta\smallsetminus\{1\}$, and put $r:=\phi_\beta(u)$,
which is finite by axiom V\,1 of \cite[Definition 6.1.2]{KP}; thus $u\in U_{\beta,r}$. Since
$\beta$ is a non-zero character of $\bfA_\bfT$, and the valuation identifies $\Lambda(A_\bfT)$
with $X_*(\bfA_\bfT)$ by \cite[V.2.6]{Ren}, there is $s\in A_\bfT$ with $\beta(s)<0$; choose
$n\in\Z_{>0}$ with $n\beta(s)\leq r$. The filtration $(U_{\beta,r})_{r\in\R}$ is decreasing
by \Cref{not:root.sub}\eqref{not:root.filtration}, so
$$
u\in U_{\beta,r}\subset U_{\beta,n\beta(s)}
\overset{\text{Lem \ref{lem:root.conj}}}{=}s^nU_{\beta,0}s^{-n}\subset G^0 ,
$$
the last inclusion because $U_{\beta,0}<G^\sharp_{\fx,0}<G_{\fx,0}<G^0$ and $G^0$ is normal
in $G$. Thus $U_\beta\subset G^0$ for every $\beta\in\Phi_Q$, and therefore $U_Q\subset G^0$
by \Cref{lem:U.filtration.alg}.
}
\end{enumerate}
\end{proof}

\begin{lemma}\label{lem:U.filtration}
There is $\mdef{c_U}\in\Z_{>0}$ with the following property. Let \RamiC{$\bfQ<\bfG$} be a
standard parabolic or the opposite of one, let $\RamiC{\bfU_\bfQ}$ be its unipotent radical. \RamiC{Let $Q:=\bfQ(F)$ and $U_Q:=\bfU_\bfQ(F)$. Fix}
an ordering of $\Phi_Q$. Then, with the notation of \Cref{not:root.sub}:
\begin{enumerate}
\item\label{lem:U.filtration:1} 
\RamiF{For every $n\in\Z_{\geq 0}$ the multiplication map
$$
\prod_{\beta\in\Phi_Q}U_{\beta,n}\longrightarrow U_Q
$$
is injective, and its image satisfies
$$
U_Q\cap G_{-2n}\subset\prod_{\beta\in\Phi_Q}U_{\beta,n}\subset U_Q\cap G_{-n}.
$$
}
\item\label{lem:U.filtration:2} For every $m\in\Z_{>0}$,
$$
U_Q\cap G_m\subset\prod_{\beta\in\Phi_Q}U_{\beta,-c_Um} ,
$$
where the product is taken in the fixed order.
\item\label{lem:U.filtration:3}
\RamiF{Let $s\in T$ satisfy $\beta(s)\geq0$ for every $\beta\in\Phi_Q$. Then for every
$n\in\Z_{\geq0}$
$$
s\,(U_Q\cap G_{-n})\,s^{-1}\subset U_Q\cap G_{-n} .
$$}
\end{enumerate}
\end{lemma}
\begin{proof}
$ $
\begin{itemize}
    \item [\eqref{lem:U.filtration:1}]
 For $n=0$ the assertion is \Cref{lem:U.level0}; assume from now on that $n\geq1$. By \S\ref{ssec:not}, $G_{-n}=G_{\fx,n}$ is the Moy--Prasad subgroup
of \cite[Definition 13.2.1]{KP}. 
\RamiF{Without loss of generality assume that $Q$ is standard.}
Applying
\cite[Proposition 13.2.5(3)]{KP} with $r=n$, 
\RamiF{we obtain a bijection:
$$\left(\prod_{\alpha\in \Phi_{\bfB^-}}U_{\alpha,n}U_{2\alpha,n} \right)\times T_{-n} \times \left(\prod_{\alpha\in\Phi_{\bfB}}U_{\alpha,n}U_{2\alpha,n} \right)\to G_{-n}$$
where the order can be chosen arbitrary.
Chose a total ordering on $\Phi_B$ such that the elements of $\Phi_Q$ are bigger than the rest (so they are at the end). We obtain a bijection

$$\left(\prod_{i}\Gamma_i\right)\times 
\left(\prod_{\beta\in\Phi_Q}U_{\beta,n}U_{2\beta,n} \right)\to G_{-n},$$
where $\Gamma_i<Q^-$ are some subgroups.

Since $U_Q\cap Q^-=1$ we obtain a bijection
}

\RamiF{
}
\begin{equation}\label{eq:U.filtration.exact}
\RamiF{\prod_{\beta\in\Phi_Q}U_{\beta,n}U_{2\beta,n}
\to
U_Q\cap G_{-n},}
\end{equation}

\RamiF{
As
$U_{\beta,n}\subset U_{\beta,n}U_{2\beta,n}$, the map of the assertion is a restriction of that
bijection. It is therefore injective, and its image is contained in $U_Q\cap G_{-n}$. For the
remaining inclusion, \Cref{not:root.sub}\eqref{not:root.div} gives
$U_{2\beta,2n}=U_{2\beta}\cap U_{\beta,n}\subset U_{\beta,n}$ and hence
$U_{\beta,2n}U_{2\beta,2n}\subset U_{\beta,n}$ for every $\beta\in\Phi_Q$. Multiplying these
inclusions in the fixed order and applying \eqref{eq:U.filtration.exact} with $2n$ in place of $n$
gives $U_Q\cap G_{-2n}\subset\prod_{\beta\in\Phi_Q}U_{\beta,n}$.}
\item[\eqref{lem:U.filtration:3}]
\RamiF{
 The bijection \eqref{eq:U.filtration.exact} holds for $n=0$ as well: its
left hand side is then $\prod_{\beta\in\Phi_Q}U_{\beta,0}$, because
$U_{2\beta,0}=U_{2\beta}\cap U_{\beta,0}$ by \Cref{not:root.sub}\eqref{not:root.div}, and the
identity itself is \Cref{lem:U.level0}. It is therefore enough to show that the conjugate by $s$ of
each factor of \eqref{eq:U.filtration.exact} is contained in that factor, i.e. that
$sU_{\beta,n}s^{-1}\subset U_{\beta,n}$ and $sU_{2\beta,n}s^{-1}\subset U_{2\beta,n}$ for every
$\beta\in\Phi_Q$. By \Cref{lem:root.conj} these two conjugates are $U_{\beta,\,n+\beta(s)}$ and
$U_{2\beta,\,n+2\beta(s)}$, so both inclusions follow from $\beta(s)\geq0$ and the monotonicity of
the two filtrations.}

\item[\eqref{lem:U.filtration:2}] 
There are only finitely many parabolic subgroups containing $T$, so it
is enough to fix $Q$ and to produce a constant $c_Q$; taking the maximum over the finitely many
$Q$ then gives $c_U$. If $\Phi_Q=\varnothing$ then $U_Q=1$ and there is nothing to prove, so
assume $\Phi_Q\neq\varnothing$; then $Q$ is proper and, $\bfA_{\bfM_Q}$ being central in
$\bfM_Q$, there is $\lambda_0\in\Lambda(A_{\bfM_Q})$ with $\beta(\lambda_0)>0$ for every
$\beta\in\Phi_Q$. Put $a:=\varpi^{\lambda_0}\in T$ and fix $D\in\Z_{>0}$ with $\beta(a)\leq D$
for every $\beta\in\Phi_Q$.

\emph{The grading defined by $\lambda_0$.} As $\bfA_\bfT$ is split, $\cW$ is the direct sum of its
$\bfA_\bfT$-weight spaces; let $\cW_j$ be the sum of those whose weight takes the value $j\in\Z$ at
$\lambda_0$, so that $\sigma(a)$ acts on $\cW_j$ by $\varpi^{j}$. Since $\beta(\lambda_0)>0$ for
every $\beta\in\Phi_Q$, the group $\bfU_Q$ is contained in the unipotent radical of the parabolic
subgroup of $\bfG$ attached to $\lambda_0$; hence $\sigma(U_Q)$ preserves the flag
$\cW_{\geq j}:=\bigoplus_{j'\geq j}\cW_{j'}$ and acts trivially on
$\bigoplus_j\cW_{\geq j}/\cW_{\geq j+1}$. In other words, for $v\in U_Q$ the component of
$\sigma(v)-1$ from $\cW_j$ to $\cW_{j'}$ vanishes unless $j'\geq j+1$. Fix an $\cO_F$-basis of $\cW$
adapted to the decomposition $\cW=\bigoplus_j\cW_j$, let $\cL'$ be its $\cO_F$-span, and let
$c_2\in\Z_{>0}$ satisfy
$\varpi^{c_2}\End_{\cO_F}(\cL)\subset\End_{\cO_F}(\cL')\subset\varpi^{-c_2}\End_{\cO_F}(\cL)$; such a
$c_2$ exists because $\End_{\cO_F}(\cL)$ and $\End_{\cO_F}(\cL')$ are $\cO_F$-lattices in
$\End_F(\cW)$.

\emph{Conjugating by a power of $a$.} Let $m\in\Z_{>0}$ and $u\in U_Q\cap G_m$. By
\Cref{cor:ov.sigma.G} there is $c_1\in\Z_{>0}$, independent of $m$ and $u$, with
$ov_\sigma(g)\leq c_1(m+1)$ for every $g\in G_m$; hence
$$
\sigma(u)-1\in\varpi^{-c_1(m+1)}\End_{\cO_F}(\cL)\subset\varpi^{-c_1(m+1)-c_2}\End_{\cO_F}(\cL') .
$$
Conjugation by $\sigma(a)^k$ multiplies the component from $\cW_j$ to $\cW_{j'}$ by
$\varpi^{\,k(j'-j)}$, and $k(j'-j)\geq k$ for every non-vanishing component of $\sigma(u)-1$, so
for every $k\in\Z_{>0}$
\begin{align*}
\sigma(a^kua^{-k})-1&=\sigma(a)^k\bigl(\sigma(u)-1\bigr)\sigma(a)^{-k}\\
&\in\varpi^{\,k-c_1(m+1)-c_2}\End_{\cO_F}(\cL')
\subset\varpi^{\,k-c_1(m+1)-2c_2}\End_{\cO_F}(\cL).
\end{align*}

\emph{Conclusion.} Take $k:=c_{nbh}+c_1(m+1)+2c_2$, with $c_{nbh}$ the constant of
\Cref{lem:sigma.nbhd}. Then $\sigma(a^kua^{-k})\in1+\varpi^{c_{nbh}}\End_{\cO_F}(\cL)$, so
$a^kua^{-k}\in G_{-1}$ by \Cref{lem:sigma.nbhd}. As $a\in T$ normalizes $U_Q$, this gives
$$
\RamiF{a^kua^{-k}\in U_Q\cap G_{-1}\subset U_Q\cap G_{0}
\overset{\text{Lem \ref{lem:U.level0}}}{=}\prod_{\beta\in\Phi_Q}U_{\beta,0}},
$$
the product being taken in the fixed order. Conjugating back by $a^{-k}$ and using
\Cref{lem:root.conj},
$$
\RamiF{u\in\prod_{\beta\in\Phi_Q}U_{\beta,\,-k\beta(a)}\subset\prod_{\beta\in\Phi_Q}U_{\beta,-kD}},
$$
the last inclusion because the filtration is decreasing and $\beta(a)\leq D$. Finally
$k\leq(c_{nbh}+2c_1+2c_2)m$ because $m\geq1$, so
$c_Q:=D(c_{nbh}+2c_1+2c_2)$ satisfies the assertion.
\end{itemize}
\end{proof}

\begin{proof}[Proof of \Cref{lem:U.contract}]
$ $

\eqref{lem:U.contract:1} By \Cref{lem:U.filtration}\eqref{lem:U.filtration:2} we have
$U_Q\cap G_k\subset\prod_{\beta\in\Phi_Q}U_{\beta,-c_Uk}$, so conjugation by $a$ carries
it into $\prod_{\beta\in\Phi_Q}U_{\beta,-c_Uk+\beta(a)}$, and
$$
-c_Uk+\beta(a)\geq -c_Uk+c_{cnj}k\geq k
$$
provided $c_{cnj}\geq 1+c_U$. The image therefore
lies in $\prod_{\beta\in\Phi_Q}U_{\beta,k}$, which \RamiF{is contained in}
 $U_Q\cap G_{-k}$ by
\Cref{lem:U.filtration}\eqref{lem:U.filtration:1}.

\eqref{lem:U.contract:2} 
\RamiC{
\begin{enumerate}[{Case} 1.]
    \item $g\in T$.\\
    \RamiF{By \Cref{lem:U.filtration}\eqref{lem:U.filtration:3}, applied to $s:=ga\in T$, it is
enough to prove that $\beta(a)+\beta(g)\geq0$ for every $\beta\in\Phi_Q$.}

The elements of $\Phi_Q$ are finitely many $\Q$-valued functionals
on $\Lambda(T)$.
So \Cref{lem:ov.G.comp} provides
$c_1\in\Z_{>0}$,  with
$$
|\beta(t)|\leq c_1(ov_G(t)+1),
$$
for any $t\in T$ and $\beta\in \Phi_Q$.
Take $c_{cnj}=\max(1+c_U,2c_1)\geq 2c_1$.
The assertion follows from $\beta(a)\geq-\beta(g)$,
which  foll\EitanE{ows} from \EitanE{the} previous inequality and the \EitanE{assumption.}

Taking $c_{cnj}:=\max(1+c_U,2c_1)$ proves both assertions.
    \item The general case.\\ 
    Take $c_{cnj}$ be as in the previous case.
Let $Q,k,a,g, n$ be as \Cref{lem:U.contract}\eqref{lem:U.contract:2}.    
    
    By the Cartan decomposition in $M_Q$ write $g=k_1tk_2$ with $k_1,k_2\in M_0$ and $t\in T$. As
$M_0<G_0$, the elements $g$ and $t$ lie in the same $G_0\times G_0$ double coset, so
$ov_G(t)=ov_G(g)\leq k$ by \Cref{lem:ov.cartan}. The group $U_Q\cap G_{-n}$ is normalized by $M_0$, since $M_0<G_0$ normalizes
$G_{-n}$ and $M_Q$ normalizes $U_Q$; and $a$ is central in $M_Q$, so $a(U_Q\cap G_{-n})a^{-1}$ is
normalized by $M_0$ as well. The as\EitanE{s}ertion follows now from the \EitanE{previous} case:
\begin{align*}
a\,(U_Q\cap G_{-n})\,a^{-1}
&=
k_2^{-1} a\,(U_Q\cap G_{-n})\,a^{-1} k_2 \\
&<
k_2^{-1} t^{-1} (U_Q\cap G_{-n}) t k_2 \\
&= k_2^{-1} t^{-1} k_1^{-1}(U_Q\cap G_{-n})k_1 t k_2 \\
&=
g^{-1}(U_Q\cap G_{-n})\,g .
\end{align*}
\end{enumerate}
}
\end{proof}

\subsection{Cosets of the groups $G_{-n}$ inside balls}\label{ssec:app.ball}
In this subsection we prove \Cref{lem:coset.ball} that relates balls of the rectified variety $\bfG^{rss}$ to cosets of the groups $G_{-n}$.

Its proof goes through the embedding $\sigma$ of \Cref{not:sigma}:
\Cref{lem:cong} bounds the matrix entries of $\sigma$ on $G_{-n}$, and the coordinates of a
rectification of $\bfG^{rss}$ are polynomials in these matrix entries and in $\fD^{-1}$.

\begin{lemma}\label{lem:cong}
There is $\mdef{c_{cng}}\in\Z_{\geq0}$ such that for every $n\in\Z_{>0}$
$$
\sigma(G_{-n})\subset 1+\varpi^{\,n-c_{cng}}\End_{\cO_F}(\cL).
$$
\end{lemma}

\begin{proof}
Let $\RamiC{F^{un}}$ be the maximal unramified extension of $F$ and $\cO_{\RamiC{F^{un}}}$ its \RamiC{ring of integers.}

Let $\cG$ be the $\cO_F$-group scheme attached to the point $\fx$ and to the constant concave
function $0$, that is, the descent to $\cO_F$ --- \cite[Proposition 9.8.3]{KP} --- of the
$\cO_{\RamiC{F^{un}}}$-group scheme of \cite[Theorem 8.5.2]{KP}.
\RamiC{$\cG$ satisfies the following:
\begin{itemize}
    \item It is affine, smooth and of finite type.
    \item Its
generic fiber is $\bfG$.
\item $\cG(\cO_{\RamiC{F^{un}}})=G({\RamiC{F^{un}}})_{\fx,0}$, see 
\cite[Theorem 8.5.2]{KP}.
\item $\cG(\cO_F)=\RamiC{G_{\fx,0}\subset} G_0$, see 
\cite[Proposition 9.8.3(2)]{KP}.
\end{itemize}
}
By \cite[Proposition 8.5.16]{KP}, applied over $\cO_{\RamiC{F^{un}}}$ to the
constant concave function $0$ and to the integer $n$, the group scheme attached to the constant
concave function $n$ is the $n$-th congruence group scheme of $\cG_{\cO_{\RamiC{F^{un}}}}$, so by
\cite[Definition A.5.12(1) and Lemma A.5.13(3)]{KP}
$$
G({\RamiC{F^{un}}})_{\fx,n}=\ker\bigl(\cG(\cO_{\RamiC{F^{un}}})\to \cG(\cO_{\RamiC{F^{un}}}/\varpi^n\cO_{\RamiC{F^{un}}})\bigr).
$$
By \cite[Proposition 9.8.3(2)]{KP} again, $G_{-n}=G_{\fx,n}=G({\RamiC{F^{un}}})_{\fx,n}\cap G$. Thus for
$g\in G_{-n}\RamiC{\subset \cG(\cO_F)}$ the two $\cO_{\RamiC{F^{un}}}$-algebra homomorphisms $\cO_{\RamiC{F^{un}}}[\cG_{\cO_{\RamiC{F^{un}}}}]\to\cO_{\RamiC{F^{un}}}$ corresponding
to $g$ and to $1$ agree modulo $\varpi^n$.

Let $\sigma_{ij}\in F[\bfG]$ be the matrix entries of $\sigma$ in the fixed $\cO_F$-basis of
$\cL$. Since the generic fiber of $\cG$ is $\bfG$ we have $F[\bfG]=\cO_F[\cG]\otimes_{\cO_F}F$,
so, there is $c_{cng}\in\Z_{\geq0}$ with
$\varpi^{c_{cng}}\sigma_{ij}\in\cO_F[\cG]$ for all $i,j$. Applying the previous paragraph to the
elements $\varpi^{c_{cng}}\sigma_{ij}$ of $\cO_F[\cG]\subset\cO_{\RamiC{F^{un}}}[\cG_{\cO_{\RamiC{F^{un}}}}]$ we get, for
every $g\in G_{-n}$,
$$
\varpi^{c_{cng}}\bigl(\sigma(g)_{ij}-\delta_{ij}\bigr)\in\varpi^{n}\cO_{\RamiC{F^{un}}}\cap F=\varpi^n\cO_F ,
$$
since $\sigma_{ij}(1)=\delta_{ij}$. Hence
$\sigma(g)_{ij}-\delta_{ij}\in\varpi^{\,n-c_{cng}}\cO_F$, which is the assertion.
\end{proof}
\RamiC{
\begin{lemma}\label{lem:perturb}
     For every
$P\in F[\bfG]$ there is $c_P\in\Z_{>0}$ such that for every 
$s\in \Z_{>0}$ and 
$g,u\in G$ with 
$\sigma(u)\in1+\varpi^s\End_{\cO_F}(\cL)$
we have:
\begin{equation}
|P(g)|\leq \ell^{c_P(1+ov_G(g))}\qquad\text{and}\qquad
|P(gu)-P(g)|\leq \ell^{c_P(1+ov_G(g))-s}.
\end{equation}
\end{lemma}
}
\RamiC{
\begin{proof}
Since $1+\varpi^s\End_{\cO_F}(\cL)$ is a subgroup of $\GL(\cW)$, as noted in \Cref{not:N1}, we
have $\sigma(u)^{-1}\in1+\varpi^s\End_{\cO_F}(\cL)$ as well. Thus $\sigma(u)$ and $\sigma(u)^{-1}$
preserve $\cL$, so $ov_\sigma(u)=0$, and $N_1(u)\leq\ell^{-s}$ by \eqref{eq:N1.cong}.

We first show that there are $a,b\in\Z_{>0}$, depending only on $P$, with
$$
|P(g)|\leq \ell^{a+b\,ov_\sigma(g)}\qquad\text{and}\qquad
|P(gu)-P(g)|\leq \ell^{a+b\,ov_\sigma(g)-s}
$$
for every $s$, $g$ and $u$ as in the statement.

The first inequality is \cite[(18.2.2)]{Kot} for the norm $\ell^{ov_\sigma}$ of
\Cref{lem:norm.sigma}. For the second, apply \cite[Lemma 19.2]{Kot} to the morphism
$\bfG\times\bfG\to\bA^1$, $(g,u)\mapsto P(gu)-P(g)$, which carries the closed subvariety
$\bfG\times\{1\}$ to $\{0\}$: it gives a norm $\|\cdot\|$ on $G\times G$ with
$$
|P(gu)-P(g)|\leq \|(g,u)\|\cdot N_{\bfG\times\{1\}}(g,u).
$$
Here $N_{\bfG\times\{1\}}(g,u)=N_1(u)$, as in \cite[(19.3.3)]{Kot}: the generators fixed in
\Cref{not:N1} pull back, along the second projection, to generators of the ideal of
$\bfG\times\{1\}$ in $F[\bfG\times\bfG]$. Moreover $\|\cdot\|$ is equivalent to
$\max(\ell^{ov_\sigma(g)},\ell^{ov_\sigma(u)})$ by \cite[Proposition 18.1(4)]{Kot} and
\Cref{lem:norm.sigma}, and $ov_\sigma(u)=0$, so $\|(g,u)\|\leq\ell^{a+b\,ov_\sigma(g)}$ for
suitable $a,b$; together with $N_1(u)\leq\ell^{-s}$ this gives the second inequality.

By \Cref{cor:ov.sigma.G} there is $c\in\Z_{>0}$ with $ov_\sigma(g)\leq c\bigl(ov_G(g)+c\bigr)$ for
every $g\in G$. Hence, for $c_P:=a+bc^2+bc$,
$$
a+b\,ov_\sigma(g)\leq a+bc^2+bc\,ov_G(g)\leq c_P\bigl(1+ov_G(g)\bigr) ,
$$
and the two inequalities of the statement follow.
\end{proof}
}
\RamiC{
\begin{proof}[Proof of \Cref{lem:coset.ball}]
Let $i=(i_1,\dots,i_M):\bfG^{rss}\to\bA^M$ be the closed
embedding of the rectification.

As $\bfG^{rss}$ is the principal open subset of $\bfG$ defined by $\fD$,  there are $N\in\Z_{\geq0}$ and $P_1,\dots,P_M\in F[\bfG]$ with
$i_j=P_j\fD^{-N}$. 

By \Cref{lem:perturb} there is $c_2\in \Z_{>0}$ such that 
for every 
$s\in \Z_{>0}$, $1\leq j\leq M$, and 
$g,u\in G$ with 
$\sigma(u)\in 1+\varpi^s\End_{\cO_F}(\cL)$ we have: 
\begin{itemize}
\item $|P_j(g)|\leq\ell^{c_2(ov_G(g)+1)}$ 
\item $|P_j(gu)-P_j(g)|\leq\ell^{c_2(ov_G(g)+1)-s}$  
\item $|\fD(g)|\leq\ell^{c_2(ov_G(g)+1)}$
\item $|\fD(gu)-\fD(g)|\leq\ell^{c_2(ov_G(g)+1)-s}$
\end{itemize}
The proof will proceed in the following steps.
\begin{enumerate}[Step 1.]
    \item For any 
    $g\in G^{rss}$, $s>(c_2+1)(ov_{G^{rss}}(g)+1)$ and    
    $u\in G$ with
$\sigma(u)\in 1+\varpi^s\End_{\cO_F}(\cL)$ we have 
$gu\in G^{rss}$:\\
Put $R:=ov_{G^{rss}}(g)$. By \Cref{not:settings}\eqref{it:Gn} and \Cref{not:rss} we have
$ov_G(g)\leq R$ and $\operatorname{val}(\fD(g))\leq R$, that is $|\fD(g)|\geq\ell^{-R}$. Hence, by
the fourth bound above and $s>(c_2+1)(R+1)=c_2(R+1)+R+1$, \EitanE{we obtain:}
$$
|\fD(gu)-\fD(g)|\leq\ell^{c_2(R+1)-s}<\ell^{-R}\leq|\fD(g)| .
$$
Therefore $|\fD(gu)|=|\fD(g)|\neq0$, and in particular $gu\in G^{rss}$.
\item There is $c_3\in\Z_{>0}$ s.t. for every $g\in G^{rss}$, $m\in \Z_{>0}$, $s>c_3(ov_{G^{rss}}(g)+1)+m$ and $u$ with $\sigma(u)\in 1+\varpi^s\End_{\cO_F}(\cL)$ we have   
$gu\in B_{-m}^{\bfG^{rss}}(g)$:\\
Put $c_3:=(N+2)c_2+2N+2$, let $g$, $m$, $s$ and $u$ be as in the statement, and write
$R:=ov_{G^{rss}}(g)$. As $c_3\geq c_2+1$ and $m>0$ we have
$s>(c_2+1)(R+1)$, so Step 1 gives $|\fD(gu)|=|\fD(g)|$; moreover $ov_G(g)\leq R$,
$\ell^{-R}\leq|\fD(g)|\leq\ell^{c_2(R+1)}$ and $|\fD(g)-\fD(gu)|\leq\ell^{c_2(R+1)-s}$, as in Step 1.

The rectification being simple,
$B^{\bfG^{rss}}_{-m}(g)=i^{-1}\bigl(i(g)+\varpi^{m}\cO_F^{M}\bigr)$ by
\Cref{def:simp balls}\eqref{def:simp balls:2}, so it is enough to show that
$|i_j(gu)-i_j(g)|\leq\ell^{-m}$ for $1\leq j\leq M$.
From $|\fD(g)^{N}-\fD(gu)^{N}|\leq|\fD(g)|^{N-1}|\fD(g)-\fD(gu)|$, $|\fD(g)|^{N}\leq\ell^{Nc_2(R+1)}$ and
$|\fD(g)|^{-1}\leq\ell^{R}$ we get
\begin{align*}
\bigl|\fD(gu)^{-N}-\fD(g)^{-N}\bigr|=\frac{\bigl|\fD(g)^{N}-\fD(gu)^{N}\bigr|}{|\fD(g)\fD(gu)|^{N}}
&\leq\ell^{Nc_2(R+1)+R}\cdot\ell^{2NR}\cdot\ell^{c_2(R+1)-s}\\
&=\ell^{(N+1)c_2(R+1)+(2N+1)R-s} .
\end{align*}
Since $i_j=P_j\fD^{-N}$ we have
$$
i_j(gu)-i_j(g)=\bigl(P_j(gu)-P_j(g)\bigr)\fD(gu)^{-N}+P_j(g)\bigl(\fD(gu)^{-N}-\fD(g)^{-N}\bigr) ,
$$
and $|\fD(gu)^{-N}|\leq\ell^{NR}$, so by the first two bounds above both summands are at most
$\ell^{(N+2)c_2(R+1)+(2N+1)R-s}\leq\ell^{c_3(R+1)-s}$. As $s>c_3(R+1)+m$, this is smaller than
$\ell^{-m}$, as required.
\item Completion of the proof:\\
Let $x\in G^{rss}$ and $m\in\Z_{>0}$, put $R:=ov_{G^{rss}}(x)\in\Z_{\geq0}$ and
$$
c_{ball}:=2(c_{cng}+c_3+1) ,
$$
and let $n:=c_{ball}(m+R)$ and $u\in G_{-n}$. By \Cref{lem:cong} we have
$\sigma(u)\in1+\varpi^{s}\End_{\cO_F}(\cL)$ for $s:=n-c_{cng}$. Since $m+R\geq1$ we have
$2c_{cng}(m+R)-c_{cng}\geq0$ and $2(m+R)\geq m+R+1$, so
$$
s=2c_{cng}(m+R)-c_{cng}+2(c_3+1)(m+R)\geq(c_3+1)(m+R+1)>c_3(R+1)+m .
$$
In particular $s\in\Z_{>0}$, so Step 2 applies and gives $xu\in B^{\bfG^{rss}}_{-m}(x)$. As $u$ was
an arbitrary element of $G_{-c_{ball}(m+R)}$, this is the assertion of \Cref{lem:coset.ball}.
\end{enumerate}
\end{proof}
}

\subsection{\RamiA{\RamiC{Compatibility}  of the scale with Levi subgroups}}\label{ssec:cGM}
\begin{proof}[Proof of \Cref{lem:cGM}]
$ $
\begin{enumerate}
    \item
$\bfM$ is a closed subgroup of $\bfG$, so the restriction of $\|\cdot\|_G$ to $M$ is a norm on
$M$ by \Cref{lem:nrm.cl.emb}, and $m\mapsto\ell^{ov_M(m)}$ is a norm on $M$ by
\Cref{rem:app.Levi}. Two norms on the same set being equivalent, there is a constant $c_1$ with
$$
ov_G(m)\leq c_1(ov_M(m)+1)\qquad (m\in M).
$$

Choose $c^G_M\RamiC{=2c_1}\in \Z_{>0}$ \RamiC{so}
$$
c_1(n+1)\leq c^G_M n
$$
for every $n\in \Z_{>0}$. If $m\in M_n$, then $ov_M(m)\leq n$, hence
$$
ov_G(m)\leq c_1(n+1)\leq c^G_M n.
$$
Thus $m\in G_{c^G_M n}$, proving the first assertion.
\item By \cite[Proposition 13.2.5(6)]{KP}, applied at the point $\beta_M(\fx)$, we have
$G_{-n}\cap M=M_{-n}$. 

\RamiA{Let $P$ be the standard parabolic whose Levi is $M$, let $U$ be its unipotent radical and let $U^{-}$ be the unipotent radical of its opposite parabolic.}
By \Cref{lem:iva.decomp} every $g\in G_{-n}$ \RamiC{can be written as} $u^-mu$ with
$u^-\in U^-\cap G_{-n}$, $m\in M\cap G_{-n}$ and $u\in U\cap G_{-n}$; if $g\in P$ then
$u^-=g(mu)^{-1}\in U^-\cap P=1$, so
$$
G_{-n}\cap P=(M\cap G_{-n})(U\cap G_{-n}).
$$
Applying $pr_M$, which is the identity on $M$ and trivial on $U$, gives
$(G_{-n})_M\RamiC{=pr_M(G_{-n} \cap P)}=M\cap G_{-n}=M_{-n}$.
\end{enumerate}
\end{proof}

\RamiA{
\subsection{Effective Iwasawa decomposition}\label{ssec:eff.Iwa}
\begin{proof}[Proof of \Cref{lem.Iv}]
\RamiC{Set $U:=U_\lambda$ and $M:=M_\lambda$}.
Let $i_U:U_\lambda\to G, \,i_M:M_\lambda\to G$ be the standard embeddings. Let  $i_{U\times M}:U\times M\to G$ be defined by $i(u,m)=um$. This is also a closed embedding. Hence, by \Cref{lem:nrm.cl.emb},  \Cref{lem:nrm.prd} and \Cref{prop:G.norms}\eqref{prop:G.norms:1}, the abstract norms 
$i_U^*(\|\cdot\|_G)\times i_M^*(\|\cdot\|_G)$ and
$i_{U\times M}^*(\|\cdot\|_G)$ are norms on $U\times M$.
Thus we have a  constant $c\in \Z_{>1}$ s.t.  for any $(u,m)\in U\times M$ we have
\begin{equation}\label{eq:Iwa.NDP}
\RamiC{\max(ov_G(u), ov_G(m))}\leq c(ov_G(um)+1).    
\end{equation}

Take $c_I:=\RamiC{2c}$. Let $n\in \Z_{>1}$ and let $g\in G_n$.
By the Iwasawa decomposition $G=U_\lambda M_\lambda G_0$ (see e.g. \cite[V.5.1]{Ren}). 
Write $g=umk$ with $u\in U_\lambda$, $m\in M_\lambda$, $k\in G_0$.
It is enough to show that $u,m\in G_{nc_I}$.
We have:
\begin{align*}
\RamiC{\max(ov_G(u), ov_G(m))}
&\overset{\eqref{eq:Iwa.NDP}}{\leq} c(ov_G(um)+1)=c(ov_G(gk^{-1})+1)
\overset{\text{Lem~\ref{lem:ov.cartan}}}{\DimaB{=}}\\
&\DimaB{=} c(ov_G(g)+1)\leq {nc_I}.
\end{align*}
as required.
\end{proof}
}
\begingroup
 \let\clearpage\relax
  \let\cleardoublepage\relax  
  \printindex
\endgroup

\bibliographystyle{alpha}
\bibliography{Ramibib}
\end{document}